\input epsf
\documentstyle{amsppt}
\magnification=\magstep1
\pagewidth{6.48truein}
\pageheight{9.0truein}
\TagsOnRight
\NoRunningHeads
\catcode`\@=11
\def\logo@{}
\footline={\ifnum\pageno>1 \hfil\folio\hfil\else\hfil\fi}
\topmatter
\title Plane partitions I: a generalization of MacMahon's formula
\endtitle
\author Mihai Ciucu\endauthor
\subjclass
Primary 05A15, 05B45, 05A17. 
Secondary 52C20, 11P81
\endsubjclass
\keywords
Plane partitions, lozenge tilings, non-intersecting
lattice paths, tiling enumeration, perfect matchings
\endkeywords
\thanks 
The author was supported by a Membership at the Institute for Advanced Study.
\endthanks
\affil
  Georgia Institute of Technology\\
  School of Mathematics\\
  Atlanta, GA 30332-0160
\endaffil
\abstract
The number of plane partitions contained in a given box was shown by MacMahon
\cite{8} to be given by a simple product formula. By a simple bijection, this
formula also enumerates lozenge tilings of hexagons of side-lengths $a,b,c,a,b,c$ (in
cyclic order) and angles of 120 degrees.
We present a generalization in the case $b=c$ by giving simple product formulas
enumerating lozenge tilings of regions obtained from a hexagon of 
side-lengths $a,b+k,b,a+k,b,b+k$ (where $k$ is an arbitrary non-negative integer) and
angles of 120 degrees by removing certain triangular regions along its symmetry axis.  
\endabstract
\endtopmatter
\document

\def\mysec#1{\bigskip\centerline{\bf #1}\message{ * }\nopagebreak\par}

\def\myref#1{\item"{[{\bf #1}]}"} 
 
\def\pf{{\it Proof.\ }} 

\def\cite#1{\relaxnext@
  \def\nextiii@##1,##2\end@{[{\bf##1},\,##2]}%
  \in@,{#1}\ifin@\def\next{\nextiii@#1\end@}\else
  \def\next{[{\bf#1}]}\fi\next}
\def\proclaimheadfont@{\smc}

\def\pf{{\it Proof.\ }}
\define\N{{\bold N}}
\define\Z{{\bold Z}}

\define\w{\operatorname{w}}
\define\M{\operatorname{M}}
\define\twoline#1#2{\line{\hfill{\smc #1}\hfill{\smc #2}\hfill}}

\def\mypic#1{\epsffile{#1}}

\mysec{1. Introduction}
 
\medskip
A plane partition is a rectangular array of nonnegative integers with the property
that all rows and columns are weakly decreasing. A plane partition
$\pi=(\pi_{ij})_{0\leq i<a,0\leq j<b}$ can be identified with its three dimensional 
diagram
$D_\pi=\{(i,j,k):0\leq k<\pi_{ij}\}$, and hence can be viewed as an order ideal of $\N^3$
(an order ideal of a partially ordered set is a subset $I$ such that $x\in I$ and
$y\leq x$ imply $y\in I$). 

Replacing each point of $D_\pi$ by a unit cube
centered at it and parallel to the coordinate axes, the diagram of $\pi$ becomes a
stack of unit cubes justified onto the three coordinate planes. By projection on a
plane normal to the vector $(1,1,1)$ one obtains that plane partitions whose diagrams
fit inside an $a\times b\times c$ box are identified with tilings of a hexagon of
side-lengths $a,b,c,a,b,c$ (in cyclic order) and angles of 120 degrees by unit rhombi
with angles of 60 and 120 degrees (see \cite{3} or \cite{7} for the details of this 
bijection); we call such rhombi {\it lozenges} and such tilings
{\it lozenge tilings}. 

MacMahon showed \cite{8} that the number of plane partitions contained in an 
$a\times b\times c$ box (equivalently, the number of lozenge tilings of the
corresponding hexagon), with $a\leq b$, equals

$$\frac{(c+1)(c+2)^2\cdots(c+a)^a(c+a+1)^a\cdots(c+b)^a(c+b+1)^{a-1}\cdots(c+b+a-1)}
{1\cdot2^2\cdots a^a\cdot(a+1)^a\cdots b^a\cdot(b+1)^{a-1}\cdots(b+a-1)}.$$

The elegance of this result strongly invites one to find generalizations of it. Two natural 
directions in which attempts have been made are the following. First, one may consider
enumerating order ideals of the product of $d$ finite chains, with $d\geq4$ (the case
$d=4$ was already considered by MacMahon \cite{9}). Numerical evidence suggests that
there is no simple product formula for the answer to this problem, and very little is known
about it (see \cite{11,I \S{11}}). Second, taking the viewpoint of tilings, one
may consider the question of counting tilings by unit rhombi of $2d$-gons with parallel
and congruent opposite sides. Again, there seem to be very few results on this in the
literature for $d\geq4$ (see \cite{4}).

We generalize MacMahon's result in the case $b=c$ by presenting a general family of
hexagonal regions with holes on the triangular lattice, such that the number of
tilings of each member of this family is given by a simple product formula (examples of
regions from this family are shown in Figures 1.2--1.4).

Consider the tiling of the plane by unit
equilateral triangles. Define a {\it region} to be any region of the plane that can be
obtained as a finite union of these unit triangles. 

Consider the hexagonal region $H(a,b,k)$ with side-lengths $a,b+k,b,a+k,b,b+k$, where $k$ is a 
nonnegative integer, and suppose $H(a,b,k)$ is drawn such that its base is the edge of length
$a+k$ (see Figure 1.1 for an example). Then $H(a,b,k)$ has a vertical symmetry axis $\ell$.

\topinsert
\centerline{\mypic{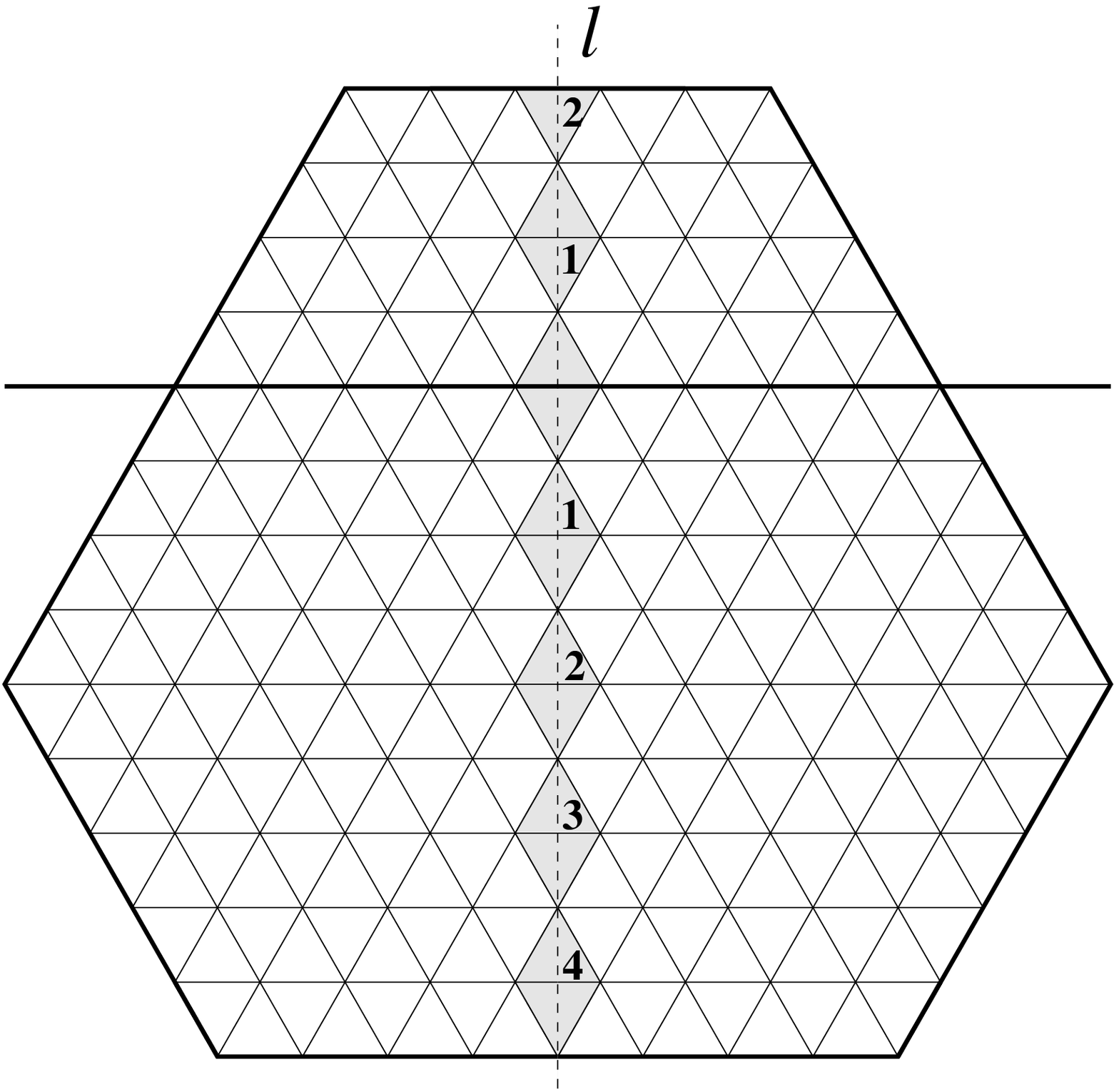}}
\centerline{{\smc Figure~1.1.} {\rm The hexagon $H(5,5,3)$; its vertebrae are shaded}}
\centerline{{\rm and labeled starting from the indicated horizontal reference line.}}
\endinsert

It is easy to see
that $H(a,b,k)$ contains precisely $k$ more unit triangles pointing upward than unit triangles 
pointing downward. Since a lozenge always covers one unit triangle of each kind, it
follows that $H(a,b,k)$ does not have any tilings for $k\ne0$.

This can be fixed by removing some unit triangles from $H(a,b,k)$, so as to restore the
balance between up- and down-pointing unit triangles. One natural way to do this, which 
was in fact our motivating example, is to remove an up-pointing triangular 
region of 
side-length $k$ from the center of $H(a,b,k)$ (the instance $k=1$ of this is treated
in \cite{2} and corresponds to a problem posed independently by Propp \cite{10} and 
Kuperberg (private communication)). At least from the viewpoint of our proof, it turns 
out that it is easier to enumerate tilings of the more general regions defined below.

We call a triangular subregion of $H(a,b,k)$ which is symmetric with respect to $\ell$ 
a {\it window}. If the vertex opposite the base of a window is above the base
we call the window a {\it $\Delta$-window}; otherwise we call it a {\it
$\nabla$-window}. A window is called even or odd according as its side-length is even or
odd.

\topinsert
\centerline{\mypic{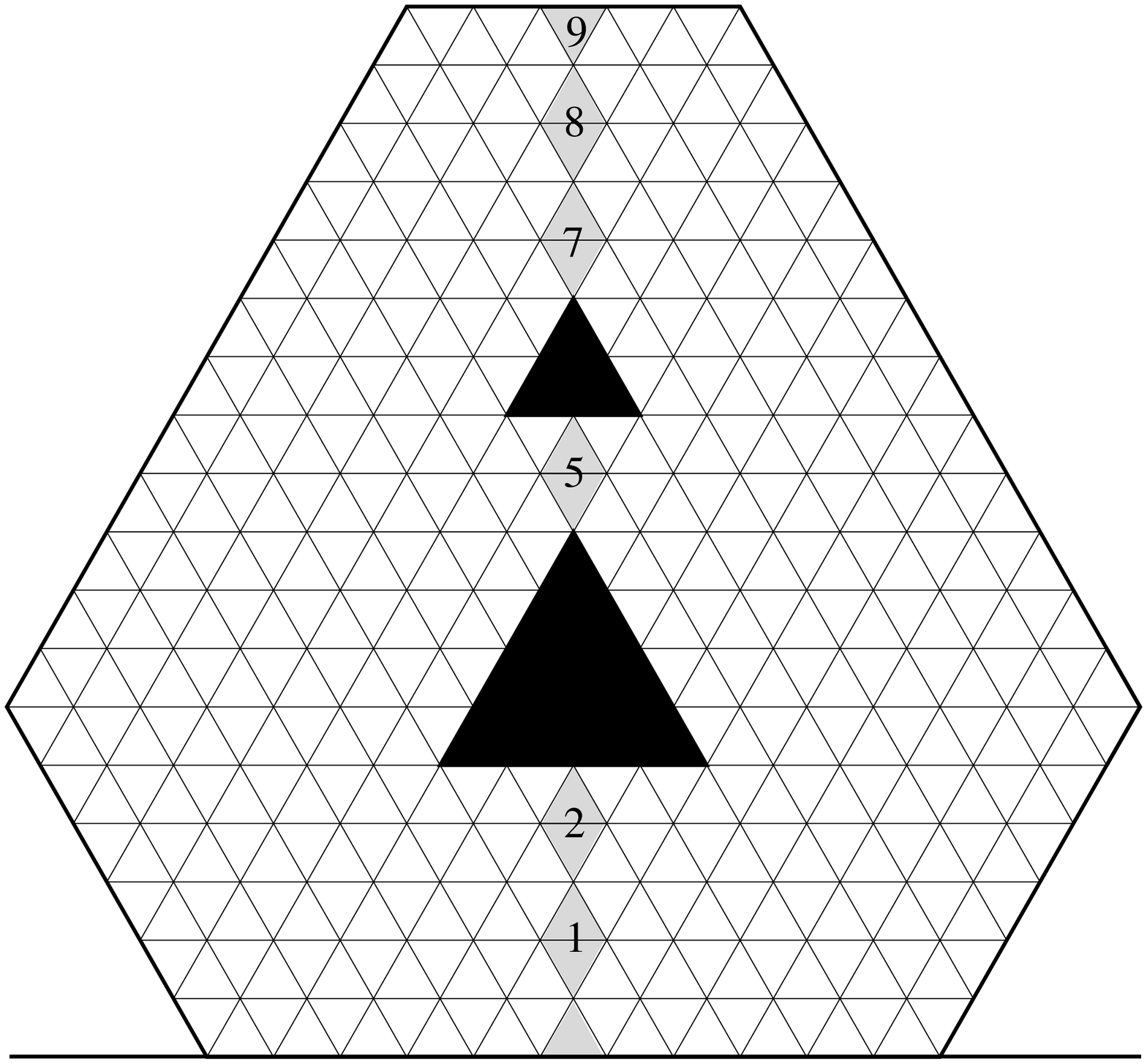}}
\centerline{{\smc Figure~1.2.} {\rm The region $H_{(1,2,5,7,8,9)}(5,6,6)$. The vertebrae}}
\centerline{{\rm contained in it are shaded, and their labels and the reference}}
\centerline{{\rm line for the labeling are shown.}}
\endinsert

The unit triangles of $H(a,b,k)$ along $\ell$ can be grouped in pairs forming lozenges --- 
we call them {\it rhombic vertebrae} ---
except possibly for those such triangles touching the boundary of $H(a,b,k)$, which we call 
{\it triangular vertebrae} (instances of both kinds of vertebrae occur in Figure 1.1; 
the vertebrae are indicated by a shading). 
The definitions below involve the concept of labeling the
vertebrae starting from a reference line --- a line of the triangular grid 
perpendicular to $\ell$: by this we mean that we label the vertebrae that fully lie 
on each
side of the reference line starting with 1 for the closest vertebra, 2 for the next,
and so on; see Figure 1.1 (in case the reference line is the base of $H(a,b,k)$ and the
vertebra touching it is triangular, we skip this vertebra and start by labeling the 
next one by 1).

Suppose $k$ is even. Remove from $H(a,b,k)$ some number of even $\Delta$-windows of total 
size $k$ (the total size of a collection of windows is the sum of the side-lengths of 
the windows). Since we removed only even windows, each vertebra is either fully removed
or fully left in place. Label the vertebrae of $H(a,b,k)$ starting from its base, and suppose 
the vertebrae in the leftover region have, in increasing order, labels $l_1,\dotsc,l_n$
(see Figure 1.2 for an example). Let $\bold l=(l_1,\dotsc,l_n)$. The region obtained
from $H(a,b,k)$ by removing windows as above may have forced lozenges (e.g., if some of the
windows touch). Remove all such forced lozenges. Denote the 
leftover region by $H_{\bold l}(a,b,k)$.

\topinsert
\centerline{\mypic{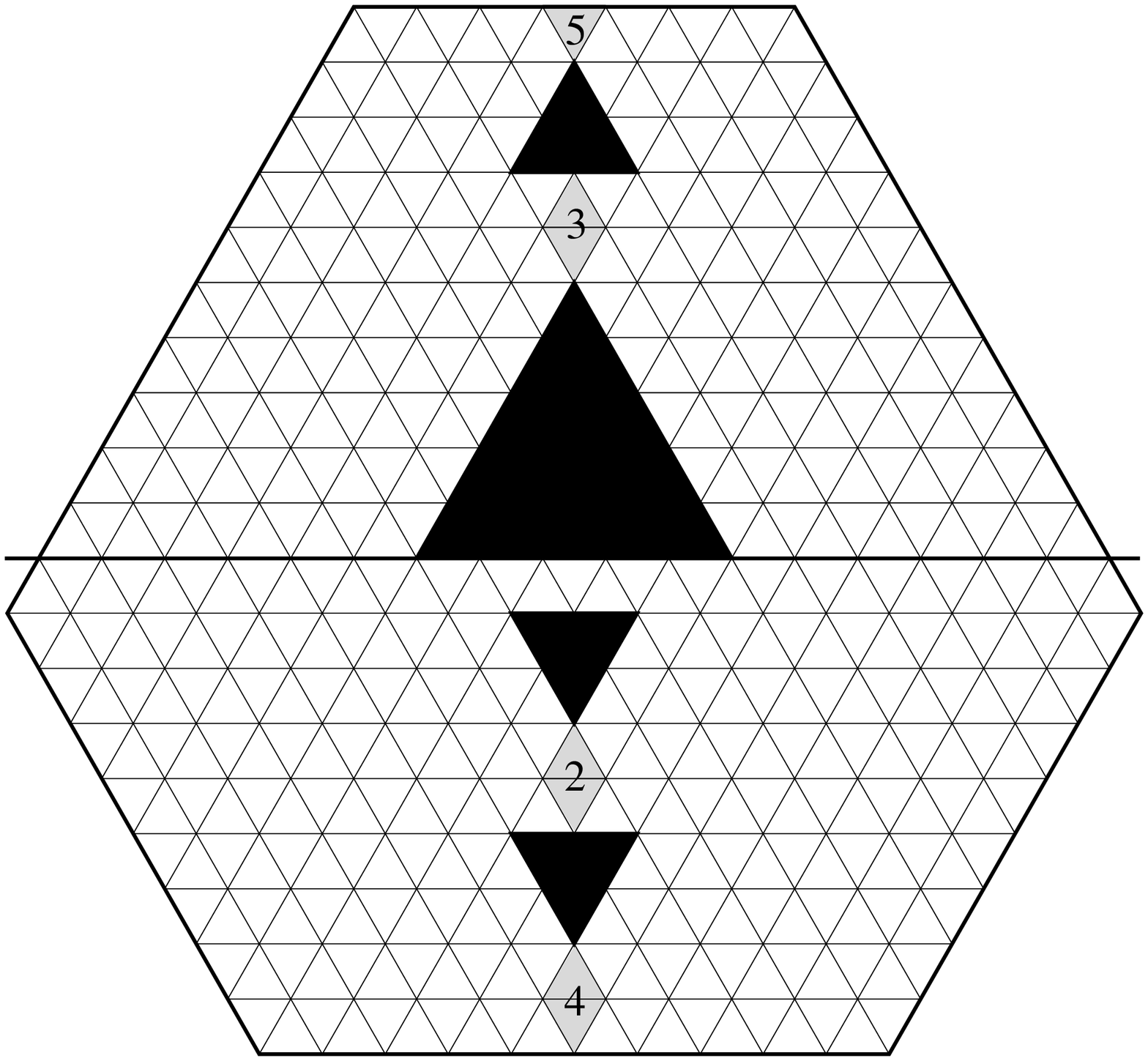}}
\centerline{{\smc Figure~1.3.} {\rm The region $H_{(2,4),(3,5)}(7,8,3)$. The vertebrae}}
\centerline{{contained in it, their labels and the reference line are shown.}}
\endinsert

\topinsert
\centerline{\mypic{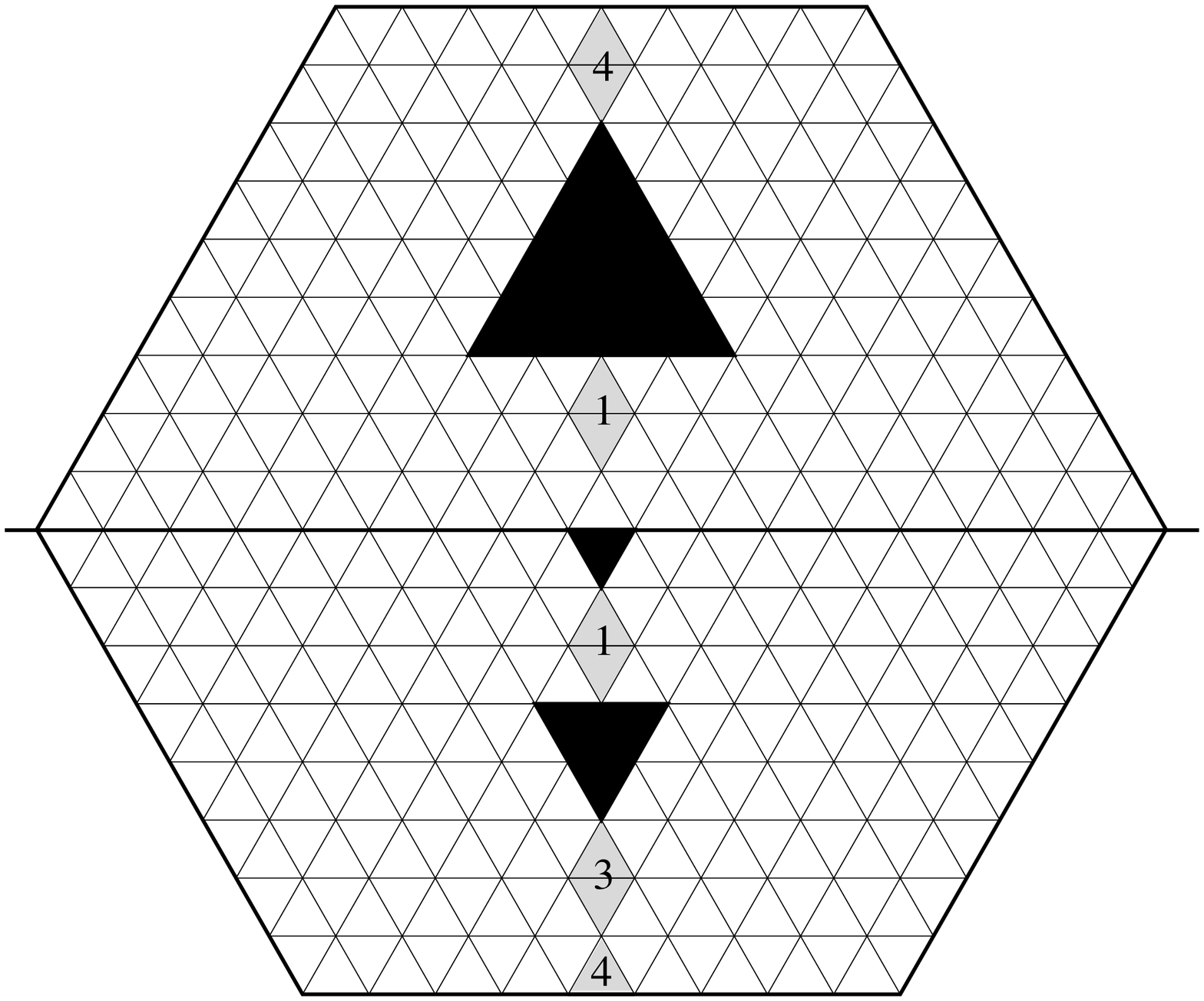}}
\centerline{{\smc Figure~1.4.} {\rm The region $\bar{H}_{(1,3,4),(1,4)}(8,8,1)$,}}
\centerline{{with the vertebrae it contains.}}
\endinsert

If $k$ is odd, start, from top to bottom, by removing a number of even
$\Delta$-windows from $H(a,b,k)$, then remove a single odd window (which can be either a
$\Delta$- or a $\nabla$-window) and continue downwards by removing some even
$\nabla$-windows (the set of removed even windows above the odd window, or the set of 
removed windows below it, or possibly both, may be empty). In this process, obey 
the restriction that the total size of the 
$\Delta$-windows is $k$ larger than the total size of the $\nabla$-windows. Choose the
base of the single odd window as our reference line, and label the vertebrae on both
sides of it as indicated above. Suppose the vertebrae below the reference line 
in the leftover region have, in increasing order, labels $l_1,\dotsc,l_m$, and the 
ones above it have, in increasing order, labels
$q_1,\dotsc,q_n$ (note that one or both of these lists may be empty). Let 
$\bold l=(l_1,\dotsc,l_m)$, $\bold q=(q_1,\dotsc,q_n)$. Remove any forced lozenges,
and denote the leftover
region by $H_{{\bold l},{\bold q}}(a,b,k)$ or $\bar{H}_{{\bold l},{\bold q}}(a,b,k)$,
according as the removed odd window was a $\Delta$- or a $\nabla$-window (Figures 1.3
and 1.4 show an example of each kind).

The regions $H_{\bold l}(a,b,k)$, $H_{{\bold l},{\bold q}}(a,b,k)$, and 
$\bar{H}_{{\bold l},{\bold q}}(a,b,k)$ are the regions whose number of tilings we will
determine, thus generalizing the case $b=c$ of MacMahon's result (which is obtained for 
$k=0$). 

Before giving the formulas expressing the number of tilings of these regions, we need to 
define two families of polynomials, both indexed by pairs of lists of positive integers. 

For $m,n\geq0$ define monic polynomials $B_{m,n}(x)$ and $\bar{B}_{m,n}(x)$ by the 
formulas

$$\align
B_{m,n}(x)&=2^{-mn-m(m-1)/2}(x+n+1)_m(x+n+2)_m\\
&(x+2)(x+3)^2\cdots(x+n-1)^2(x+n)\\
&\left(x+\frac{3}{2}\right)\left(x+\frac{5}{2}\right)^2\cdots
\left(x+\frac{2n-1}{2}\right)^2\left(x+\frac{2n+1}{2}\right)\\
&\prod_{i=1}^{n}\frac{(x+i)_m}{(x+i+1/2)_m}
\prod_{i=1}^{m}(2x+n+i+2)_{n+i-1}\tag1.1\\
\bar{B}_{m,n}(x)&=2^{-mn-n(n+1)/2}(x+m+1)_n\\
&(x+1)(x+2)^2\cdots(x+m-1)^2(x+m)\\
&\left(x+\frac{3}{2}\right)\left(x+\frac{5}{2}\right)^2\cdots
\left(x+\frac{2m-3}{2}\right)^2\left(x+\frac{2m-1}{2}\right)\\
&\prod_{i=1}^{m}\frac{(x+i)_n}{(x+i+1/2)_n}
\prod_{i=1}^{n}(2x+m+i+1)_{m+i},\tag1.2
\endalign$$
where $(a)_k:=a(a+1)\dotsc(a+k-1)$ is the shifted factorial (in the middle two lines of 
(1.1) and (1.2), the bases grow by one from each factor to the next, while the exponents 
increase by one midway through and then decrease by one to the end).

Next, given (possibly empty) lists of strictly increasing positive integers 
${\bold l}=(l_1,\dotsc,l_m)$ and ${\bold q}=(q_1,\dotsc,q_n)$, define constants 
$c_{{\bold l},{\bold q}}$ and $\bar{c}_{{\bold l},{\bold q}}$ by setting

$$\align
c_{{\bold l},{\bold q}}=2^{{n-m\choose2}-m}&\prod_{i=1}^{m}\frac{1}{(2l_i)\,!}
\prod_{i=1}^{n}\frac{1}{(2q_i-1)\,!}\frac{\prod_{1\leq i<j\leq m}(l_j-l_i)
\prod_{1\leq i<j\leq n}(q_j-q_i)}{\prod_{i=1}^{m}\prod_{j=1}^{n}(l_i+q_j)}\tag1.3\\
\bar{c}_{{\bold l},{\bold q}}=2^{{n-m\choose2}-m}&\prod_{i=1}^{m}\frac{1}{(2l_i-1)\,!}
\prod_{i=1}^{n}\frac{1}{(2q_i)\,!}\frac{\prod_{1\leq i<j\leq m}(l_j-l_i)
\prod_{1\leq i<j\leq n}(q_j-q_i)}{\prod_{i=1}^{m}\prod_{j=1}^{n}(l_i+q_j)}.\tag1.4
\endalign$$

Let $\lambda=\lambda({\bold l})$ and $\mu=\mu({\bold q})$ be the partitions having parts 
$l_m-m,l_{m-1}-m+1,\dotsc,l_1-1$ and $q_n-n,q_{n-1}-n+1,\dotsc,q_1-1$, respectively. 

We are now ready to define our polynomials $P_{{\bold l},{\bold q}}(x)$ and $\bar{P}_{{\bold l},{\bold q}}(x)$. They are given by the following formulas:

$$\align
P_{{\bold l},{\bold q}}(x-\lambda_1)&=c_{{\bold l},{\bold q}}B_{m,n}(x)\\
\prod_{c\in\lambda}(x&-h(c)+m+1)(x+h(c)+n+1)\prod_{c\in\mu}(x-h(c)+n+2)(x+h(c)+m)\tag1.5\\
\bar{P}_{{\bold l},{\bold q}}(x-\lambda_1)&=\bar{c}_{{\bold l},{\bold q}}\bar{B}_{m,n}(x)\\
\prod_{c\in\lambda}(x&-h(c)+m+1)(x+h(c)+n)\prod_{c\in\mu}(x-h(c)+n+1)(x+h(c)+m),\tag1.6
\endalign$$
where for a cell $c=(i,j)$ of a partition diagram, $h(c)$ is defined to be $i+j$ (in all 
formulas above empty products are considered equal to 1).

Denote by ${{\bold l}^{(i)}}$ the list obtained from ${\bold l}$ by omitting its $i$-th
element, and by ${\bold l}-1$ the list obtained from ${\bold l}$ by subtracting 1 from 
each of its elements (except in case $l_1=1$, when ${\bold l}-1$ is defined to consist of 
$l_2-1,\cdots,l_m-1$).

\medskip
Our generalization of MacMahon's enumeration of plane partitions contained in a symmetric box
is given in Theorem 1.1 below. Its proof is deduced in Section 3 from Proposition 2.1, which
interprets the above-defined polynomials $P_{{\bold q},{\bold l}}(x)$ and 
$\bar{P}_{{\bold q},{\bold l}}(x)$ as tiling enumerators of two families of lattice regions 
defined in Section 2. Proposition 2.1 and Theorem 1.1 are the main
results of this paper.
 
Let $\M(R)$ denote the number of lozenge tilings of the region $R$ (the ``M'' stands for
``matchings,'' as lozenge tilings can be identified with perfect matchings on the dual graph).

\proclaim{Theorem 1.1} The number of lozenge tilings of the regions  
$H_{\bold l}(a,b,k)$, $H_{{\bold l},{\bold q}}(a,b,k)$, and 
$\bar{H}_{{\bold l},{\bold q}}(a,b,k)$ is given by the following formulas:

$(a)$ For $k$ even, we have

$$
2^{-m}\M(H_{\bold l}(a,b,k))=\left\{
\aligned &P_{\emptyset,{\bold l}}((a+k-2)/2)
\bar{P}_{{\bold l}-1,\emptyset}(a/2),\ \ \ \ \ \ \ \ \ \ \ \ \text{$a$ even, $l_1=1$},\\
&P_{\emptyset,{\bold l}}((a+k-2)/2)
P_{{\bold l}-1,\emptyset}(a/2),\ \ \ \ \ \ \ \ \ \ \ \ \text{$a$ even, $l_1>1$},\\
&P_{\emptyset,{\bold l}^{(m)}}((a+k-1)/2)
\bar{P}_{{\bold l},\emptyset}((a-1)/2),\ \ \ \ \!\,\,\text{$a$ odd}.
\endaligned
\right.\tag1.7
$$

$(b)$ For $k$ odd, we have

$$
2^{-m-n}\M(H_{{\bold l},{\bold q}}(a,b,k))=\left\{
\aligned &\bar{P}_{{\bold l},{\bold q}}((a+k-1)/2)
P_{{\bold q},{\bold l}^{(m)}}(a/2),\ \ \ \ \ \text{$a$ even},\\
&\bar{P}_{{\bold l},{\bold q}^{(n)}}((a+k)/2)
P_{{\bold q},{\bold l}}((a-1)/2),\ \ \ \,\text{$a$ odd},
\endaligned
\right.\tag1.8
$$

and

$$
2^{-m-n}\M(\bar{H}_{{\bold l},{\bold q}}(a,b,k))=\left\{
\aligned &P_{{\bold l},{\bold q}}((a+k-1)/2)
\bar{P}_{{\bold q},{\bold l}^{(m)}}(a/2),\ \ \ \ \ \text{$a$ even},\\
&P_{{\bold l},{\bold q}^{(n)}}((a+k)/2)
\bar{P}_{{\bold q},{\bold l}}((a-1)/2),\ \ \ \,\text{$a$ odd}.
\endaligned
\right.\tag1.9
$$

\endproclaim

\medskip
\flushpar
{\smc Remark 1.2.} For $k\ne0$, the bijection between tilings of our regions and
stacks of unit cubes breaks down. Indeed, trying to lift such a tiling to three
dimensions one is lead to an ``impossible'' construction, similar to Penrose's 
impossible triangle (see Figure 1.5).

\medskip
In Section 3 we show how to reduce the statement of Theorem 1.1 to the problem of 
finding the tiling generating function of certain simply-connected regions, defined 
in Section 2. In Sections 4 and 5
we prove, using an inductive argument, that these tiling generating functions are 
given by the above-defined polynomials $P$ and $\bar{P}$. The idea of the proof is very 
simple: using combinatorial considerations, we deduce recurrence relations for the 
tiling generating functions of our regions, and we show that the polynomials $P$ and 
$\bar{P}$ satisfy the same recurrences. This approach, however, raises the question of 
how could one guess that the tiling generating functions of our regions are given by the 
intricate formulas (1.1)--(1.6). We explain in Section 6 how this can be done in a 
conceptual way, using some computer assistance.

\topinsert
\centerline{\mypic{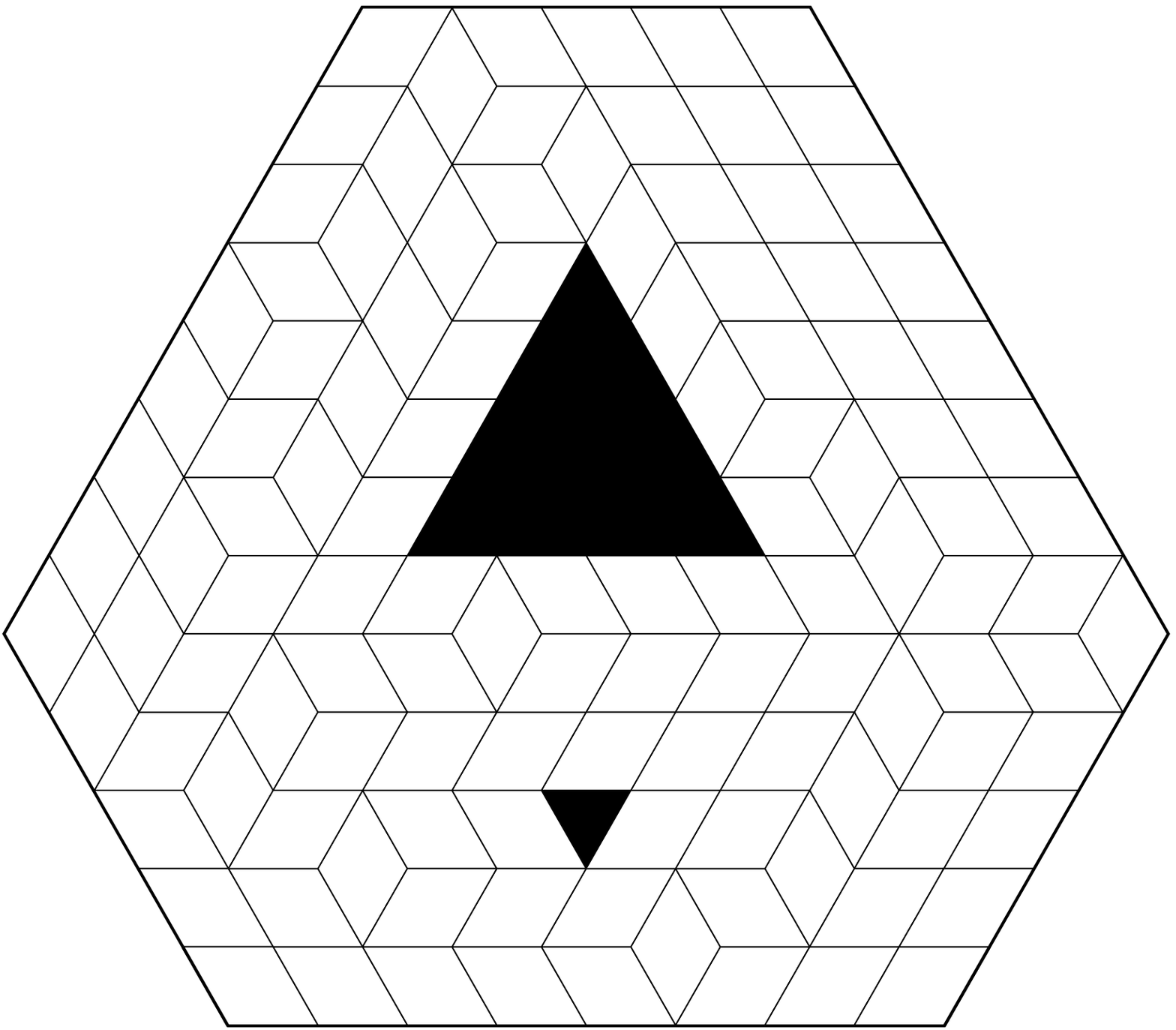}}
\centerline{{\smc Figure~1.5.} {\rm A tiling of $H_{(1),(1,4,5)}(5,5,3)$.}}
\endinsert

\mysec{2. Two families of regions}

\medskip
We now define two families of regions that will turn out to be closely related to the
polynomials $P$ and $\bar{P}$ given by (1.5) and (1.6).  

Consider the triangular lattice drawn so that one of the directions of the lattice lines 
is horizontal. From each point of the lattice, there are six possible steps one can take
on the lattice. We call the four non-horizontal steps, according to their approximate 
cardinal direction, northeast, northwest, southwest and southeast.

\topinsert
\centerline{\mypic{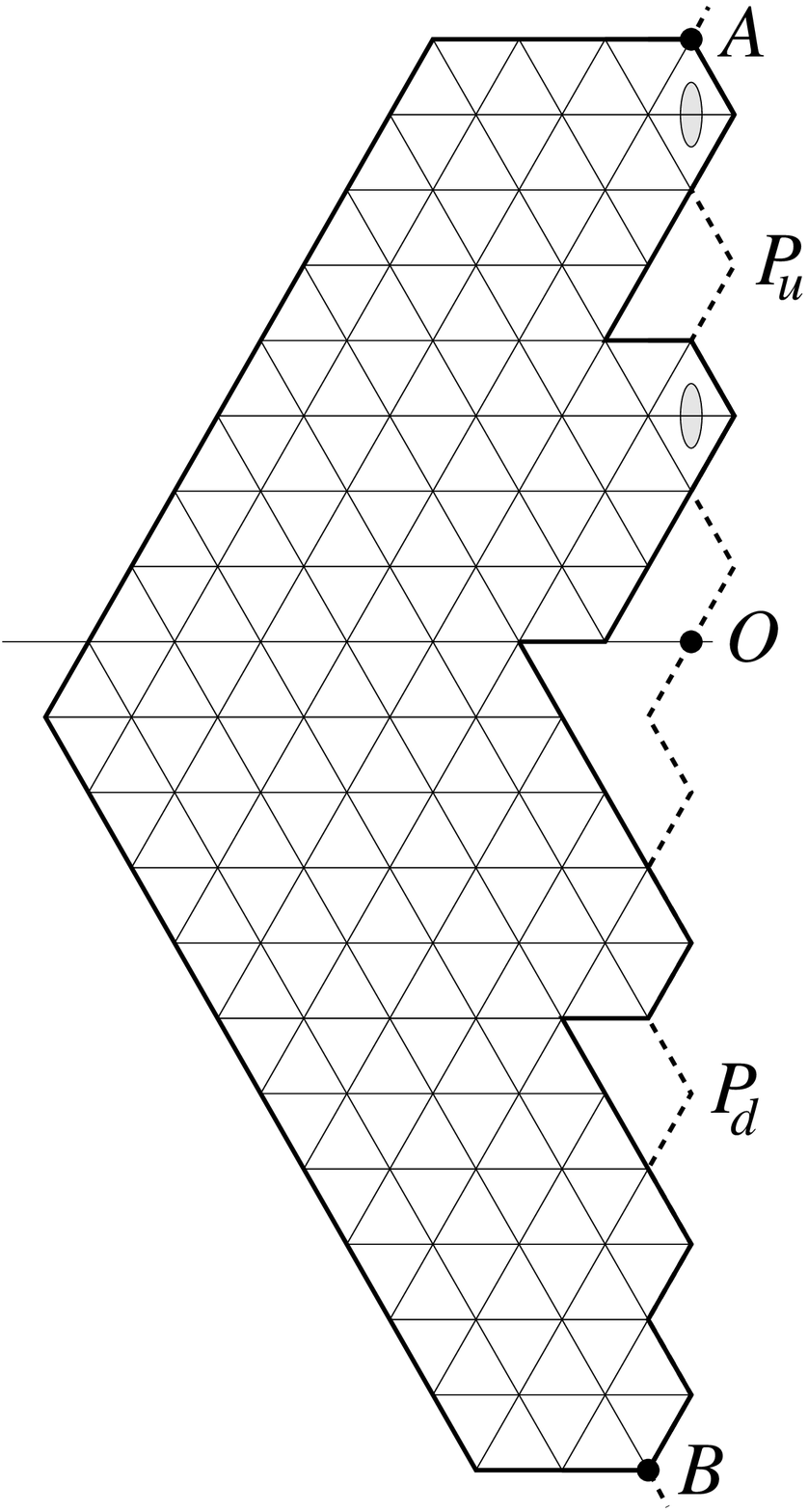}}
\centerline{{\smc Figure~2.1{\rm (a). $R_{(2,4,5),(2,4)}(2)$.}}}
\endinsert

Let $O$ be a point of the triangular lattice and consider the infinite lattice paths 
$P_u$ and $P_d$ that start from $O$ and take alternate steps northeast and northwest,
respectively southwest and southeast (these paths are shown in dotted lines in 
Figure 2.1(a)). We say that two consecutive 
steps of $P_u$, the first going
northeast and the second northwest, form a {\it bump}; similarly, a step southeast
followed by a step southwest on $P_d$ form a bump. Label the bumps of $P_u$ and $P_d$, 
starting with the ones closest to $O$, consecutively by $1,2,\dotsc$. Select 
arbitrary finite subsets of the bumps of $P_d$ and $P_u$, and suppose their labels, 
in increasing order, form (possibly empty) lists ${\bold l}=(l_1,\dotsc,l_m)$ and 
${\bold q}=(q_1,\dotsc,q_n)$, respectively. 
Assume that at least one of ${\bold l}$ and ${\bold q}$ is non-empty. 

Starting from these selected bumps, we define the ``right boundary'' of a region 
$R_{{\bold l},{\bold q}}(x)$
as follows (here $x$ is a non-negative integer on which a lower bound will be imposed 
below). Call the uppermost and lowermost points on a bump the top and bottom
of that bump, respectively. For each selected bump on $P_u$, draw a horizontal ray to the
left of its top, and a ray going southwest from its bottom. For each selected bump on 
$P_d$, draw a ray going left from its bottom and a ray going northwest from its top. 

By following these rays, we join the selected bumps on $P_u$ in a connected piece: follow
the horizontal ray from the top of bump $q_i$ to its intersection point with the ray 
starting from the bottom $b$ of bump $q_{i+1}$, and then use the latter ray to reach $b$. 
The selected bumps on $P_d$ are connected into a single piece in a similar way. 

To join the above two pieces together, we draw a horizontal ray left of $O$. We follow the
ray going southwest from the bottom of bump $q_1$ of $P_u$ until its intersection
with the horizontal ray originating at $O$. Then we follow this latter ray until it hits
the ray going northwest from the top $a$ of the $l_1$-th bump of $P_d$, which we then track
until we reach $a$. We have now included all selected bumps in a connected path $Q$ joining
the top $A$ of the last selected bump on $P_u$ to the bottom $B$ of the last selected bump
on $P_d$ (in case ${\bold l}$ is empty, $B$ is taken to be the intersection of the ray going
southwest from the bottom of the first selected bump on $P_u$ with the horizontal through 
$O$; for ${\bold q}=\emptyset$, $A$ is chosed analogously; see Figures 2.1(b) and 2.1(c)). 
It is easy to see that this path 
has exactly $2l_m-m+n+1$ steps going southeast (the one exception is the case
${\bold l}=\emptyset$, when it has $2l_m-m+n=n$ steps).

\topinsert
\twoline{\mypic{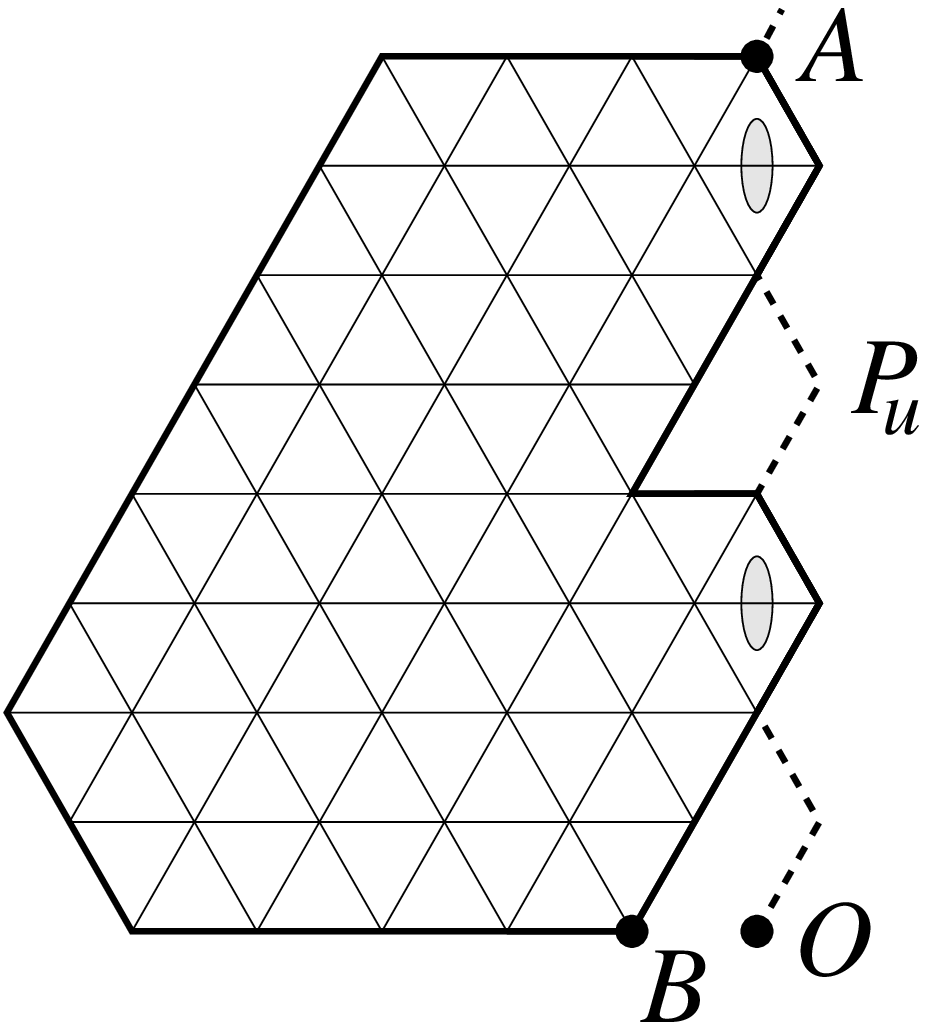}}{\mypic{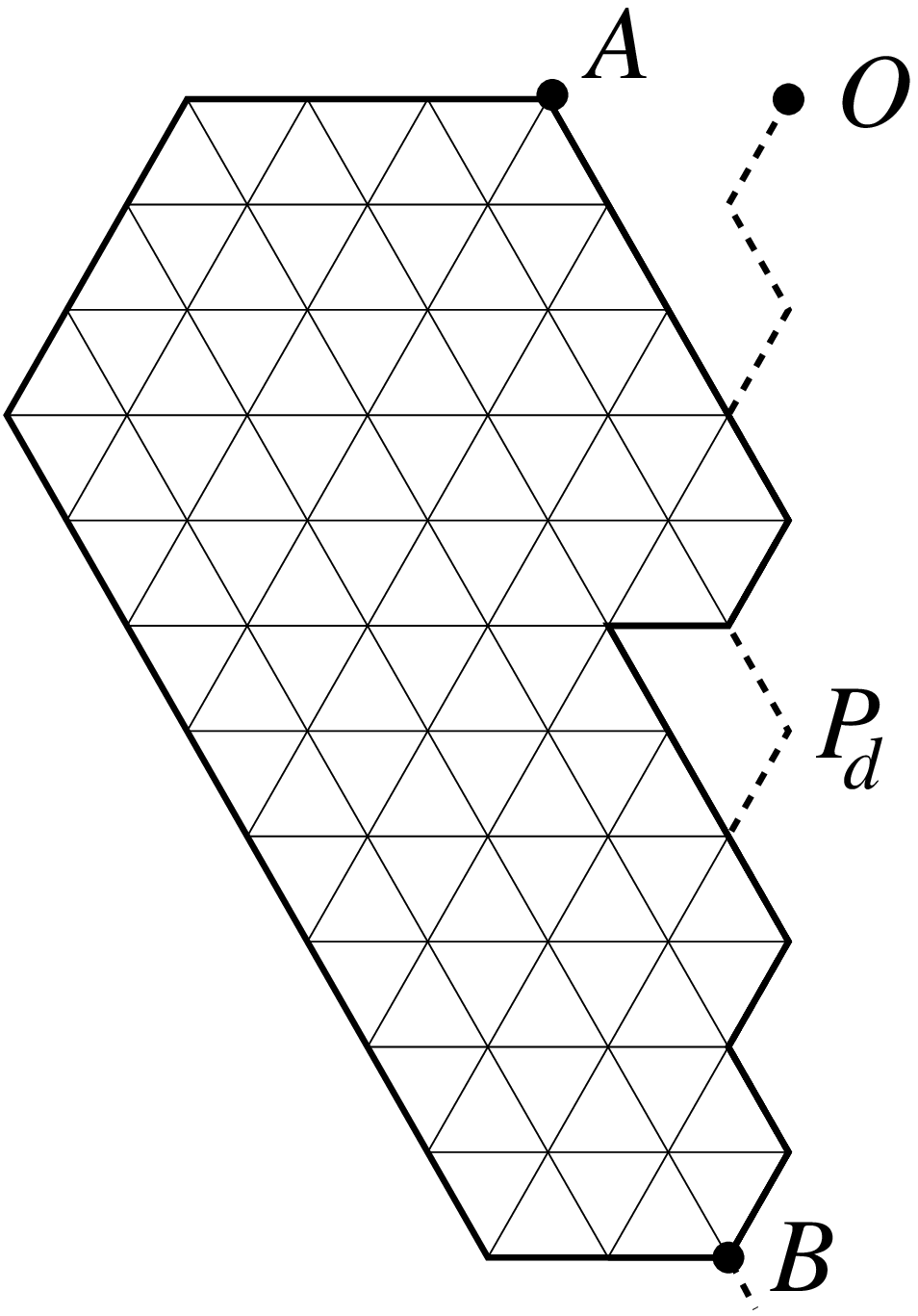}}
\medskip
\twoline{Figure~2.1{\rm (b). $R_{\emptyset,(2,4)}(4)$.}}{Figure~2.1{\rm (c). 
$R_{(2,4,5),\emptyset}(2)$.}}
\endinsert

To close the boundary of $R_{{\bold l},{\bold q}}(x)$, we continue from $B$ by taking
$x$ steps left (in case ${\bold l}$ is empty, we move left until we are $x+1$ units from $O$).
Then we head northwest for $2l_m-m+n+1$ steps (one fewer when ${\bold l}=\emptyset$), 
after which we turn and go northeast 
until we hit the horizontal ray going left from $A$ (see Figure 2.1(a)). Finally, take 
$x+(l_m-m)-(q_n-n)+1$ steps east on this ray, until we reach $A$.
It is not hard to see that the length of the southwest portion of this
boundary must be chosen as indicated (i.e.,
to match the number of southeast steps on $Q$) in order for the enclosed 
region to admit lozenge tilings (this follows for example by encoding tilings as families of
non-intersecting lattice paths, as described at the beginning of Section 4).

In general we will allow the lozenges in a tiling to be weighted. The weight of a tiling 
is the product of the weights on its lozenges. The sum of weights of all tilings of a 
region $R$ is denoted $\M(R)$ and is called the tiling generating function of $R$.

We define $R_{{\bold l},{\bold q}}(x)$ to be the region enclosed by the above-defined
boundary, requiring in addition
that the tile positions fitting in the selected bumps on $P_u$ are weighted
by $1/2$ (i.e., if one of these positions is occupied by a lozenge in some tiling, that
lozenge gets weight $1/2$ --- this will be indicated in our figures by a shaded oval 
placed on such tile positions; all other tiles get weight 1). Since the top and bottom 
sides have
non-negative length, $x$ must satisfy
$x\geq\max\{0,q_n-l_m-n+m-1\}$ (and $x\geq\max\{-1,q_n-n-1\}=q_n-n-1$, for 
${\bold l}=\emptyset$). Define $R_{\emptyset,\emptyset}(x)=\emptyset$ for all $x$.

\topinsert
\centerline{\mypic{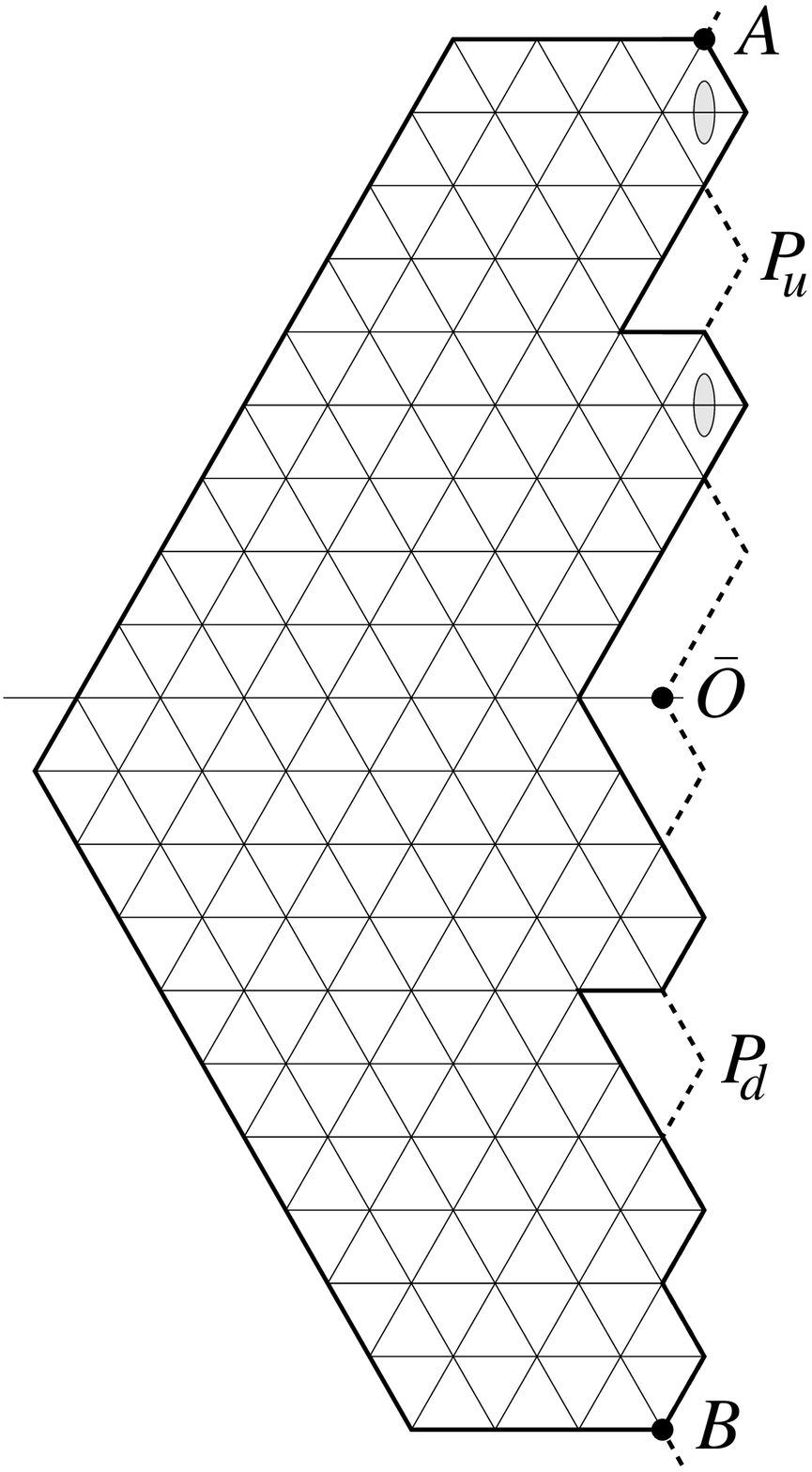}}
\bigskip
\centerline{{\smc Figure~2.2{\rm (a). $\bar{R}_{(2,4,5),(2,4)}(3)$.}}}
\endinsert

\topinsert
\bigskip
\twoline{\mypic{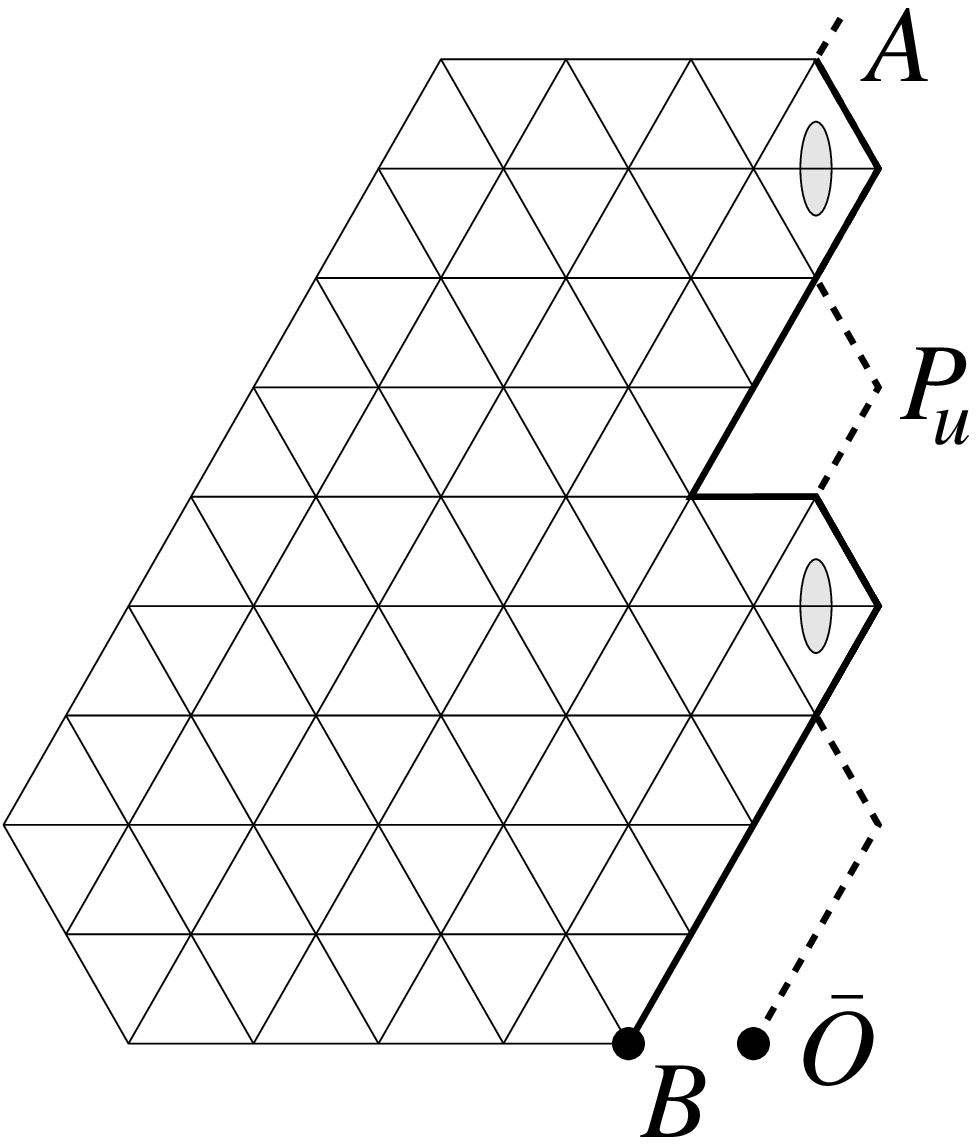}}{\mypic{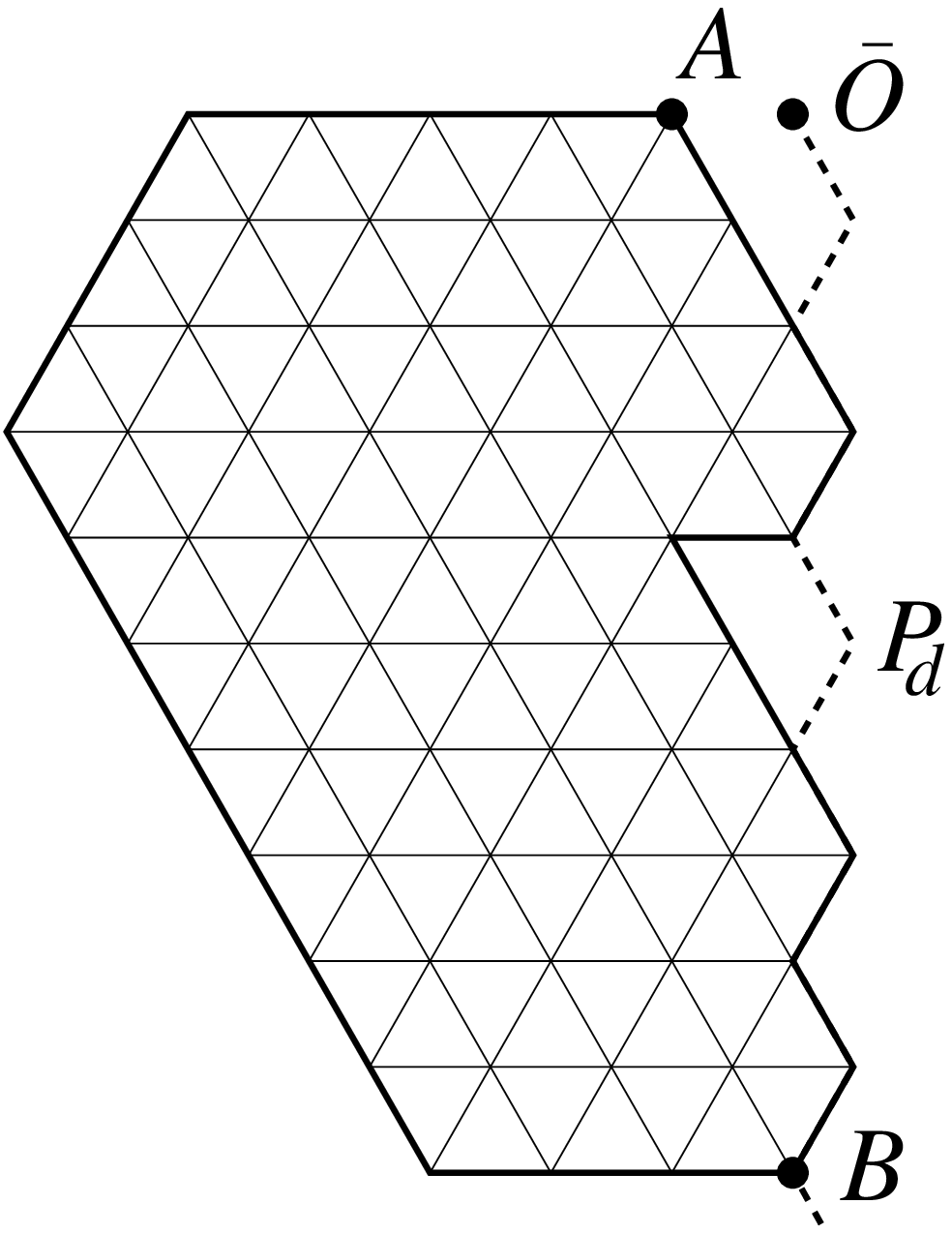}}
\bigskip
\twoline{Figure~2.2{\rm (b). $\bar{R}_{\emptyset,(2,4)}(5)$.}}{Figure~2.2{\rm (c). 
$\bar{R}_{(2,4,5),\emptyset}(3)$.}}
\endinsert

Our second family of regions, denoted $\bar{R}_{{\bold l},{\bold q}}(x)$, is defined
almost identically (see Figures 2.2(a)--(c)). The only difference is that, in the 
construction of the path $Q$ that
connects the selected bumps, instead of the horizontal ray starting from $O$, we use the
one starting from the lattice point $\bar{O}$, one step southwest of $O$ (and, in
case ${\bold l}$ is empty, when defining the bottom side of 
$\bar{R}_{\emptyset,{\bold q}}(x)$, we move left until we are $x$ units away from 
$\bar{O}$). As a consequence,
this path has now one fewer southeast steps than before, and hence we define the southwest
boundary of $\bar{R}_{{\bold l},{\bold q}}(x)$ to have length $2l_m-m+n$. Also, the top
side of the enclosed region has now length $x+(l_m-m)-(q_n-n)$. The parameter $x$
satisfies therefore $x\geq\max\{0,q_n-l_m-n+m\}$ (this time the condition on $x$ remains 
the same for ${\bold l}=\emptyset$). Again, we define 
$\bar{R}_{\emptyset,\emptyset}(x)=\emptyset$ for all $x$.

The close connection between the above-defined regions and the polynomials $P$ and 
$\bar{P}$ given by (1.5) and (1.6) is expressed by the following result.

\proclaim{Proposition 2.1} For all lists ${\bold l}$ and ${\bold q}$ of strictly 
increasing positive integers we have

$$\align
\M(R_{{\bold l},{\bold q}}(x))&=P_{{\bold l},{\bold q}}(x)\tag2.1\\
\M(\bar{R}_{{\bold l},{\bold q}}(x))&=\bar{P}_{{\bold l},{\bold q}}(x),\tag2.2
\endalign$$
for all $x$ for which the regions on the left hand side of these equalities are defined.
\endproclaim

In the next section we show how the above Proposition implies Theorem 1.1. The proof
of Proposition 2.1 is presented in Sections 4 and 5. 
 
\mysec{3. Reduction to simply-connected regions}

\medskip
One useful way to approach certain tiling enumeration problems (lozenge tilings 
in particular) is to identify them with families of non-intersecting lattice paths, and 
then use the 
Gessel-Viennot theorem (see e.g. \cite{12, Theorem 1.2} or \cite{5}) to express the 
number we seek as a 
determinant. While this could be carried out directly for the regions on the left 
hand sides of (1.7)--(1.9) that we are concerned with, it turns out that it will be
more effective to reduce first our problem to the case of certain simply connected subregions,
for the number of tilings of which we can easily obtain recurrences. This reduction can be 
achieved using the Factorization Theorem for perfect matchings presented in 
\cite{1,Theorem 1.2}.

Let $R$ be a region of the triangular lattice, and suppose $R$ has a vertical symmetry 
axis $\ell$ (an example is shown in Figure 3.1). Let $P_l$ and $P_r$ be the two zig-zag 
lattice paths taking alternate steps northeast and northwest and touching $\ell$ from
left and right, respectively. Group the unit triangles of $R$ crossed by $\ell$ 
in sequences of contiguous triangles. These are separated by sequences of 
triangles of the lattice not contained in $R$; we call these sequences gaps.

\topinsert
\centerline{\mypic{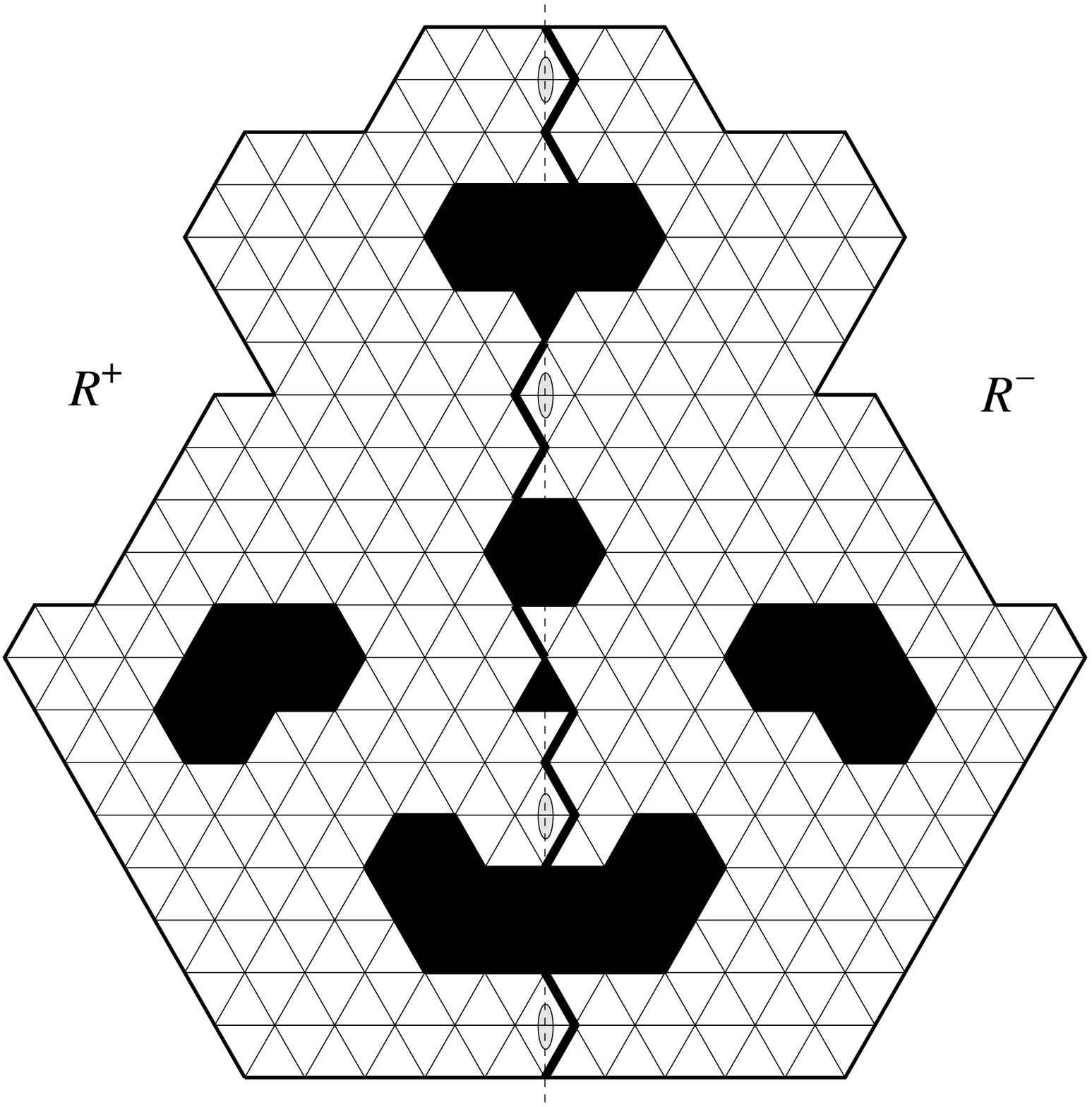}}
\centerline{{\smc Figure~3.1.} {\rm }}
\endinsert

We cut $R$ into two subregions $R^+$ and $R^-$ as follows. Consider the paths $P_l$ 
and $P_r$. Start at top by following $P_r$ to the first gap. After each gap, keep 
following the path along which we reached it or switch to the other path, according as 
the gap contains an even or an odd number of unit triangles. Continue this until we 
reach the bottom of $R$. Finally, assign weight $1/2$ to all tile positions in $R$ that 
have two sides contained in our cutting path. Define 
$R^+$ and $R^-$ to be the regions obtained to the left and right of our cutting path, 
respectively, and having the above-mentioned lozenge positions weighted by $1/2$.  

Since we are interested in counting lozenge tilings, we
may assume that $R$ has an even number of unit triangles crossed by $\ell$ (indeed, 
if this number was odd, the symmetry of $R$ would imply that $R$ contains an odd number of 
unit triangles, so it would have no lozenge tilings); define the {\it width} of $R$, denoted 
$\w(R)$, to be half this number.

Then the Factorization Theorem for perfect matchings of \cite{1,Theorem 1.2} implies

$$\M(R)=2^{\w(R)}\M(R^+)\M(R^-).\tag3.1$$

Indeed, lozenge tilings of $R$ can be identified with perfect matchings of the dual graph
$G$ of $R$, i.e., the graph whose vertices are the unit triangles of $R$ and whose edges
connect precisely those pairs of triangles that share an edge. Since $R$ is symmetric with 
respect to $\ell$, so is $G$, and all conditions in the hypothesis of 
\cite{1, Theorem 1.2} are easily seen to be met. We obtain that $\M(G)=2^k\M(G^+)\M(G^-)$,
where $G^+$ and $G^-$ are certain precisely defined subgraphs,
and $\M(G)$ denotes the matching generating function of $G$. However, it is 
easily checked that $G^+$ and $G^-$ are exactly the dual graphs of the regions $R^+$ and 
$R^-$, respectively, thus proving (2.1).

\medskip
{\it Proof of Theorem 1.1.} All equalities in the statement of Theorem 1.1 follow directly
by applying (2.1) to the regions on the left hand side of (1.7)--(1.9) and using 
Proposition~2.1. More precisely, in all cases, the absolute value of 
the exponent in the power of 2 in (1.7)--(1.9)
turns out to be the width of the corresponding region, and, by Proposition 2.1, the two factors in the product
on the right hand side turn out to be equal to the tiling generating 
functions of the two regions that arise by applying the Factorization Theorem to our regions.

(a) Take $R$ to be the region $H_{\bold l}(a,b,k)$, and apply (2.1). Since all the gaps
are even, the cutting path determining $R^+$ and $R^-$ lies fully on the right side of 
$\ell$ (this is illustrated in Figures 3.2(a)--(c)). Suppose $a$ is even. Since $a+k$ is 
also even, this cutting path starts at the
midpoint of the top side of $R$ and ends at the midpoint of its base (see Figures 3.2(a) and 
(b)). Consequently, all
unit triangles of $R$ crossed by $\ell$ can be paired to form rhombic vertebrae. By our
definition, these are precisely the $l_1$-th,$\,\dotsc$,$\,l_m$-th vertebrae of the hexagon 
$H$ enclosed by the outer boundary of $R$. Therefore, the region $R^+$ to the left of the 
cutting path is precisely $R_{\emptyset,{\bold l}}((a+k-2)/2)$.

As for the region $R^-$ to the right of the cutting path, suppose first that $l_1>1$ (see 
Figure 3.2(a)). Then
it is easy to check that the region obtained from $R^-$ after removing the forced lozenges
is precisely $\bar{R}_{(l_1-1,\dotsc,l_m-1)}(a/2)$, rotated by 180 degrees. 
On the other hand, if $l_1=1$, because of the different pattern of forced lozenges in $R^-$
(extra forced lozenges on the bottom, as shown in Figure 3.2(b)), the subregion left after 
removing the forced 
lozenges is isomorphic to $R_{(l_2-1,\dotsc,l_m-1)}(a/2)$. Since the width of $R$ is $m$, 
this proves, by Proposition 2.1, the first two equalities in (1.7).

\topinsert
\centerline{\mypic{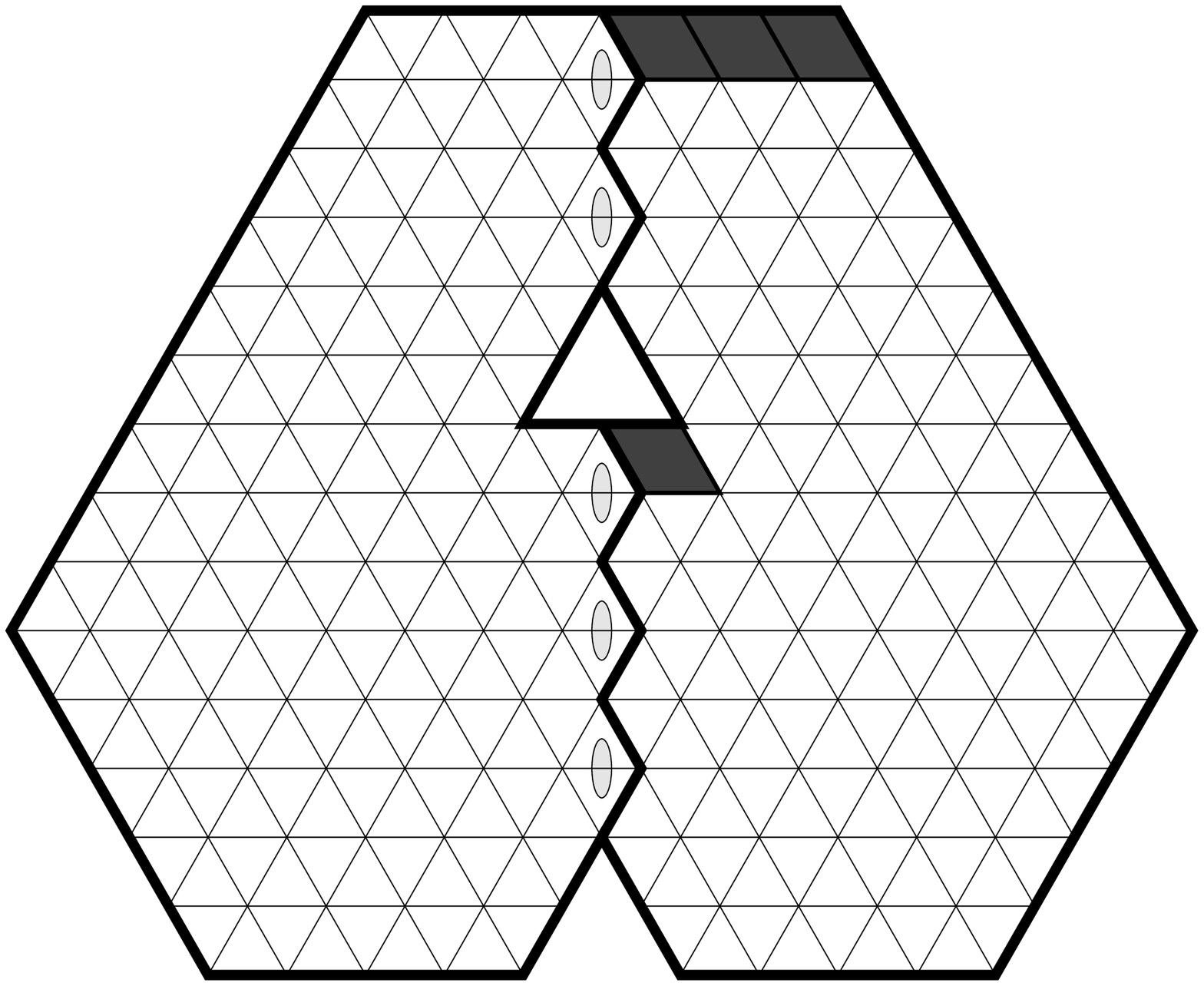}}
\centerline{{\smc Figure~3.2{\rm (a).}}}
\centerline{$H_{(2,3,4,6,7)}(6,5,4)$ reduces to}
\centerline{$R_{\emptyset,(2,3,4,6,7)}(4)$ and $R_{(1,2,3,5,6),\emptyset}(3)$.}
\endinsert

\topinsert
\twoline{\mypic{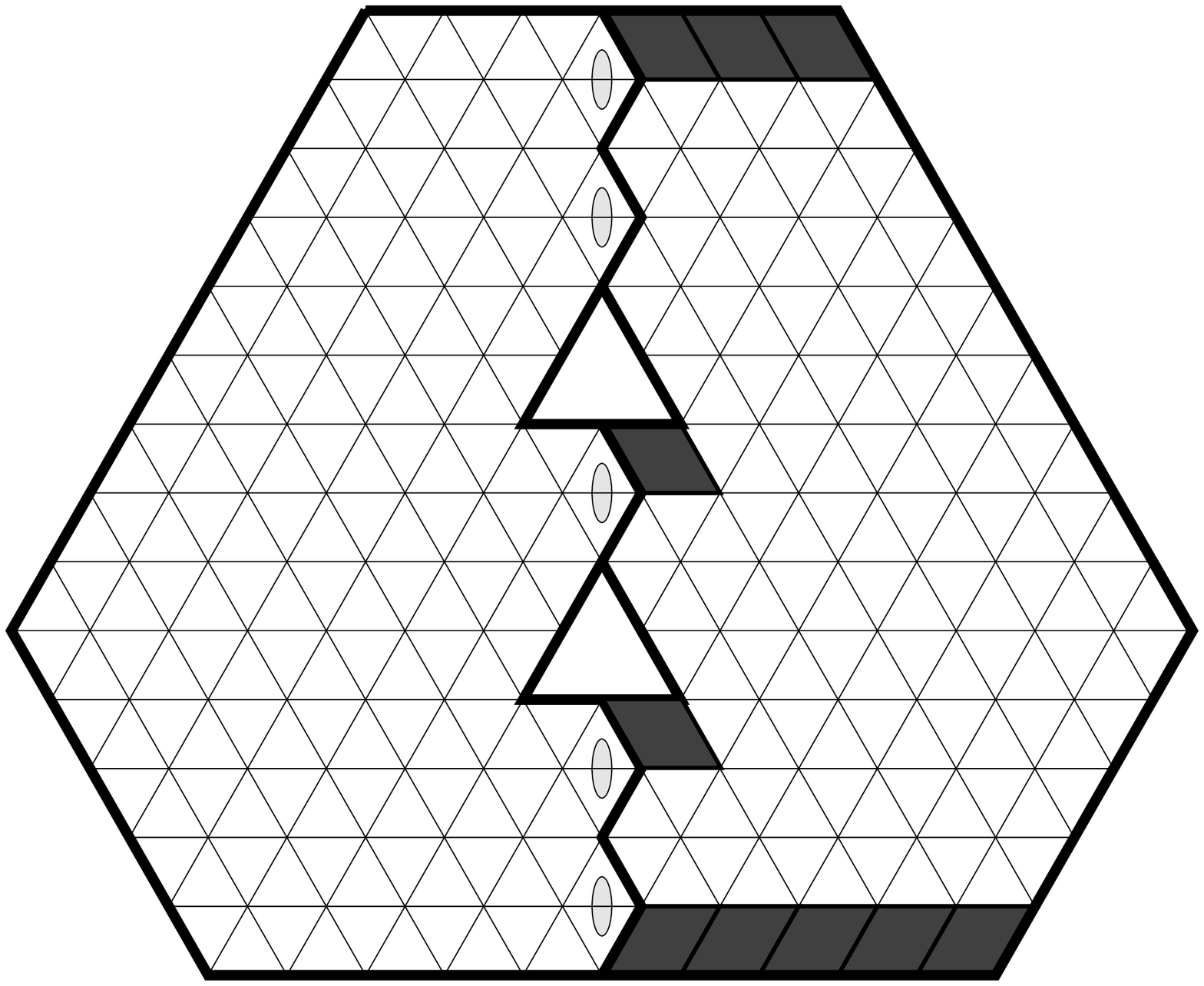}}{\mypic{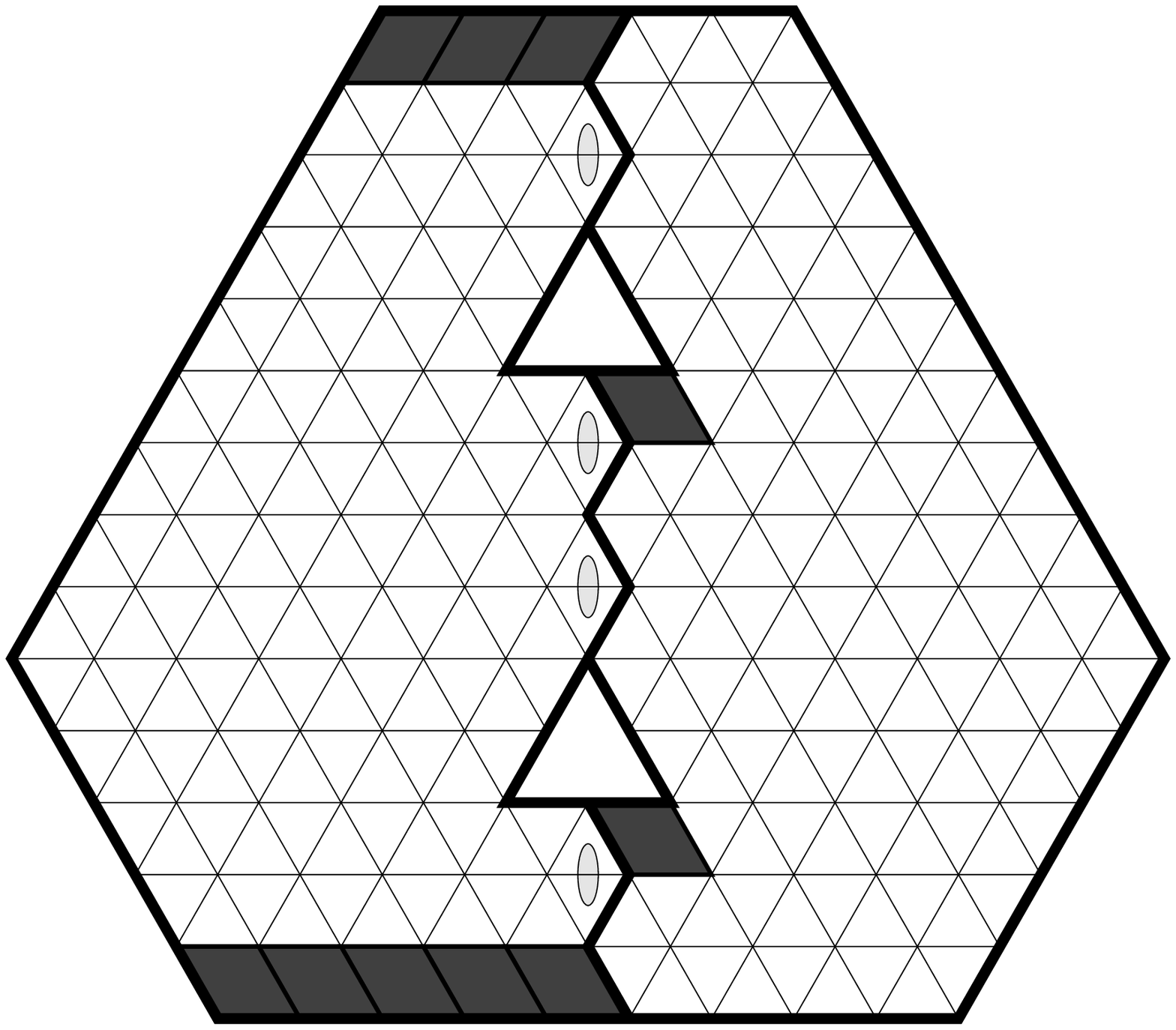}}
\twoline{Figure~3.2{\rm (b).}}{Figure~3.2{\rm (c).}} 
\twoline{{\rm $H_{(1,2,4,6,7)}(6,5,4)$ reduces to}}{{\rm $H_{(1,3,4,6,7)}(5,5,4)$ reduces to}}
\twoline{{\rm $R_{\emptyset,(1,2,4,6,7)}(4)$ and $\bar{R}_{(1,3,5,6),\emptyset}(3)$.}}
{{\rm $R_{\emptyset,(1,3,4,6)}(4)$ and $\bar{R}_{(1,3,4,6,7),\emptyset}(2)$.}}
\endinsert

For $a$ odd we proceed similarly. 
Now the last labeled vertebra of $H$ is a triangular vertebra, and it remains 
present for all choices of removing even $\Delta$-windows from $H$ (see Figure 3.2(c)). 
By definition, its label
is $l_m$. This fact and an examination of the pattern of forced lozenges in
$R^+$ and $R^-$ explains why, at the indices of the two factors on the right hand side of the 
third equality in (1.7), the list ${\bold l}$ appears once with its last entry omitted and
once in full. Indeed, because of forced lozenges at the top of $R^+$, the triangular vertebra 
labeled $l_m$ is not relevant for the region obtained from $R^+$ by removing the forced
lozenges. This leftover region is easily seen to be exactly  
$R_{\emptyset,{\bold l}^{(m)}}((a+k-1)/2)$. Furthermore, the region obtained from $R^-$
by removing the forced lozenges is isomorphic to $\bar{R}_{{\bold l},\emptyset}((a-1)/2)$,
irrespective of the value of $l_1$ (the difference in this regard from the case $a$ 
even is explained by
the fact that the pattern of forced lozenges in $R^-$ depends on $l_1$ for $a$ even, but does
not for $a$ odd). Since the width of $R$ is $m$, by Proposition 2.1 this completes the proof 
of (1.7).

(b) To prove (1.8), take $R$ to be the region $H_{{\bold l},{\bold q}}(a,b,k)$ and apply 
(2.1). The resulting cutting path stays on the right of $\ell$ until it arrives to the odd
$\Delta$-window, then it switches to the left of $\ell$ and stays there until it reaches the 
bottom of $R$ (this is illustrated in Figures~3.3(a) and (b)).

\topinsert
\twoline{\mypic{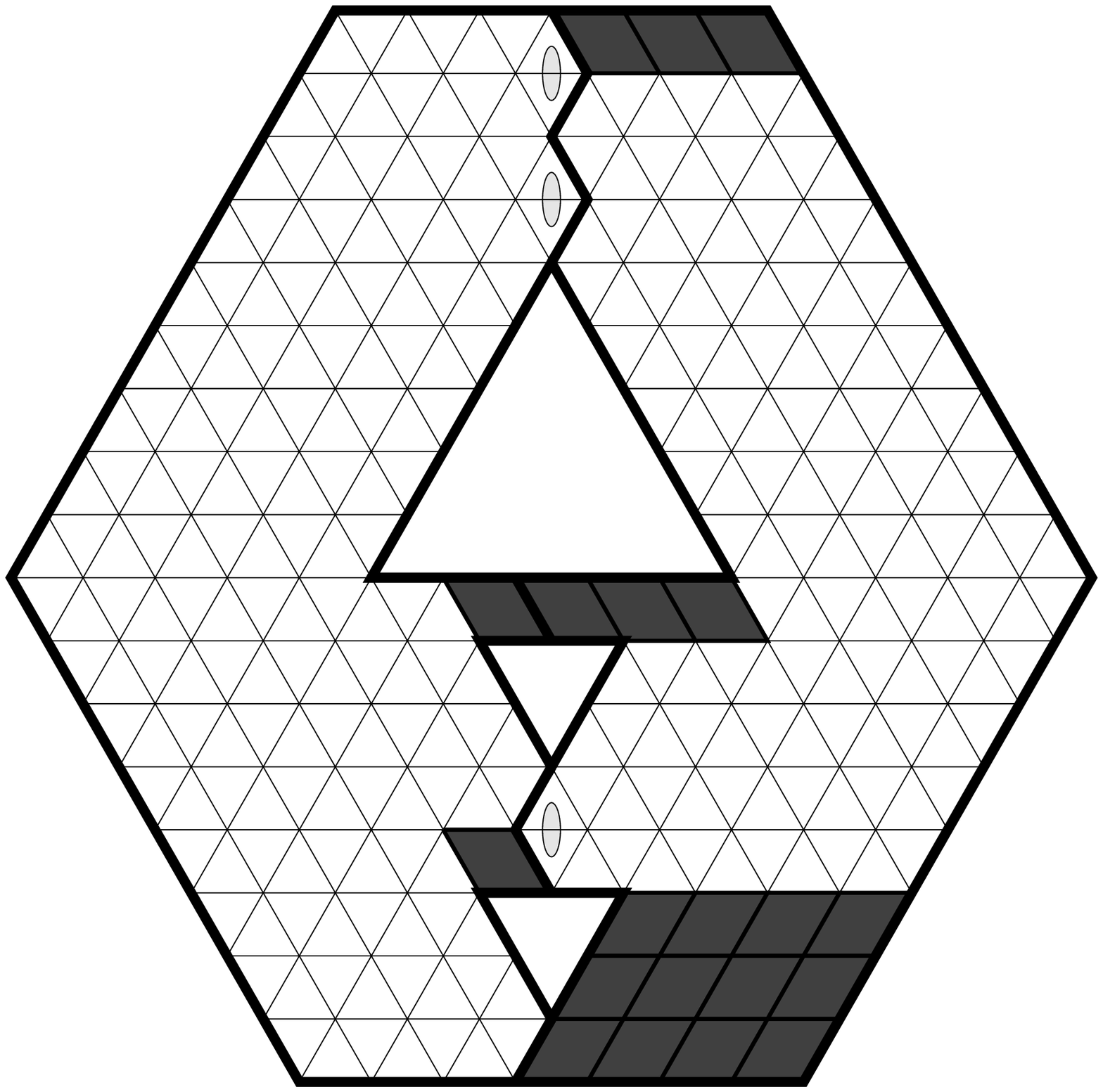}}{\mypic{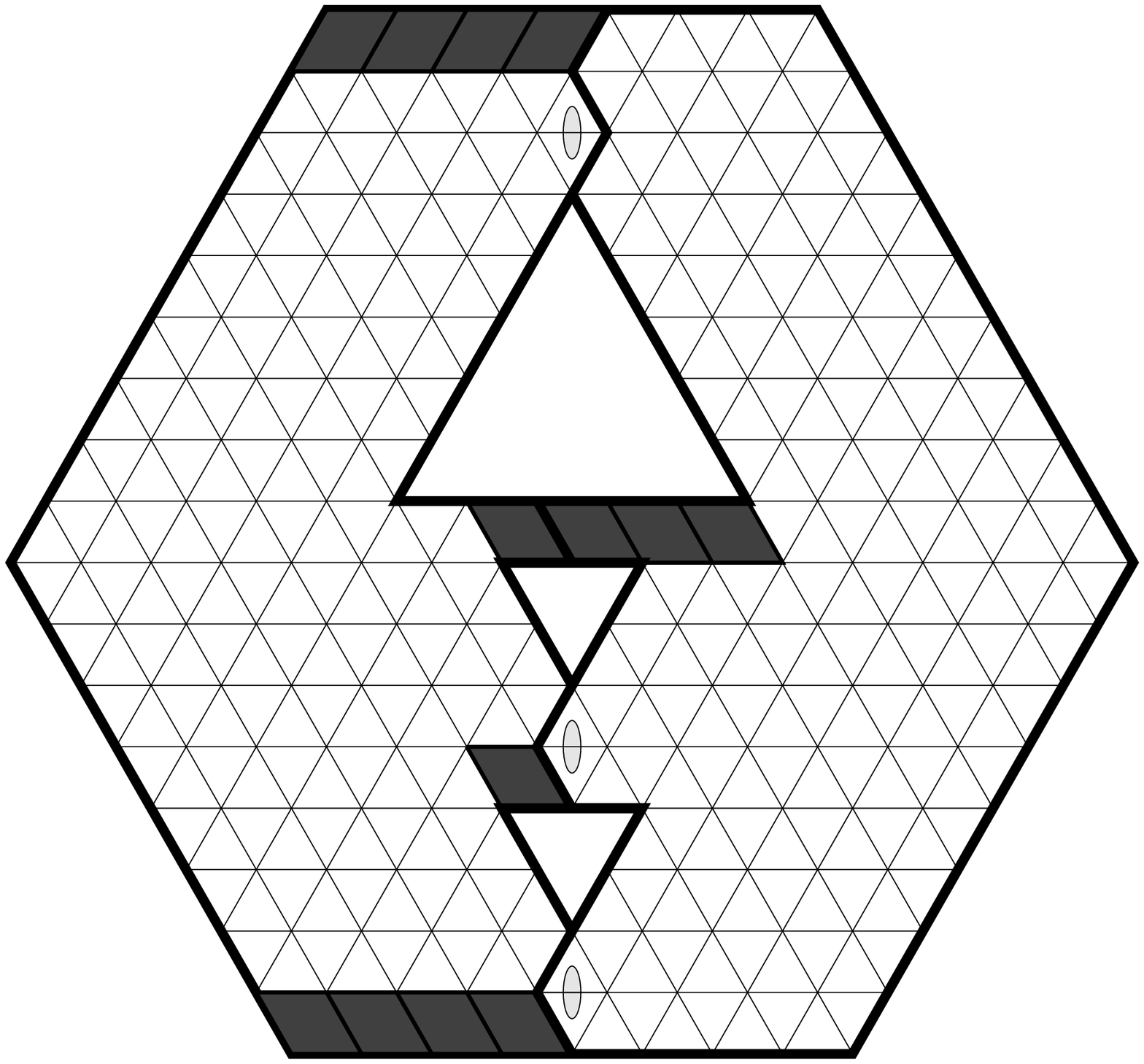}}
\twoline{Figure~3.3{\rm (a).}}{Figure~3.3{\rm (b).}}
\twoline{{\rm $H_{(2,4),(3,4)}(6,8,1)$ reduces to}}{{\rm $H_{(2,4),(3,4)}(7,8,1)$ reduces to}}
\twoline{{\rm $\bar{R}_{(2,4),(3,4)}(3)$ and $R_{(3,4),(2)}(3)$.}}
{{\rm $\bar{R}_{(2,4),(3)}(4)$ and $R_{(3,4),(2,4)}(3)$.}}
\endinsert

\topinsert
\twoline{\mypic{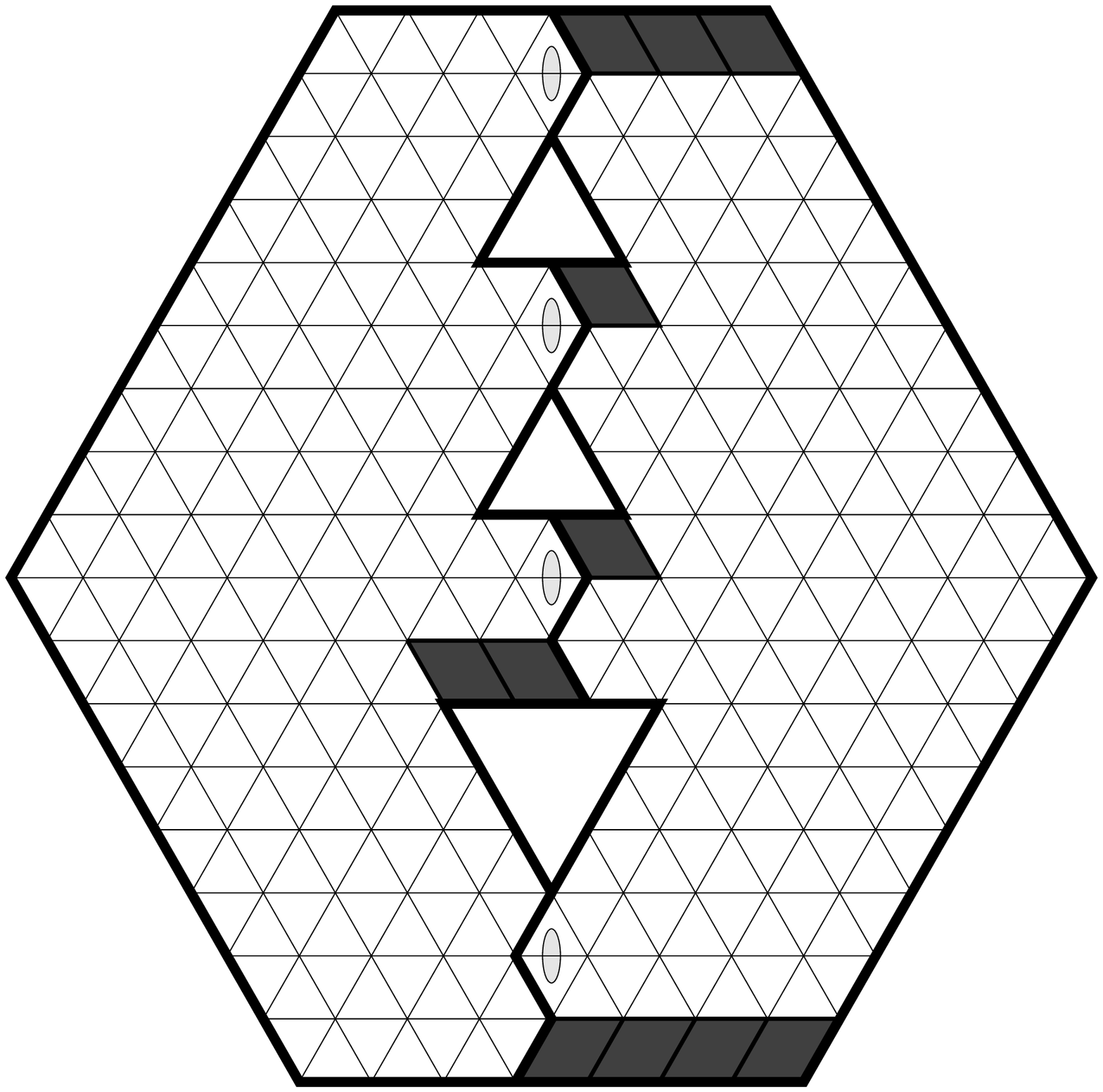}}{\mypic{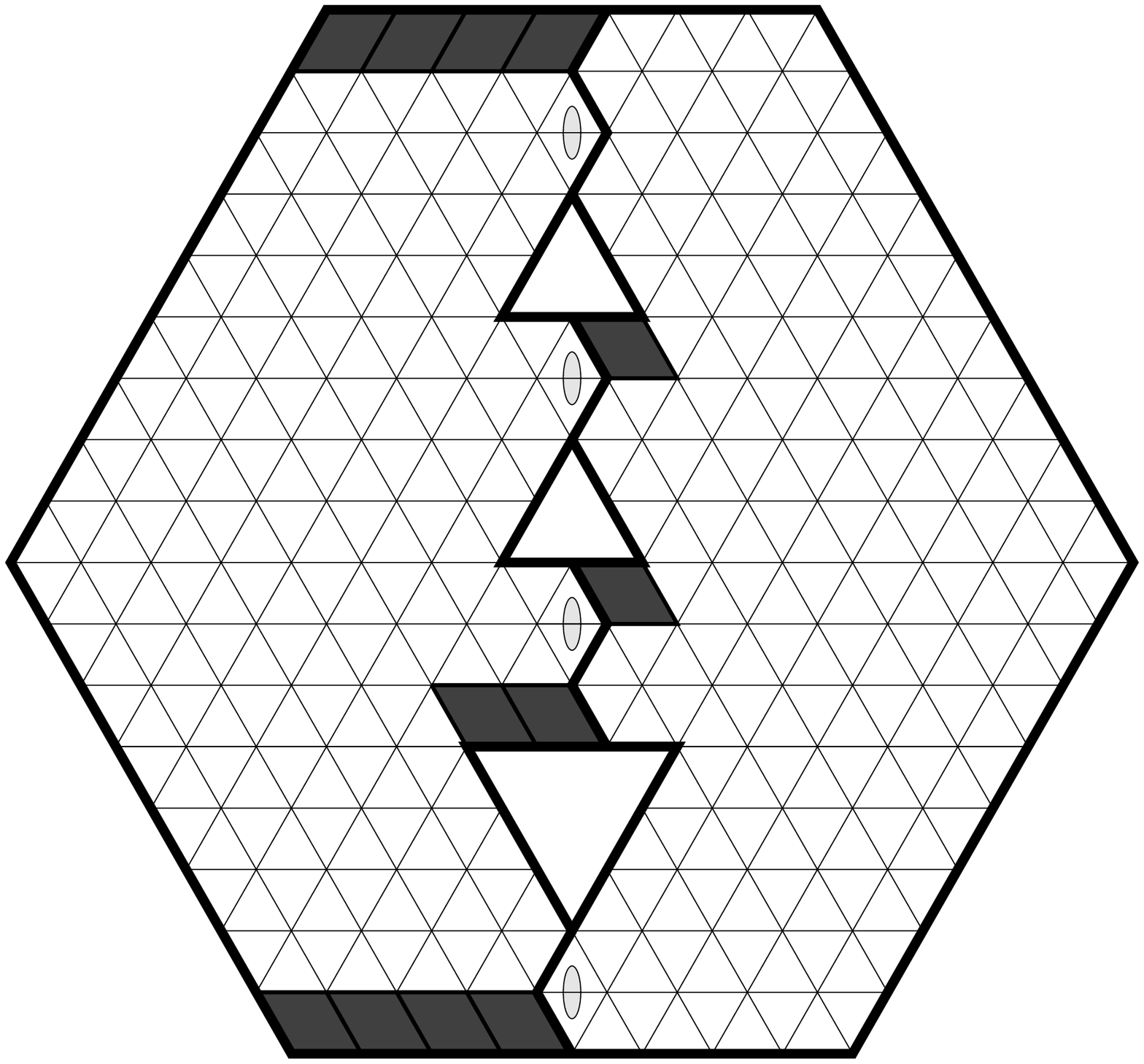}}
\twoline{Figure~3.4{\rm (a).}}{Figure~3.4{\rm (b).}}
\twoline{{\rm $\bar{H}_{(2,3),(1,3,5)}(6,8,1)$ reduces to}}
{{\rm $\bar{H}_{(2),(1,3,5,6)}(7,8,1)$ reduces to}}
\twoline{{\rm $R_{(2,3),(1,3,5)}(3)$ and $\bar{R}_{(1,3,5),(2)}(3)$.}}
{{\rm $R_{(2),(1,3,5)}(4)$ and $\bar{R}_{(1,3,5,6),(2)}(3)$.}}
\endinsert

Recall that our reference line for labeling the vertebrae of $H$ is
now the base of the odd removed window; denote it by $\Cal R$. Suppose $a$ is even. 
Then the last vertebra below
$\Cal R$ is in this case a triangular vertebra, which remains present for all choices of 
removed windows (see Figure 3.3(a)). The rest of the unit triangles of $R$ crossed by 
$\ell$ can be paired to form rhombic vertebrae. By definition, the region $R$ contains 
precisely the vertebrae
labeled $l_1,\dotsc,l_m$ below $\Cal R$, and those labeled $q_1,\dotsc,q_n$ above $\Cal R$. 
It is easy to check that the region obtained from $R^+$ after removing the forced lozenges
is precisely $\bar{R}_{{\bold l},{\bold q}}((a+k-1)/2)$, while the one obtained from $R^-$
after the same procedure is congruent to $R_{{\bold q},{\bold l}^{(m)}}(a/2)$. Since the
width of $R$ is $m+n$, we obtain by Proposition 2.1 the first equality of (1.8). 

If $a$ is odd, a similar analysis shows that $R^+$ is congruent to 
$\bar{R}_{{\bold l},{\bold q}^{(n)}}((a+k)/2)$, and $R^-$ to 
$R_{{\bold q},{\bold l}}((a-1)/2)$ (see Figure 3.3(b) for an illustration). Since the width 
of $R$ is again $m+n$, Proposition 2.1 implies the second equality of (1.8).

The regions $\bar{H}_{{\bold l},{\bold q}}(a,b,k)$ are treated in a perfectly similar way
(this is illustrated in Figures 3.4(a) and (b)).
Again, the regions $R^+$ and $R^-$ obtained by applying (2.1) turn out to be  
members of the families of regions defined in Section 2, and applying Proposition 2.1 
one arrives at the formulas (1.9). $\square$

\medskip
\flushpar
{\smc Remark 3.1.} Part (a) of Theorem 1.1 (i.e., the case $k$ even) also follows from 
Lemma~2.2 of \cite{6}.
Indeed, as shown in the next section, it is easy to express the tiling generating functions
of the regions $R^+$ and $R^-$ arising by applying the Factorization Theorem to
$H_{\bold l}(a,b,k)$ as determinants of Gessel-Viennot matrices (see (4.1)). It is also 
not hard to check that
the entries of these matrices are of the form required by  
\cite{6,Lemma 2.2}, which therefore provides explicit product formulas 
for the value of these determinants. 

For $k$ odd on the other hand, the corresponding Gessel-Viennot matrices
have a more complicated structure and their determinants cannot be computed by
\cite{6,Lemma 2.2}.

\mysec{4. Recurrences for $\M(R_{{\bold l},{\bold q}}(x))$ 
and $\M(\bar{R}_{{\bold l},{\bold q}}(x))$}

\medskip
Let ${\bold l}=(l_1,\dotsc,l_m)$ and ${\bold q}=(q_1,\dotsc,q_n)$ be two (possibly empty) 
lists of strictly increasing positive integers, and let $\Cal T$ be a tiling of 
$R_{{\bold l},{\bold q}}(x)$ (Figure 4.1 shows an example for ${\bold l}=(2,3)$, 
${\bold q}=(2,4)$ and $x=3$; a white center on a tile indicates that it is weighted by 
1/2). This tiling divides the southwestern boundary of 
$R_{{\bold l},{\bold q}}(x)$ into $2l_m-m+n+1$ unit segments. Starting from these segments,
by following the tiles of $\Cal T$, one obtains paths of rhombi that end on the 
southwest-facing unit segments on the right boundary of $R_{{\bold l},{\bold q}}(x)$ (these
paths are shaded dark in Figure 4.1). It is
clear that these paths are non-intersecting, and it is not hard to see that they completely
determine $\Cal T$. Since along the paths of rhombi the only possible steps are northeast
and east, it is clear that they can be identified with lattice paths on the square lattice
$\Z^2$ taking steps north and east. 

We will need the following special case of the Gessel-Viennot theorem for 
non-inter\-sec\-ting lattice paths (see \cite{5}). All our lattice paths 
will be on the grid graph $\Z^2$, 
directed so that its edges point east and north. We allow the edges of $\Z^2$ to be
weighted, and define the weight of a lattice path to be the product of the weights on its 
steps. The weight of an $N$-tuple of lattice paths is the product of the individual weights 
of its members. The generating function of a set of $N$-tuples of lattice paths is the sum
of the weights of its elements.

\topinsert
\centerline{\mypic{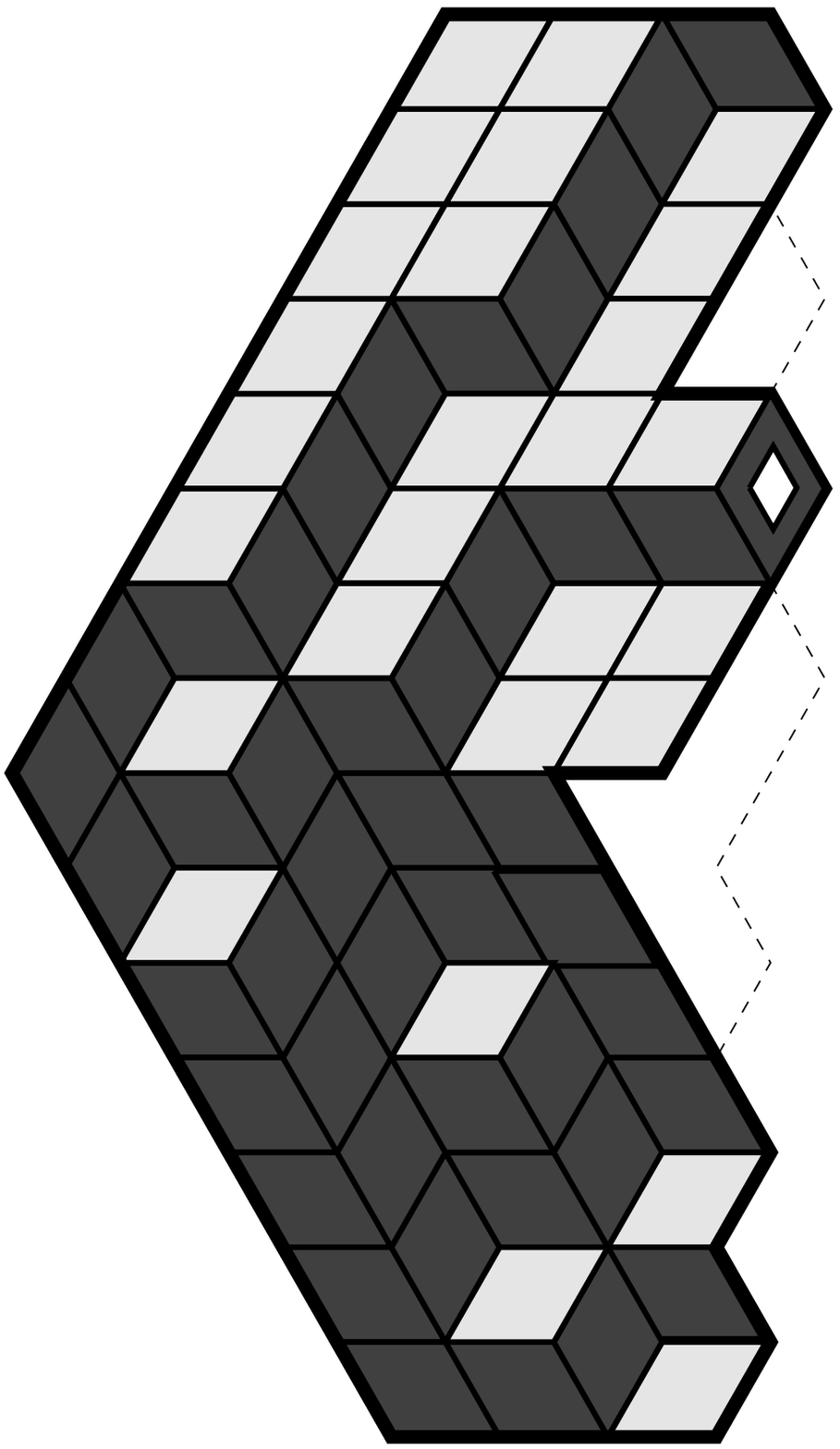}}
\centerline{{\smc Figure~4.1.} {\rm Encoding a tiling as paths of rhombi.}}
\endinsert

An $N$-tuple of starting points ${\bold u}=(u_1,\dotsc,u_N)$ is said to be 
{\it compatible} with an $N$-tuple of ending points ${\bold v}=(v_1,\dotsc,v_N)$  
if for all $i<j$ and $k<l$, every lattice path from $u_i$ to $v_l$ intersects each lattice
path from $u_j$ to $v_k$. 

\proclaim{Theorem 4.1 (Gessel-Viennot)} Let ${\bold u}$ and ${\bold v}$ be compatible 
$N$-tuples of starting 
and ending points on the graph $\Z^2$ oriented as above. Then the generating function for
$N$-tuples of  non-intersecting lattice paths joining the starting points to the ending 
points is 

$$\det\left((a_{ij})_{1\leq i,j\leq N}\right),$$ 
where $a_{ij}$ is the generating function for lattice paths from $u_i$ to $v_j$.
\endproclaim

We will find it convenient to view the paths of rhombi directly as lattice paths on $\Z^2$, 
thus bypassing the ``extra steps'' of bijecting them with lattice paths on a lattice of 
rhombi with 
angles of 60 and 120 degrees, and then deforming this to the square lattice. In this context,
the ``points'' of our lattice $\Cal L$ are the southwest-facing edges of the triangular 
lattice --- we call them {\it segments} ---, and the ``lines''of $\Cal L$  are sequences of 
adjacent rhombi extending 
horizontally or in the southwest-northeast direction --- we call these {\it rows} and 
{\it columns}, respectively. The rhombi that make up these rows and columns are called
{\it edges} of $\Cal L$. 

To coordinatize our lattice, we choose the $x$-axis to be the bottommost row, and the 
$y$-axis to be the leftmost column, intersecting $R_{{\bold l},{\bold q}}(x)$. Let 
$N=2l_m-m+n+1$, and label the starting and ending segments of the lattice paths encoding
our tiling $\Cal T$, from bottom to top, by $u_1,\dotsc,u_N$ and $v_1,\dotsc,v_N$, 
respectively. 

Weight by 1/2 the edges of $\Cal L$ corresponding to tile positions weighted by 1/2 in 
$R_{{\bold l},{\bold q}}(x)$. Weight all other edges of $\Cal L$ inside 
$R_{{\bold l},{\bold q}}(x)$ by 1, and the edges outside $R_{{\bold l},{\bold q}}(x)$ by 0.
By Theorem 4.1 and the above-mentioned bijection between tilings and lattice paths, we
obtain that 

$$\M(R_{{\bold l},{\bold q}}(x))=\det\left((a_{ij})_{1\leq i,j\leq N}\right),\tag4.1$$
where  $a_{ij}$ is the generating function for lattice paths from $u_i$ to $v_j$. 

Suppose $m\leq n$. Then expanding
this determinant along the last row turns out to give us a recurrence relation for
$\M(R_{{\bold l},{\bold q}}(x))$.

\proclaim{Lemma 4.2} For $m\leq n$, we have

$$\M(R_{{\bold l},{\bold q}}(x))=\sum_{k=1}^{n}(-1)^{n-k}C_k
\M(R_{{\bold l},{\bold q}^{(k)}}(x))+
(-1)^n\delta_{mn}\M(\bar{R}_{{\bold l},{\bold q}}(x)),\tag4.2$$
where $\delta_{mn}$ is the Kronecker symbol and

$$C_k={x+l_m+q_k\choose 2q_k+m-n}+\frac12{x+l_m+q_k\choose 2q_k+m-n-1}.\tag4.3$$
\endproclaim

\pf Let $A$ be the matrix in (4.1). Expanding the determinant in (4.1) along the last row 
we obtain

$$\M(R_{{\bold l},{\bold q}}(x))=\sum_{k=0}^{N-1}(-1)^{k}a_{N,N-k}\det(A^{N,N-k}),\tag4.4$$
where $A^{ij}$ is the matrix obtained from $A$ by deleting row $i$ and column $j$.  

We claim that at most the first $n+1$ terms in the sum (4.4) are nonzero. 
Indeed, 
recall that $a_{N,N-k}$ is the generating function for lattice paths on $\Cal L$ going
from the segment $u_N$ to the segment $v_{N-k}$. It follows from our construction of
$R_{{\bold l},{\bold q}}(x)$ that there is a unique segment $s$ on its right boundary 
contained in the row immediately below $O$, and $s$ is facing southeast. 
The $y$-coordinates of the segments 
$u_N$ and $s$ work out to be, in our coordinatization, $2l_m-m+n$ and $2l_m$, 
respectively. Since our paths on $\Cal L$ take steps northeast and east, and $m\leq n$,
we see that a path starting at $u_N$ will end, on the right boundary of 
$R_{{\bold l},{\bold q}}(x)$, at some segment above $s$, or, possibly, in case $m=n$, at 
$s$ itself (see Figures 4.2(a), 4.3(a) and 4.4(a)). Since there are precisely $n$ 
segments above $s$ on the right boundary of 
$R_{{\bold l},{\bold q}}(x)$ (one corresponding to each selected bump on $P_u$), this
proves our claim (and also implies $s=v_{N-n}$).

Furthermore, it is easy to determine the exact value of $a_{N,N-k}$ for $k\leq n$. It 
follows from the previous paragraph that $a_{N,N-n}=0$ unless $m=n$, in which case
$a_{N,N-n}=1$ (this accounts for the Kronecker symbol in (4.2)). For $1\leq k\leq n$,
$a_{N,N-n+k}$ is the generating function of lattice paths going from $u_N$ to $v_{N-n+k}$. 
In our system of coordinates it turns out that we have

$$\align
u_N&=(0,2l_m-m+n)\\
v_{N-n+k}&=(x+l_m-q_k-m+n+1,2q_k+2l_m),\ \ \ \ \text{for $1\leq k\leq n$}.
\endalign$$

Partitioning the lattice paths from $u_N$ to $v_{N-n+k}$ in two classes according to the
direction of their last step we obtain that $a_{N,N-n+k}=C_k$, with $C_k$ given by (4.3).
Replacing the upper summation limit in (4.4) by $n$ and then replacing the summation
index $k$ by $n-k$ we obtain from (4.4) that

$$\M(R_{{\bold l},{\bold q}}(x))=\sum_{k=1}^{n}(-1)^{n-k}C_k
\det(A^{N,N-n+k})+
(-1)^n\delta_{mn}\det(A^{N,N-n}).\tag4.5$$

\topinsert
\twoline{\mypic{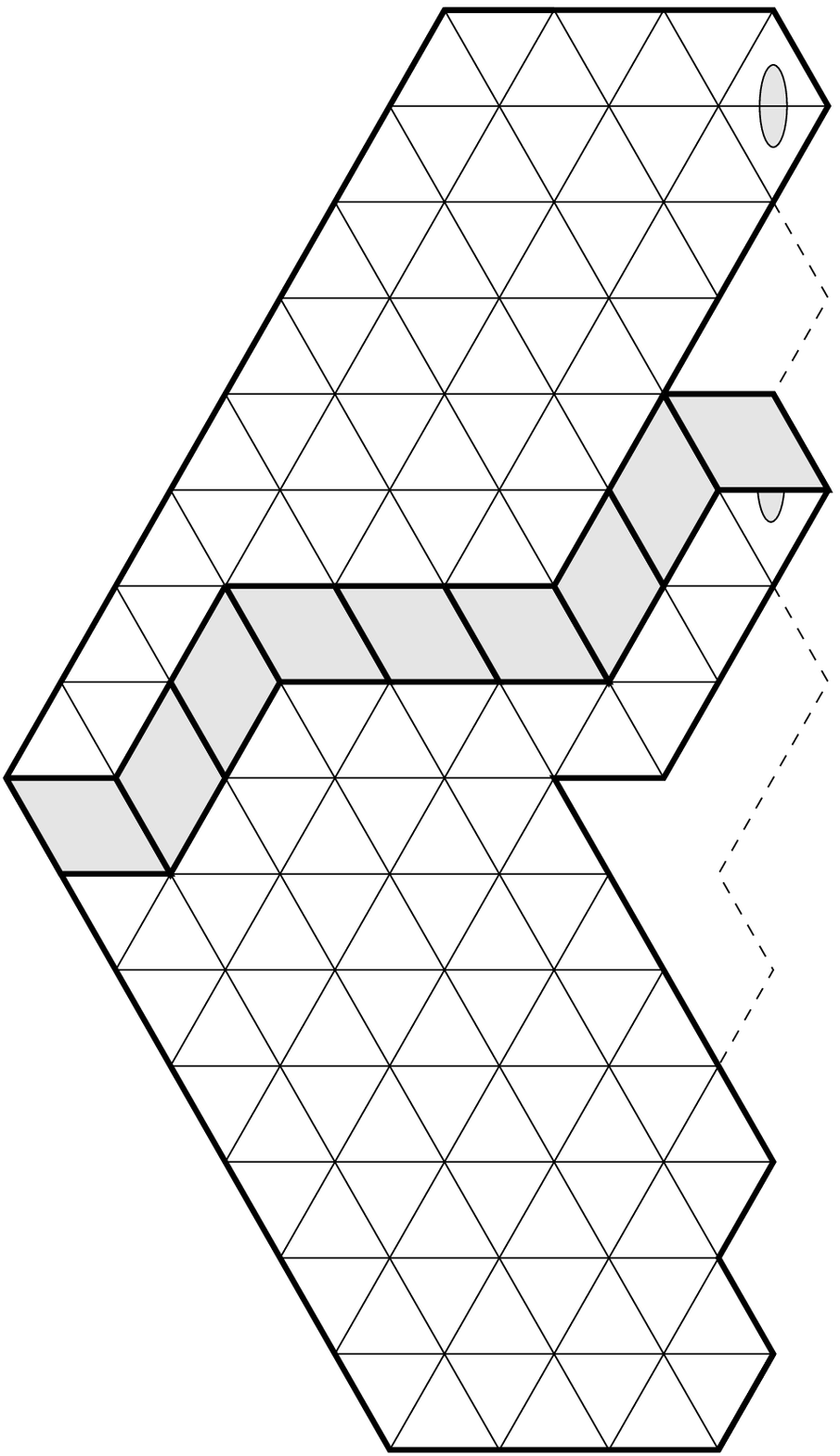}}{\mypic{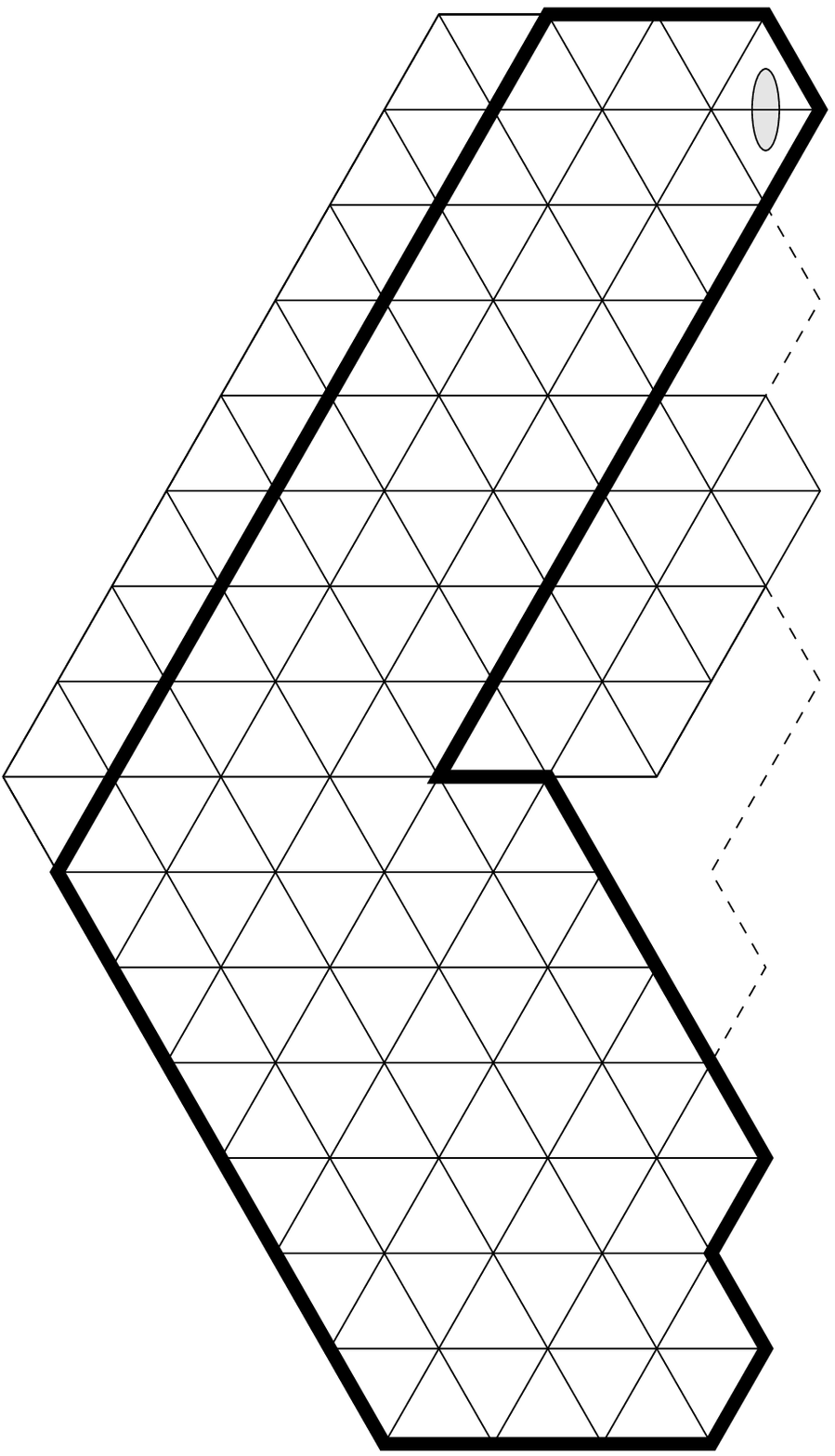}}
\twoline{Figure~4.2{\rm (a).}}{Figure~4.2{\rm (b).}}
\endinsert

\topinsert
\twoline{\mypic{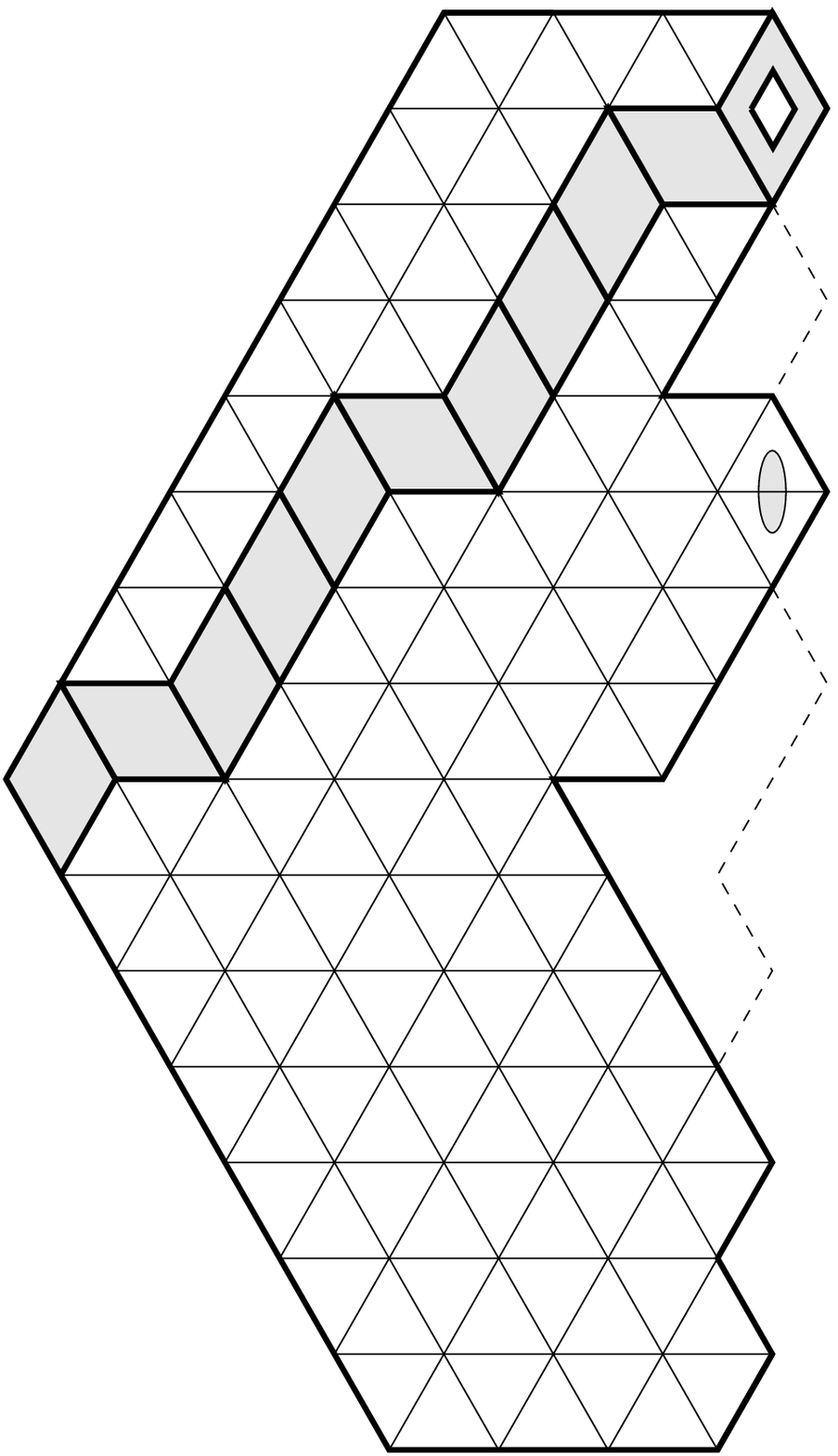}}{\mypic{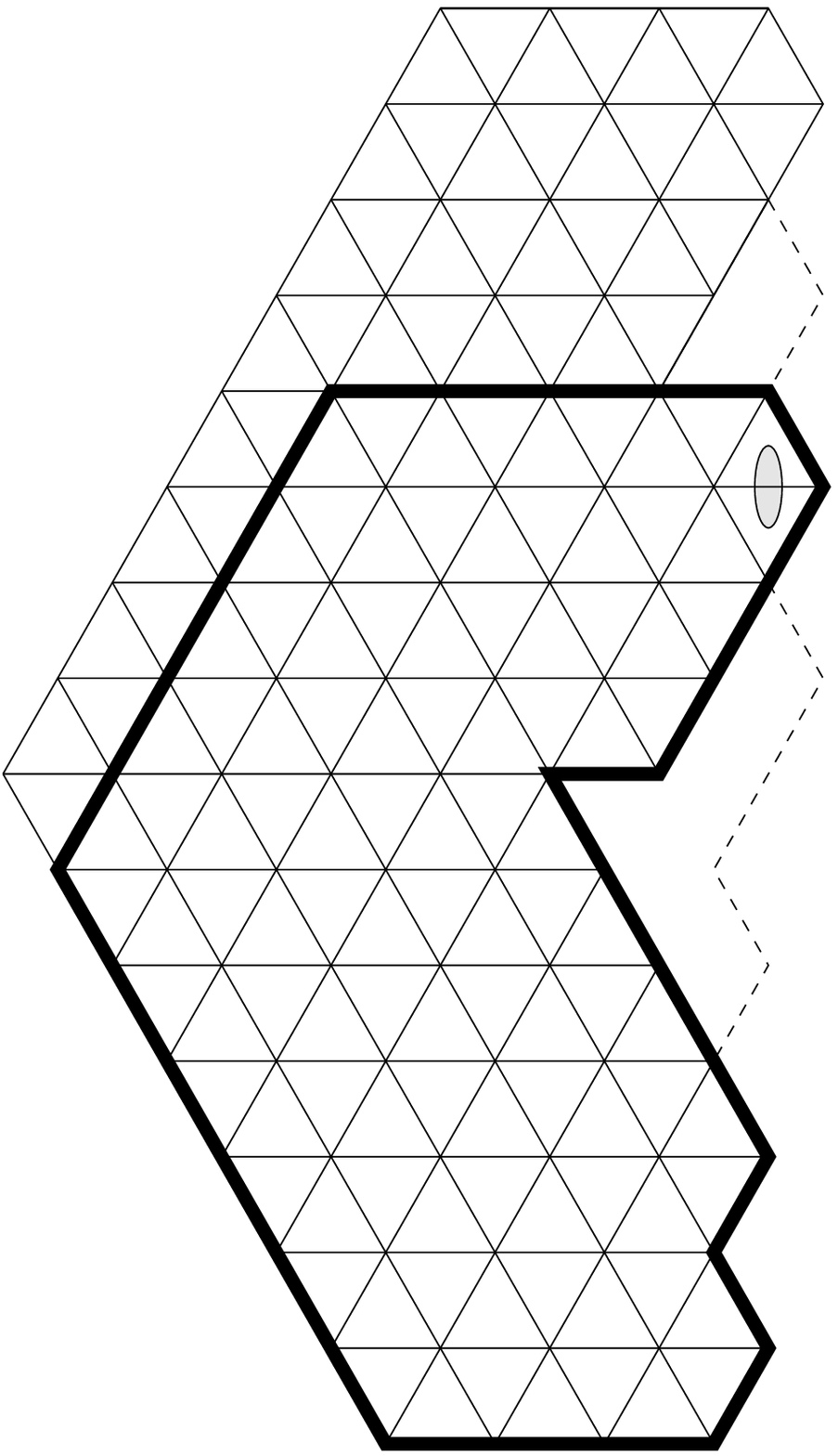}}
\twoline{Figure~4.3{\rm (a).}}{Figure~4.3{\rm (b).}}
\endinsert

\topinsert
\twoline{\mypic{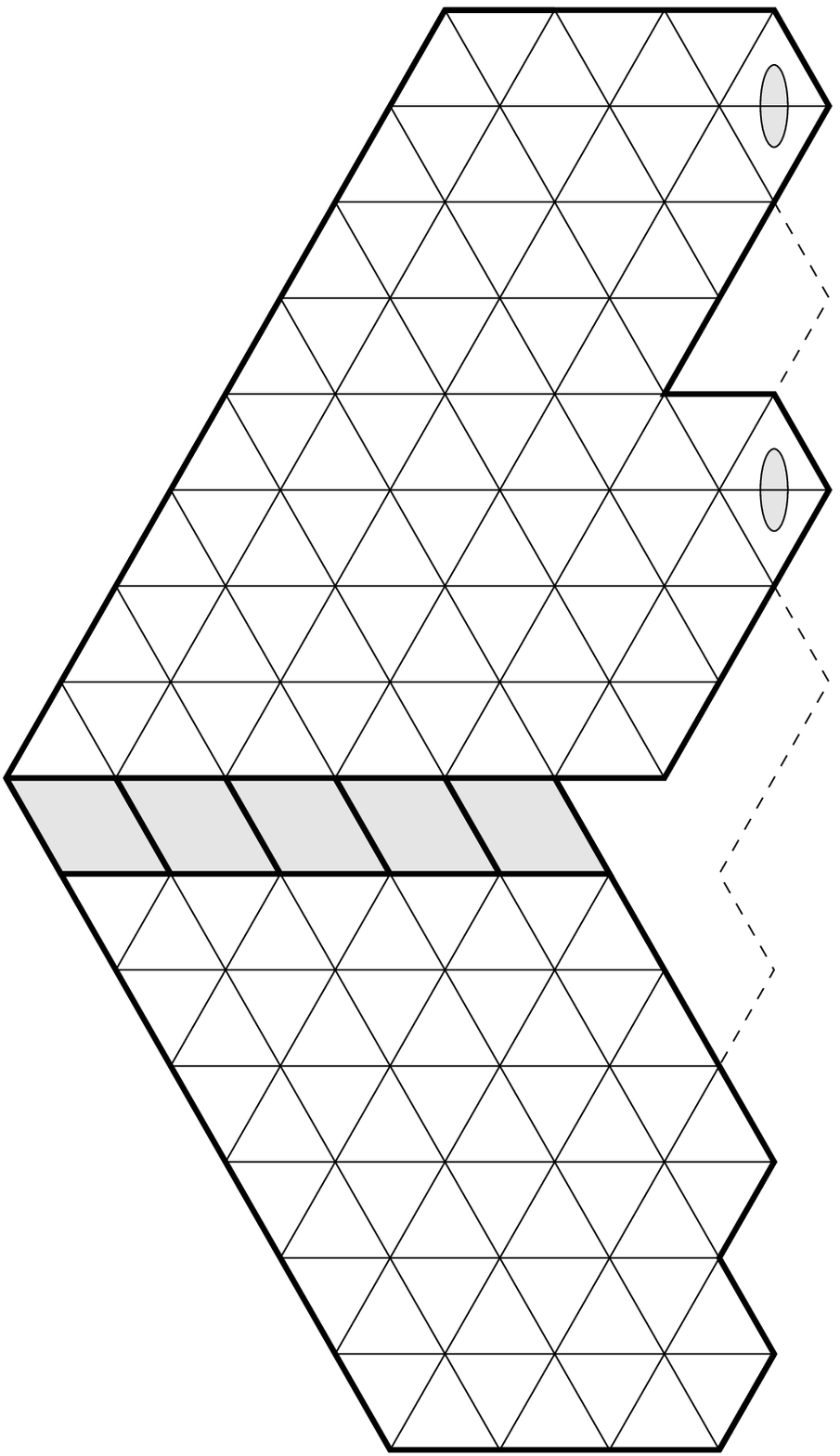}}{\mypic{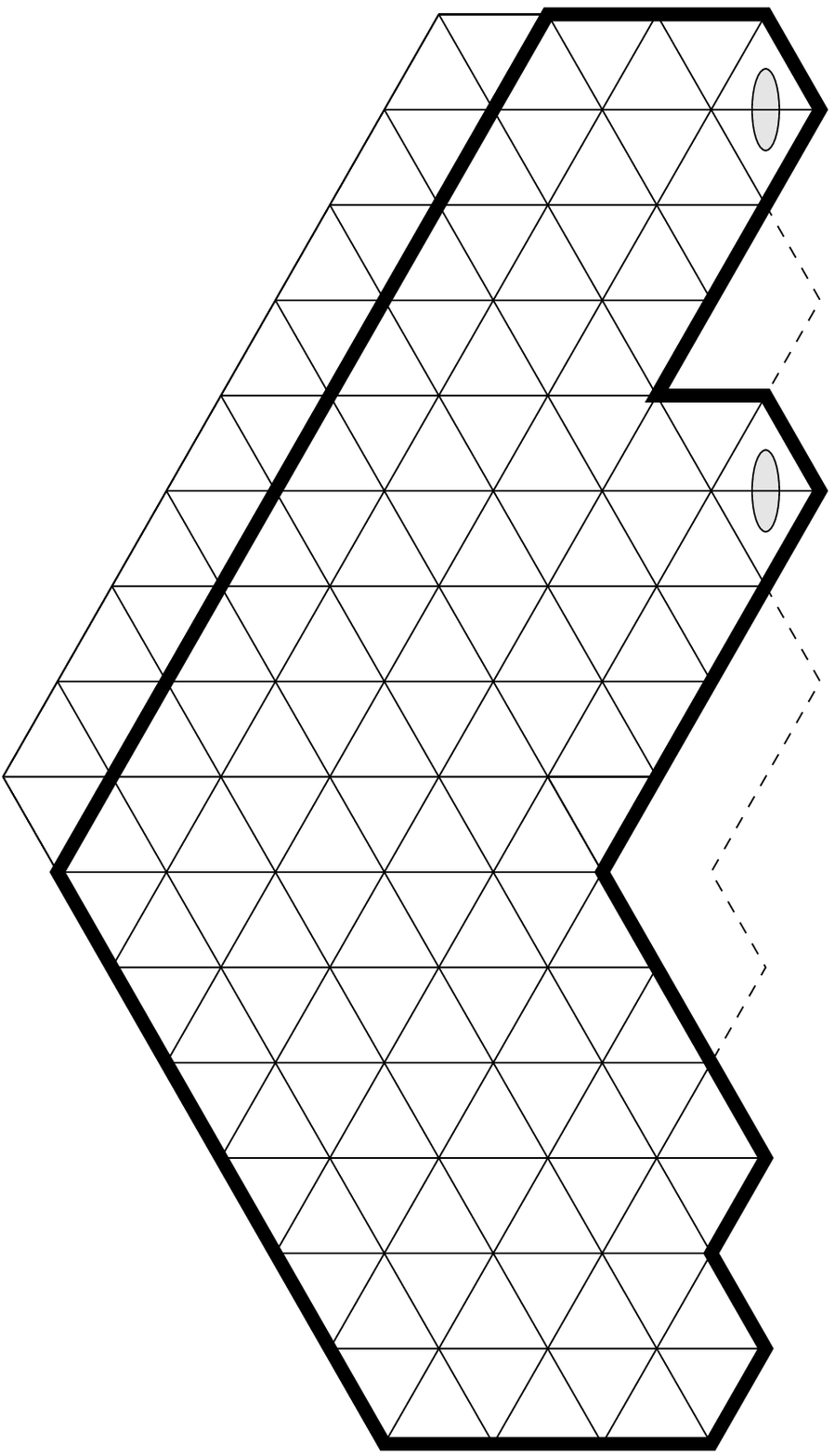}}
\twoline{Figure~4.4{\rm (a).}}{Figure~4.4{\rm (b).}}
\endinsert

To complete the proof of the Lemma, we show that the determinants on the right hand side 
of (4.5) are equal to the corresponding tiling generating functions on the right hand side 
of (4.2).

More precisely, we claim that, for $1\leq k\leq n$, we have
$\det(A^{N,N-n+k})=\M(R_{{\bold l},{\bold q}^{(k)}}(x))$. Indeed, by Theorem 4.1 it
follows that 
$\det(A^{N,N-n+k})$ is equal to the generating function of $(N-1)$-tuples of 
non-intersecting lattice paths starting at $u_1,\dotsc,u_{N-1}$ and ending at 
$v_1,\dotsc,v_{N-n+k-1},v_{N-n+k+1},\dotsc,v_N$. However, by the bijection between
tilings and lattice paths, this is precisely the generating function for tilings of
$R_{{\bold l},{\bold q}^{(k)}}(x)$ (for ${\bold l}=(2,3)$, ${\bold q}=(2,4)$, this is 
illustrated in Figures 4.2 and 4.3; the starting and ending segments of the paths on 
the left are omitted, leading to the regions on the right). 

Similarly, one sees that $\det(A^{N,N-n})=\M(\bar{R}_{{\bold l},{\bold q}}(x))$. Indeed,
$\det(A^{N,N-n})$ is the generating function for $(N-1)$-tuples of 
non-intersecting lattice paths going from $u_1,\dotsc,u_{N-1}$ to 
$v_1,\dotsc,v_{N-n-1},v_{N-n+1},\dotsc,v_N$, respectively. In this case the corresponding
region has the same sets of selected bumps as $R_{{\bold l},{\bold q}}(x)$, but its
boundary around $O$ is changed so that it becomes precisely 
$\bar{R}_{{\bold l},{\bold q}}(x))$ (see Figures 4.4(a) and (b)). This completes the proof 
of the Lemma. $\square$

\medskip
The above result gives a recurrence relation for $\M(R_{{\bold l},{\bold q}}(x))$ when
$m\leq n$. Note that this assumption was indeed necessary: for $m>n$ there are additional 
non-zero terms in the sum (4.4), and the regions they can be associated to are {\it not} 
part of our families $R$ and $\bar{R}$. To deal with the case $m>n$ we use the ``symmetry''
of $R_{{\bold l},{\bold q}}(x)$: we encode its tilings by non-intersecting lattice paths
starting on its northwestern side and taking steps southeast and east.
This leads us to the following result.

\proclaim{Lemma 4.3} For $m>n$, we have

$$\M(R_{{\bold l},{\bold q}}(x))=\sum_{k=1}^{m-1}(-1)^{m-k}D_k
\M(R_{{\bold l}^{(k)},{\bold q}}(x-1))+
D_m\M(R_{{\bold l}^{(m)},{\bold q}}(x+l_m-l_{m-1}-1)),\tag4.6$$
where

$$D_k={x+l_m+l_k-m+n+1\choose 2l_k-m+n+1}.\tag4.7$$
\endproclaim

\pf We proceed in the same fashion as in the proof of Lemma 4.2. First, we change the 
definition of the lattice $\Cal L$ so that its points (called again segments) are the
unit edges of the triangular lattice facing northwest, and the lines (still called rows 
and columns) are sequences of adjacent unit rhombi extending horizontally and in the 
northwest-southeast direction. We choose the $x$- and $y$-axis 
to be the topmost row and the
leftmost column intersecting $R_{{\bold l},{\bold q}}(x)$, respectively. Since the
northwestern side of $R_{{\bold l},{\bold q}}(x)$ has length $2q_n+m-n$, a tiling of
$R_{{\bold l},{\bold q}}(x)$ is encoded by this many non-intersecting lattice paths in 
$\Cal L$. Let $N=2q_n+m-n$ and label the starting and ending segments of these lattice
paths, from top to bottom, by $u_1,\dotsc,u_N$ and $v_1,\dotsc,v_N$, respectively.

\topinsert
\twoline{\mypic{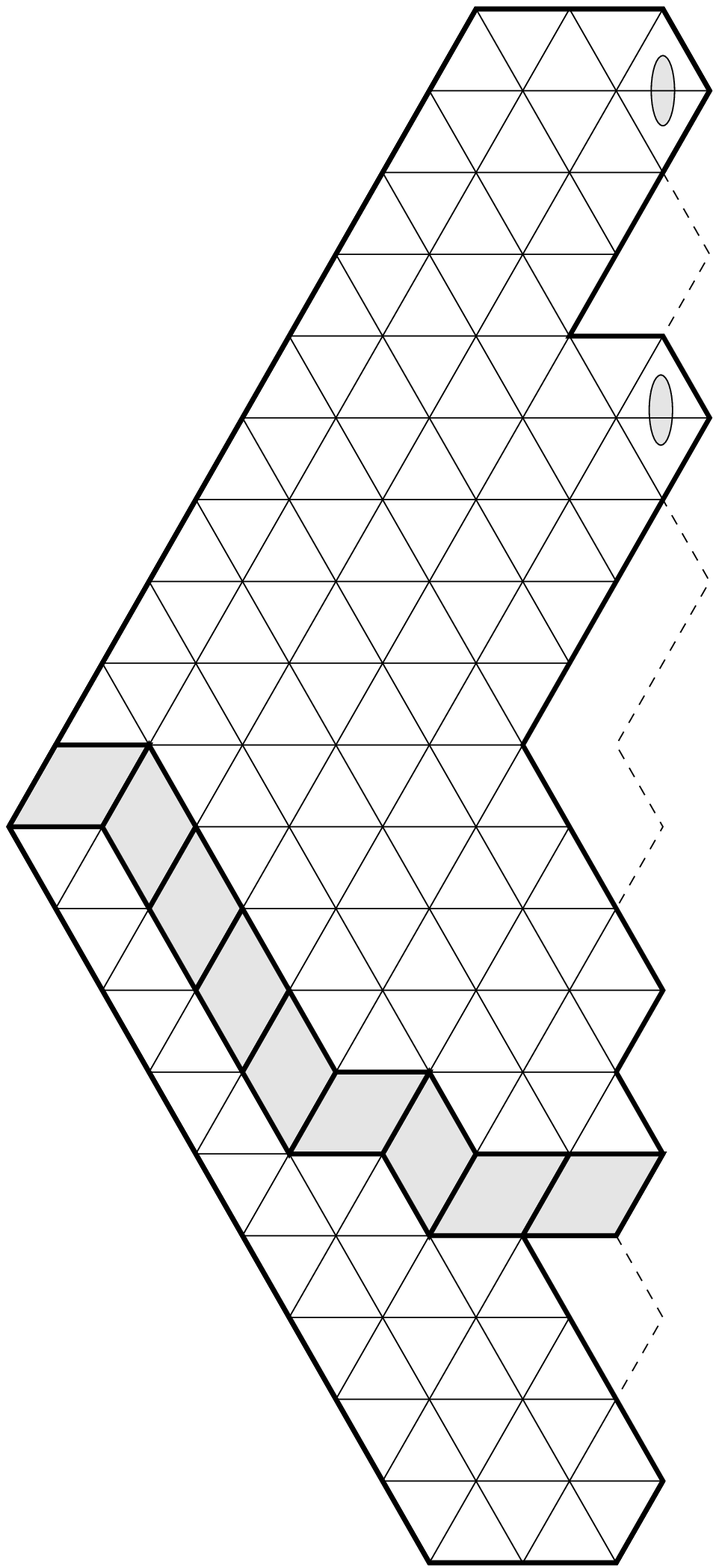}}{\mypic{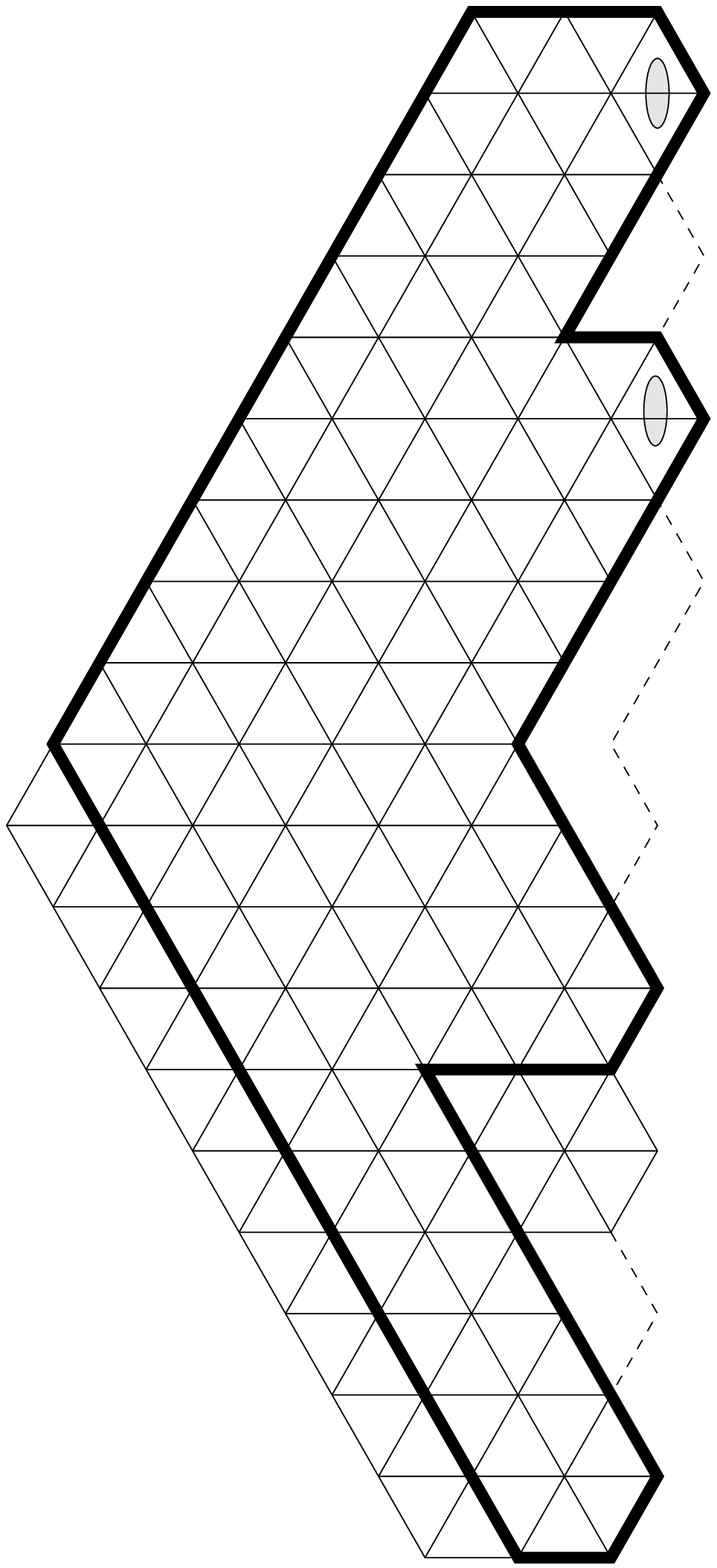}}
\twoline{Figure~4.5{\rm (a).}}{Figure~4.5{\rm (b).}}
\endinsert

\topinsert
\twoline{\mypic{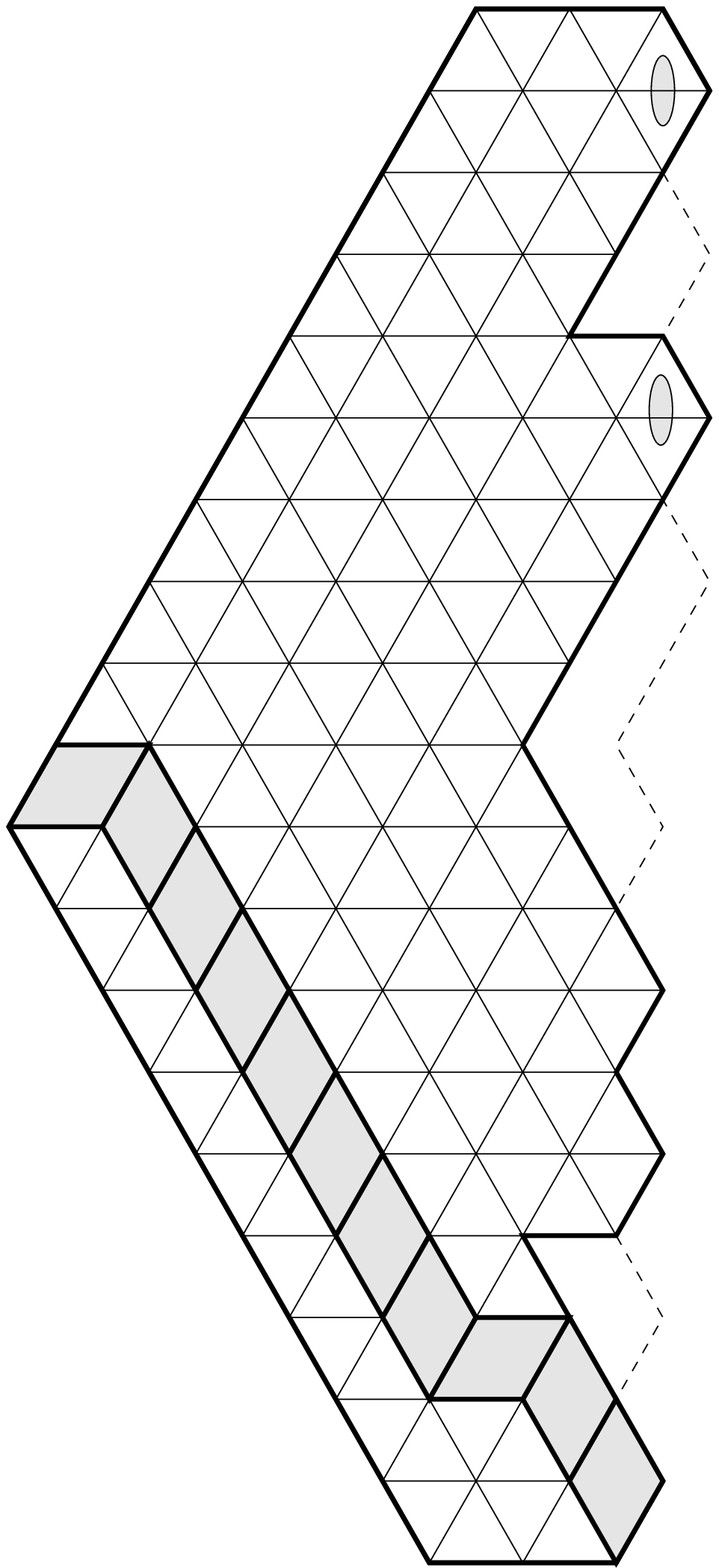}}{\mypic{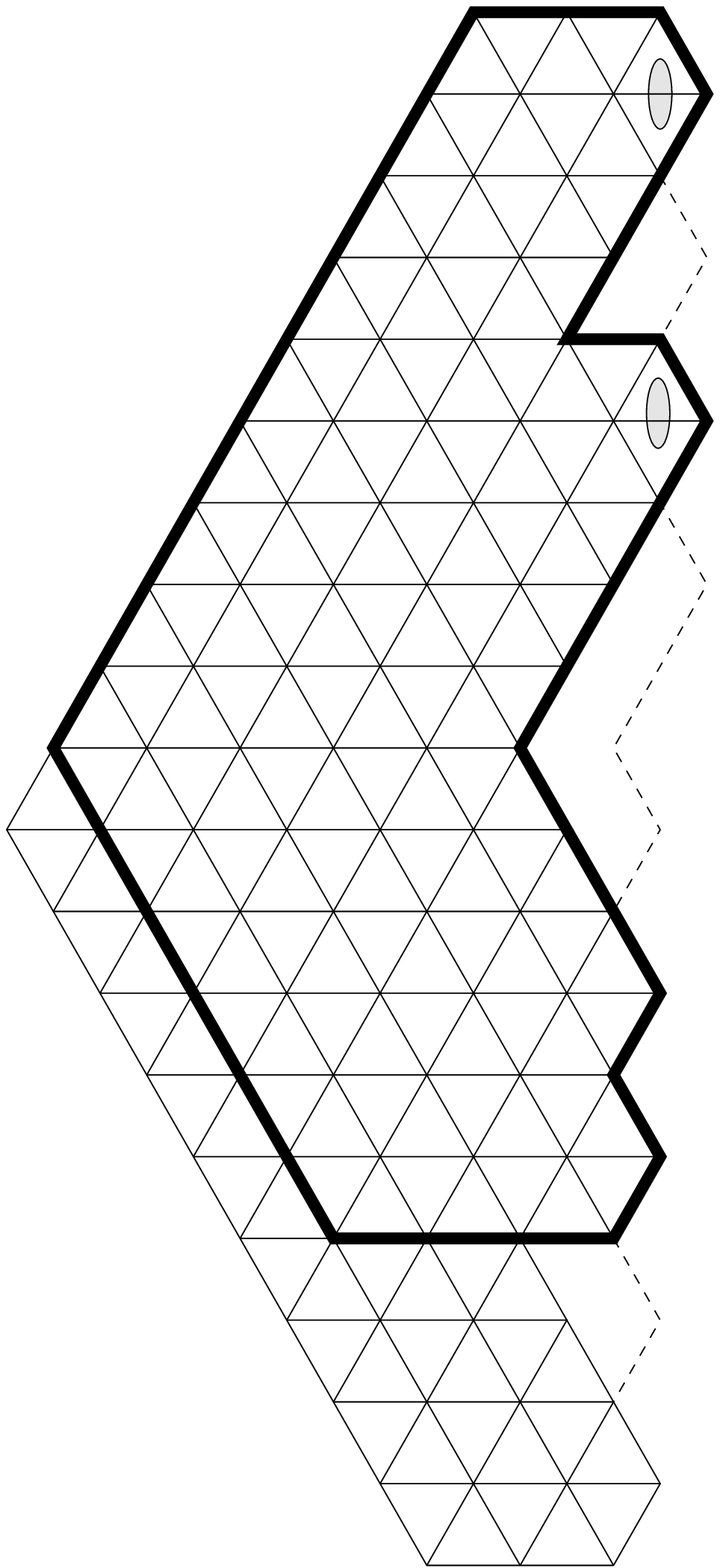}}
\twoline{Figure~4.6{\rm (a).}}{Figure~4.6{\rm (b).}}
\endinsert

Formula (4.1) holds also with this new definitions of $\Cal L$, $u_i$ and $v_i$. 
Expanding the determinant on its right hand side along the last row we obtain

$$\M(R_{{\bold l},{\bold q}}(x))=\sum_{k=0}^{N-1}(-1)^{k}a_{N,N-k}\det(A^{N,N-k}),\tag4.8$$
where the entries of $A$ are now defined in terms of the new $u_i$'s and $v_i$'s.

We claim that all terms in this sum except for the first $m$ are 0. Indeed, denote by $s$
the unique segment in the right boundary of $R_{{\bold l},{\bold q}}(x)$ contained in
the row immediately above $O$. The $y$-coordinates of $u_N$ and $s$ in $\Cal L$ are
$2q_n+m-n-1$ and $2q_n-1$, respectively. Since $m>n$, and our paths take steps southeast 
and east, it follows that the path starting at $u_N$ ends at a segment below $s$ on the
right boundary of $R_{{\bold l},{\bold q}}(x)$ (this is illustrated in Figures 4.5(a)
and 4.6(a) for ${\bold l}=(2,3,5)$, ${\bold q}=(2,4)$ and $x=2$). 
There are exactly $m$ such segments, one
for each selected bump on $P_d$. So $a_{Nk}=0$ unless $k\in\{N-m+1,\dotsc,N\}$, which
proves our claim.

For $1\leq k\leq m$, the exact value of $a_{N,N-m+k}$ follows from the coordinates
of $u_N$ and $v_{N-m+k}$. Since we have

$$\align
u_N&=(0,2q_n+m-n-1)\\
v_{N-m+k}&=(x+l_m-l_k,2l_k+2q_n),\ \ \ \ \text{for $1\leq k\leq m$},
\endalign$$
we obtain that $a_{N,N-m+k}=D_k$, where $D_k$ is given by (4.7). Therefore, replacing the
upper summation limit in (4.8) by $m-1$ and then replacing the summation index $k$ by
$m-k$, (4.8) becomes

$$\M(R_{{\bold l},{\bold q}}(x))=\sum_{k=1}^{m-1}(-1)^{m-k}D_k
\det(A^{N,N-m+k})+
D_m\det(A^{NN}).\tag4.9$$

To complete the proof of the Lemma, we show that the determinants on the right hand side 
of (4.9) are equal to the corresponding tiling generating functions on the right hand
side of (4.6). 

We claim that, for $1\leq k\leq m-1$, we have 
$\det(A^{N,N-m+k})=\M(R_{{\bold l}^{(k)},{\bold q}}(x-1))$. Indeed, by Theorem 4.1, 
$\det(A^{N,N-m+k})$ is equal to the tiling generating function for $(N-1)$-tuples of
non-intersecting lattice paths in $\Cal L$ starting at $u_1,\dotsc,u_{N-1}$ and ending at 
$v_1,\dotsc,v_{N-m+k-1},v_{N-m+k+1},\dotsc,v_N$. By the bijection between
tilings and lattice paths, this is precisely the generating function for tilings of
$R_{{\bold l}^{(k)},{\bold q}}(x-1)$. (The change in the parameter $x$ is a new ingredient
that was not present in the case $m\leq n$; this is due to the change of our 
coordinate system in the present case. Figure 4.5 illustrates this for $k=2$.)  

Similarly, one sees that 
$\det(A^{N,N})=\M(R_{{\bold l}^{(m)},{\bold q}}(x+l_m-l_{m-1}-1))$. Indeed,
$\det(A^{N,N})$ is the generating function for $(N-1)$-tuples of 
non-intersecting lattice paths going from $u_1,\dotsc,u_{N-1}$ to 
$v_1,\dotsc,v_{N-1}$, respectively. The corresponding region is the portion of
$R_{{\bold l},{\bold q}}(x)$ consisting of the unit triangles in and above the row 
containing $v_{N-1}$, and to the right of the column containing $u_N$ (see Figure 4.6). 
This is readily seen to be exactly 
$R_{{\bold l}^{(m)},{\bold q}}(x+l_m-l_{m-1}-1))$, thus completing the proof. $\square$

\medskip
We will also need the analogs of the above two Lemmas for the regions
$\bar{R}_{{\bold l},{\bold q}}(x)$. 

\proclaim{Lemma 4.4} For $m<n$, we have

$$\M(\bar{R}_{{\bold l},{\bold q}}(x))=\sum_{k=1}^{n}(-1)^{n-k}\bar{C}_k
\M(\bar{R}_{{\bold l},{\bold q}^{(k)}}(x)),\tag4.10$$
where

$$\bar{C}_k={x+l_m+q_k\choose 2q_k+m-n+1}+\frac12{x+l_m+q_k\choose 2q_k+m-n}.\tag4.11$$
\endproclaim

\pf Let $\Cal L$ be as in the proof of Lemma 4.2. The southwestern side of
$\bar{R}_{{\bold l},{\bold q}}(x)$ has length $N=2l_m-m+n$. Encode the tilings of
$\bar{R}_{{\bold l},{\bold q}}(x)$ as $N$-tuples of non-intersecting lattice paths on
$\Cal L$. Order their starting and ending points from bottom to top and express their 
number as a determinant by Theorem 4.1. Using the same 
reasoning as in the proof of Lemma 4.2, it follows that only the last $n$ entries in 
the last row
of this determinant are nonzero. These are easily seen to be equal to the numbers 
$\bar{C}_k$ given by (4.11). Expand this determinant along the last row. Using arguments
similar to those in the proof of Lemma 4.2, the order $N-1$ minors in the determinant 
expansion can be interpreted as tiling generating functions of regions from the 
$\bar{R}$-family, and one obtains (4.10). $\square$

\proclaim{Lemma 4.5} For $m\geq n$, we have

$$\align
\M(\bar{R}_{{\bold l},{\bold q}}(x))=\sum_{k=1}^{m-1}(-1)^{m-k}\bar{D}_k
\M(\bar{R}_{{\bold l}^{(k)},{\bold q}}(x-1&))+
\bar{D}_m\M(\bar{R}_{{\bold l}^{(m)},{\bold q}}(x+l_m-l_{m-1}-1))\\
+(-1)^m\delta_{mn}\M(R_{{\bold l},{\bold q}}(x-1&)),\tag4.12
\endalign$$
where $\delta_{mn}$ is the Kronecker symbol and

$$\bar{D}_k={x+l_m+l_k-m+n\choose 2l_k-m+n}.\tag4.13$$
\endproclaim

\pf Choose $\Cal L$ to be as in the proof of Lemma 4.3. The northwestern side of
$\bar{R}_{{\bold l},{\bold q}}(x)$ has length $N=2q_n+m-n+1$. Encode the tilings of
$\bar{R}_{{\bold l},{\bold q}}(x)$ as $N$-tuples of non-intersecting lattice paths on
$\Cal L$. Order their starting and ending points from bottom to top and express their 
number as a determinant by Theorem 4.1. By arguments similar 
to those in the proof of Lemma 4.3, it follows that only the last
$m+\delta_{mn}$ terms in the last row of this determinant are nonzero. More precisely,
it turns out that the $(m+1)$-st to last entry in the last row is $\delta_{mn}$, and the 
subsequent entries are given by (4.13). Expand the
determinant along the last row. In analogy to the situation in the proof of Lemma 4.3, it 
turns out that the minors of order $N-1$ in this expansion are equal to tiling 
generating functions of regions belonging to the $\bar{R}$- and $R$-families. 
Precise identification
of these regions leads to (4.12). (Note that, except for the last term on the right hand 
side of (4.12), the shift in the parameter $x$ is precisely the same as in the case of 
Lemma 4.3. The last term does not have a correspondent in (4.6) because there we didn't
allow $m=n$.) $\square$

\medskip
Finally, to successfully carry out the inductive evaluation of our tiling generating
functions, we need to pay attention to their behaviour when ${\bold l}$,
${\bold q}$ and $x$ are on the ``boundary'' (i.e., ${\bold l}$ or ${\bold q}$ is empty,
or $x$ takes on the minimum value for which our regions are defined). By definition, we take
$\M(R_{\emptyset,\emptyset}(x))=1$ and $\M(\bar{R}_{\emptyset,\emptyset}(x))=1$, for 
all $x$ (empty regions have precisely one tiling). 

In the recurrences (4.2), (4.6), (4.10) and (4.12), the value of the argument in the
terms on the right hand side is at least as large as one less the value of the
argument on the left hand side. Therefore, if at least one of  ${\bold l}$ and ${\bold q}$
is nonempty, these recurrences can be applied as long as the argument $x$ on their
left hand side is not minimum possible --- i.e., for (4.2) and (4.6), if $x$ is strictly 
larger than $\max\{0,q_n-l_m-n+m-1\}$ ($\max\{-1,q_n-n-1\}=q_n-n-1$, in case 
${\bold l}=\emptyset$), and, for (4.10) and (4.12), if $x$ is strictly 
larger than $\max\{0,q_n-l_m-n+m\}$.

When $x$ does take on its minimum possible value, we can express the tiling generating 
functions of our regions as shown in the following two Lemmas.

\proclaim{Lemma 4.6} $(a)$ Suppose $l_m-m+1\geq q_n-n$. Then for ${\bold l}\ne\emptyset$
we have

$$\M(R_{{\bold l},{\bold q}}(0))=\M(R_{{\bold l}^{(m)},{\bold q}}(l_m-l_{m-1}-1)).\tag4.14$$

$(b)$ Let $l_m-m+1\leq q_n-n$. Then we have

$$\M(R_{{\bold l},{\bold q}}(q_n-l_m-n+m-1))=
\frac12 \M(R_{{\bold l},{\bold q}^{(n)}}(q_n-l_m-n+m-1)).\tag4.15$$
\endproclaim

\pf Recall that, for ${\bold l}\ne\emptyset$, the lengths of the base and top side of
$R_{{\bold l},{\bold q}}(x)$ are $x$ and $x+(l_m-m)-(q_n-n)+1$, respectively. 
The assumption in (a) implies therefore that the base is at most as long as the 
top side. Thus,
by setting $x=0$ we obtain a valid region, which has all the lozenges along its 
southwestern side forced (see Figure 4.7(a) for an illustration). The leftover region is 
readily identified to be $R_{{\bold l}^{(m)},{\bold q}}(l_m-l_{m-1}-1)$.

\topinsert
\twoline{\mypic{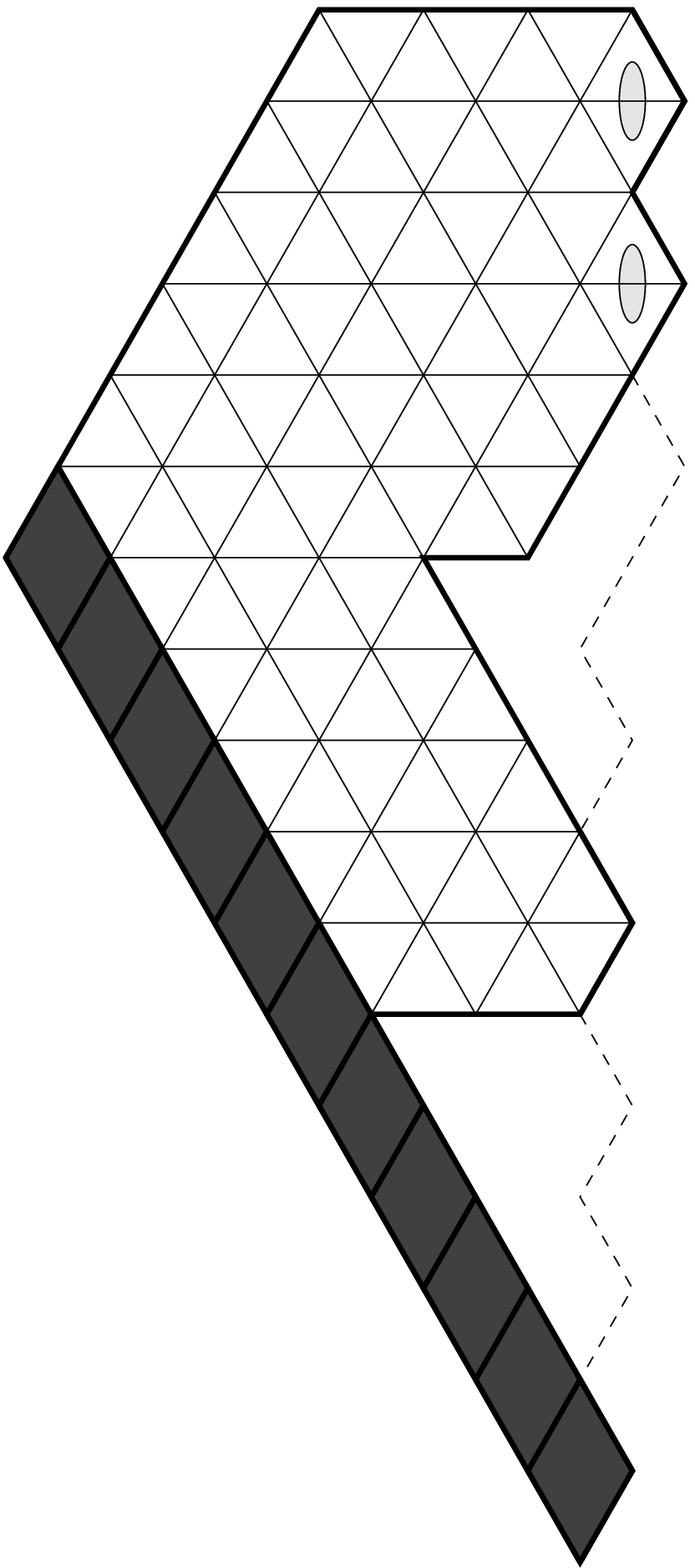}}{\mypic{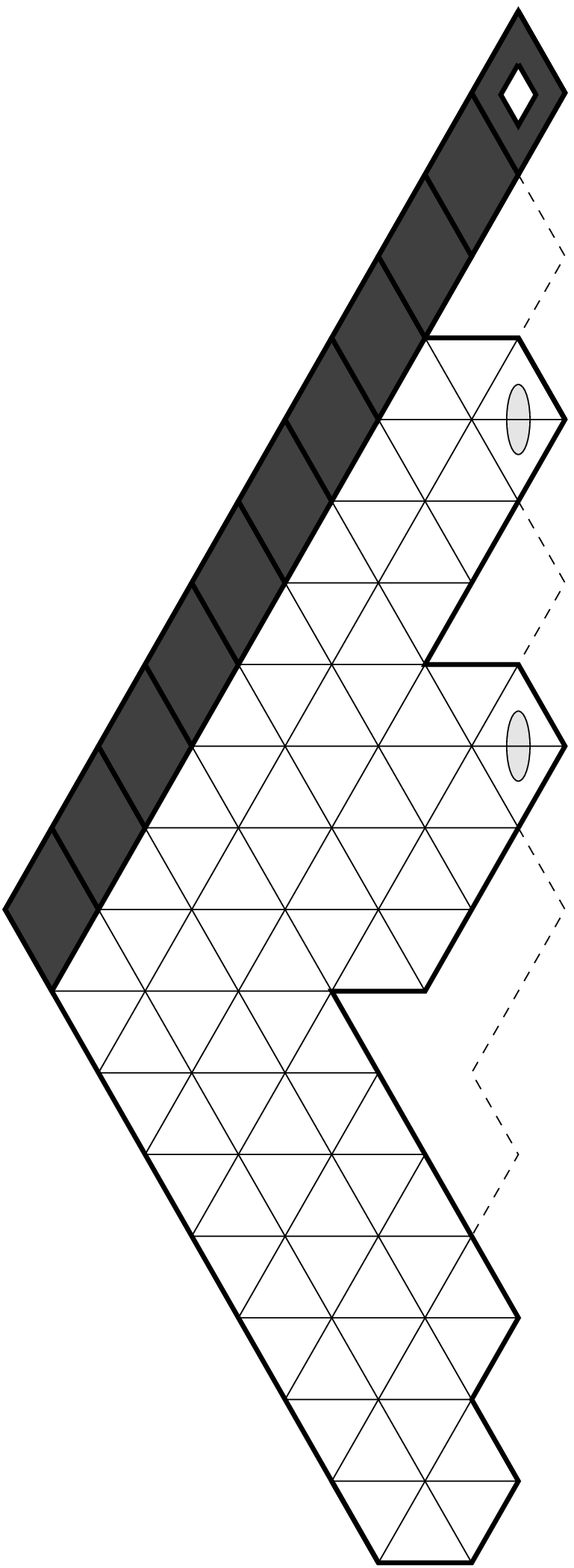}}
\twoline{Figure~4.7{\rm (a).}}{Figure~4.7{\rm (b).}}
\endinsert

Consider now part (b). Our assumption implies that the top side is shorter than or equal
to the base. Indeed, for ${\bold l}\ne\emptyset$ this follows from the above paragraph.
For ${\bold l}=\emptyset$, the top side has length $x-(q_n-n)+1$, and the base
has length $x-q_1+2$. Therefore,
it follows that the top side is at most as long as the base also in this case.

By setting $x=q_n-l_m-n+m-1$ the top side shrinks to a point and we obtain a valid
region whose lozenges along the northwestern side are forced (see Figure 4.7(b)). 
The leftover region is 
easily seen to be $R_{{\bold l},{\bold q}^{(n)}}(q_n-l_m-n+m-1)$. Since the top 
forced tile has weight 1/2 and the others have weight 1, we obtain (4.15). $\square$

\proclaim{Lemma 4.7} $(a)$ Suppose $l_m-m\geq q_n-n$. Then for ${\bold l}\ne\emptyset$
we have 

$$\M(\bar{R}_{{\bold l},{\bold q}}(0))=
\M(\bar{R}_{{\bold l}^{(m)},{\bold q}}(l_m-l_{m-1}-1)).\tag4.16$$

$(b)$ Suppose $l_m-m\leq q_n-n$. Then for ${\bold q}\ne\emptyset$ we have

$$\M(\bar{R}_{{\bold l},{\bold q}}(q_n-l_m-n+m))=
\frac12 \M(\bar{R}_{{\bold l},{\bold q}^{(n)}}(q_n-l_m-n+m)).\tag4.17$$
\endproclaim

\topinsert
\twoline{\mypic{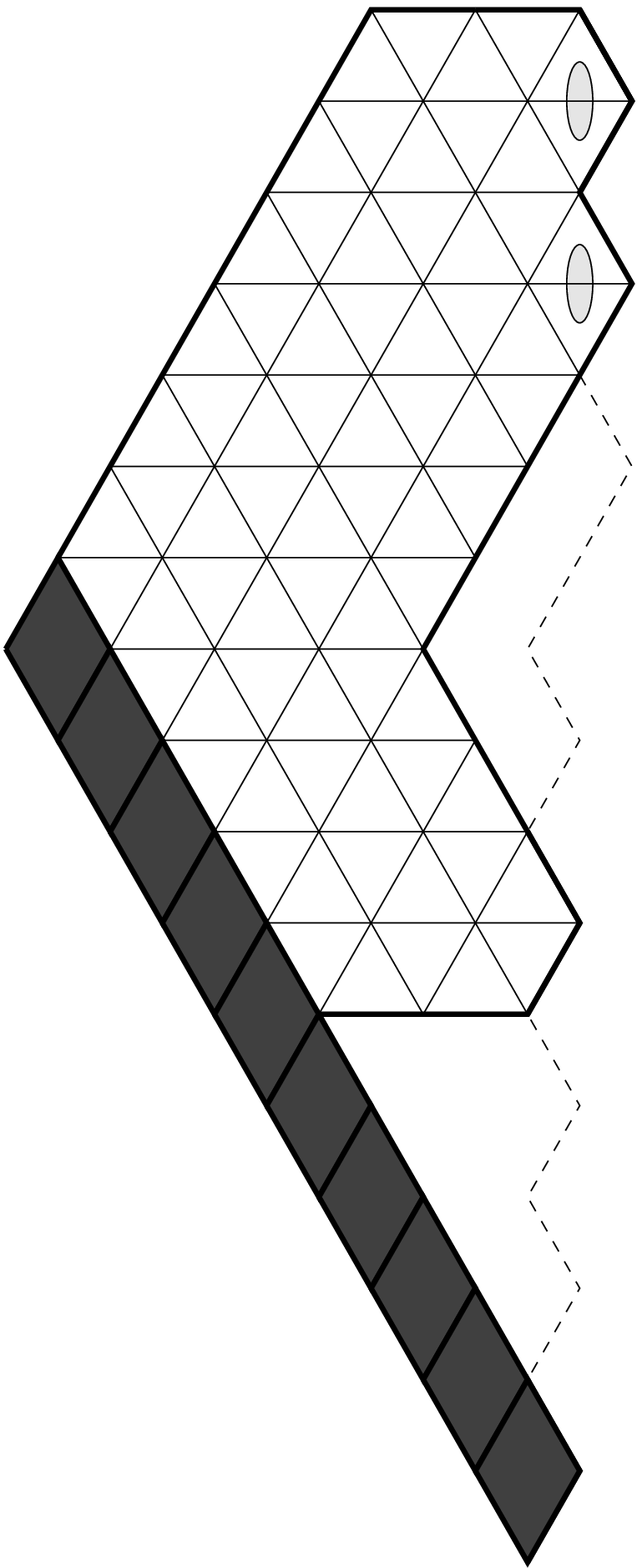}}{\mypic{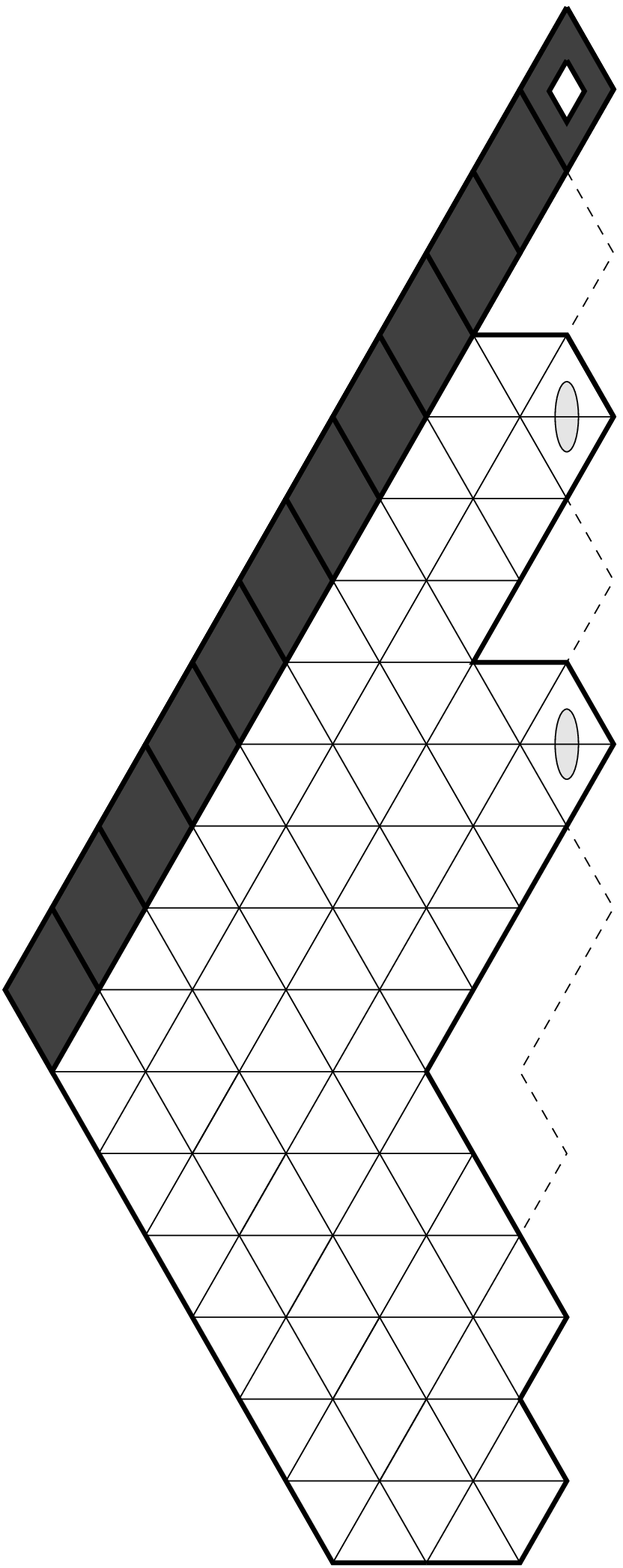}}
\twoline{Figure~4.8{\rm (a).}}{Figure~4.8{\rm (b).}}
\endinsert

\pf For ${\bold l}\ne\emptyset$, the lengths of the base and top side of
$\bar{R}_{{\bold l},{\bold q}}(x)$ are $x$ and $x+(l_m-m)-(q_n-n)$, respectively. 
Therefore, the assumption in (a) implies that the base is at most as long as the 
top side. Thus,
by setting $x=0$ we obtain a valid region, which has all the lozenges along its 
southwestern side forced (see Figure 4.8(a)). It is not hard to see that the leftover 
region is precisely $\bar{R}_{{\bold l}^{(m)},{\bold q}}(l_m-l_{m-1}-1)$.

For part (b), note that our assumption implies that the top side is at most as long as 
the base. Indeed, we have already seen this for ${\bold l}\ne\emptyset$. 
For ${\bold l}=\emptyset$, the top side and the base have lengths $x-(q_n-n)$ and 
$x-q_1+1$, respectively. Therefore, the top side is as least as long as the base in 
this case as well.

By setting $x=q_n-l_m-n+m$, the top side shrinks to a point and we obtain a valid
region whose lozenges along the northwestern side are forced (see Figure 4.8(b)). 
The leftover region is 
readily seen to be $\bar{R}_{{\bold l},{\bold q}^{(n)}}(q_n-l_m-n+m)$. Since the top 
forced tile has weight 1/2 and all others have weight 1, this proves (4.17). $\square$

\mysec{5. Proof of Proposition 2.1}

\medskip
In this section we show that the polynomials $P_{{\bold l},{\bold q}}(x)$ 
and $\bar{P}_{{\bold l},{\bold q}}(x)$ satisfy the same recurrences as
$\M(R_{{\bold l},{\bold q}}(x))$ and $\M(\bar{R}_{{\bold l},{\bold q}}(x))$. Then 
we deduce Proposition 2.1 by induction on $l_m+n+x$.

Rewrite formulas (1.5) and (1.6) as follows:

$$\align
P_{{\bold l},{\bold q}}(x)=c_{{\bold l},{\bold q}}B_{m,n}(x+l_m-m)
&\prod_{i=1}^{m}\prod_{j=i}^{l_i-1}(x+l_m-j)(x+l_m-m+n+j+2)\\
&\prod_{i=1}^{n}\prod_{j=i}^{q_i-1}(x+l_m-m+n-j+1)(x+l_m+j+1)\tag5.1\\
\bar{P}_{{\bold l},{\bold q}}(x)=\bar{c}_{{\bold l},{\bold q}}\bar{B}_{m,n}(x+l_m-m)
&\prod_{i=1}^{m}\prod_{j=i}^{l_i-1}(x+l_m-j)(x+l_m-m+n+j+1)\\
&\prod_{i=1}^{n}\prod_{j=i}^{q_i-1}(x+l_m-m+n-j)(x+l_m+j+1).\tag5.2
\endalign$$

\proclaim{Lemma 5.1} For $m\leq n$, we have

$$P_{{\bold l},{\bold q}}(x)=\sum_{k=1}^{n}(-1)^{n-k}C_k
P_{{\bold l},{\bold q}^{(k)}}(x)+
(-1)^n\delta_{mn}\bar{P}_{{\bold l},{\bold q}}(x),\tag5.3$$
where $\delta_{mn}$ is the Kronecker symbol and $C_k$ is given by $(4.3)$.
\endproclaim

\pf The outline of our proof is the following. 
Dividing both sides of (5.3) by $C_nP_{{\bold l},{\bold q}}(x)$, the terms on the 
right hand side turn out to become ratios of linear polynomials in $x$, all having the
same numerator. Dividing through further by this common numerator, the identity to be 
checked becomes the partial fraction decomposition of the rational function 
on the left hand side. Since the denominator of the latter is a product of distinct 
linear polynomials in $x$, this identity can be readily verified.

The expression on the right hand side of (4.3) giving $C_k$ can be rewritten as

$$C_k=\frac{(2x+2l_m-m+n+2)(x+l_m-q_k-m+n+2)_{2q_k+m-n-1}}{2\,(2q_k+m-n)\,!}.\tag5.4$$

{\it The left hand side}. We have

$$\align
\frac{P_{{\bold l},{\bold q}}}{C_nP_{{\bold l},{\bold q}^{(n)}}}&=
\frac{c_{{\bold l},{\bold q}}}{c_{{\bold l},{\bold q}^{(n)}}}
\frac{2\,(2q_n+m-n)\,!}{(2x+2l_m-m+n+2)(x+l_m-q_n-m+n+2)_{2q_n+m-n-1}}\\
&\ \ \ \ \ \cdot\frac{\prod_{i=1}^{m}\prod_{j=i}^{l_i-1}(x+l_m-j)(x+l_m-m+n+j+2)}
{\prod_{i=1}^{m}\prod_{j=i}^{l_i-1}(x+l_m-j)(x+l_m-m+n+j+1)}\\
&\ \ \ \ \ \cdot\frac{\prod_{i=1}^{n}\prod_{j=i}^{q_i-1}(x+l_m+j+1)(x+l_m-m+n-j+1)}
{\prod_{i=1}^{n-1}\prod_{j=i}^{q_i-1}(x+l_m+j+1)(x+l_m-m+n-j)}\\
&\ \ \ \ \ \cdot\frac{B_{m,n}(x+l_m-m)}{B_{m,n-1}(x+l_m-m)}
\endalign$$
$$\align
\phantom{\frac{P_{{\bold l},{\bold q}}}{C_nP_{{\bold l},{\bold q}^{(n)}}}}&=
\frac{2^{{n-m\choose 2}-m+1}}{2^{{n-m-1\choose 2}-m}}\,(2q_n+m-n)\,!\\
&\ \ \ \ \ \cdot\frac{{\displaystyle \prod_{i=1}^{m}\frac{1}{(2l_i)\,!}\prod_{i=1}^{n}
\frac{1}{(2q_i-1)\,!}   
\frac{\prod_{1\leq i<j\leq m}(l_j-l_i)\prod_{1\leq i<j\leq n}(q_j-q_i)} 
{\prod_{i=1}^{m}\prod_{j=1}^{n}(l_i+q_j)}  }  }
{{\displaystyle  \prod_{i=1}^{m}\frac{1}{(2l_i)\,!}\prod_{i=1}^{n-1}\frac{1}{(2q_i-1)\,!}   
\frac{\prod_{1\leq i<j\leq m}(l_j-l_i)\prod_{1\leq i<j\leq n-1}(q_j-q_i)}
{\prod_{i=1}^{m}\prod_{j=1}^{n-1}(l_i+q_j)}  }  }\\
&\ \ \ \ \ \cdot\prod_{i=1}^{m}\frac{x+l_m-m+n+l_i+1}{x+l_m-m+n+i+1}
\prod_{i=1}^{n-1}\frac{x+l_m-m+n-i+1}{x+l_m-m+n-q_i+1}\\
&\ \ \ \ \ \cdot\prod_{j=n}^{q_n-1}(x+l_m+j+1)(x+l_m-m+n-j+1)\\
&\ \ \ \ \ \cdot\frac{1}{2^{m+n}}(2x+2l_m-m+n+2)_{m+n}
\prod_{i=0}^{m-1}\frac{x+l_m+n-i+1}{x+l_m+n-i+1/2},\tag5.5
\endalign$$
where we used the (easily checked) fact that

$$\frac{B_{m,n}(x)}{B_{m,n-1}(x)}=\frac{1}{2^{m+n}}(2x+m+n+2)_{m+n}
\prod_{i=0}^{m-1}\frac{x+m+n-i+1}{x+m+n-i+1/2}.\tag5.6$$

After cancelling out common factors in the numerator and denominator of (5.5), and
then rearranging the remaining factors, we obtain

$$\align
\frac{P_{{\bold l},{\bold q}}}{C_nP_{{\bold l},{\bold q}^{(n)}}}=
&\frac{(2q_n+m-n)\,!}{(2q_n-1)\,!}
\frac{\prod_{i=1}^{n-1}(q_n-q_i)}{\prod_{i=1}^m(l_i+q_n)}\\
\cdot&\frac{(2x+2l_m-m+n+2)_{n-m}}{(2x+2l_m-m+n+2)}
\frac{\prod_{i=1}^{m}(x+l_m+l_i-m+n+1)}{\prod_{i=1}^{n-1}(x+l_m-q_i-m+n+1)}.\tag5.7
\endalign$$

{\it The terms in the sum on the right hand side.} Using (5.1) and (5.2) we obtain

$$\align
\frac{C_kP_{{\bold l},{\bold q}^{(k)}}}{C_nP_{{\bold l},{\bold q}^{(n)}}}&=
\frac{c_{{\bold l},{\bold q}^{(k)}}}{c_{{\bold l},{\bold q}^{(n)}}}
\frac{(2q_n+m-n)\,!}{(2q_k+m-n)\,!}
\frac{(x+l_m-q_k-m+n+2)_{2q_k+m-n-1}}{(x+l_m-q_n-m+n+2)_{2q_n+m-n-1}}\\
&\ \ \ \ \ \cdot\frac{\prod_{i=k}^{n-1}\prod_{j=i}^{q_{i+1}-1}(x+l_m+j+1)(x+l_m-m+n-j)}
{\prod_{i=k}^{n-1}\prod_{j=i}^{q_{i}-1}(x+l_m+j+1)(x+l_m-m+n-j)}\\
&\ \\
&=\frac{(2q_n+m-n)\,!}{(2q_k+m-n)\,!}\\
&\ \ \ \ \ \cdot\frac{{\displaystyle \prod_{i=1}^{m}\frac{1}{(2l_i)\,!}\prod_{i=1,\,i\ne k}^{n}
\frac{1}{(2q_i-1)\,!}
\frac{\prod_{1\leq i<j\leq m}(l_j-l_i)\prod_{1\leq i<j\leq n,\,i,j\ne k}(q_j-q_i)}
{\prod_{i=1}^{m}\prod_{j=1,\,j\ne k}^{n}(l_i+q_j)}  }  }
{{\displaystyle \prod_{i=1}^{m}\frac{1}{(2l_i)\,!}\prod_{i=1}^{n-1}\frac{1}{(2q_i-1)\,!}   
\frac{\prod_{1\leq i<j\leq m}(l_j-l_i)\prod_{1\leq i<j\leq n-1}(q_j-q_i)}
{\prod_{i=1}^{m}\prod_{j=1}^{n-1}(l_i+q_j)}  }  }\\
&\ \ \ \ \ \cdot\frac{(x+l_m-q_k-m+n+2)_{2q_k+m-n-1}}{(x+l_m-q_n-m+n+2)_{2q_n+m-n-1}}\\
&\ \ \ \ \ \cdot\prod_{i=k}^{n-1}\prod_{j=q_i}^{q_{i+1}-1}(x+l_m+j+1)(x+l_m-m+n-j).\tag5.8
\endalign$$

By cancelling out common factors and regrouping one obtains

$$\align
\frac{C_kP_{{\bold l},{\bold q}^{(k)}}}{C_nP_{{\bold l},{\bold q}^{(n)}}}=
&\frac{(2q_k-1)\,!}{(2q_n-1)\,!}\frac{(2q_n+m-n)\,!}{(2q_k+m-n)\,!}
\prod_{i=1}^m\frac{l_i+q_k}{l_i+q_n}\\
\cdot&\frac{\prod_{i=1,\,i\ne k}^{n-1}(q_n-q_i)}
{\prod_{i=1}^{k-1}(q_k-q_i)\prod_{i=k+1}^{n-1}(q_i-q_k)}
\frac{x+l_m-q_n-m+n+1}{x+l_m-q_k-m+n+1}.\tag5.9
\endalign$$

{\it The last term on the right hand side.} Since for $m\ne n$ this term is 0, assume
$m=n$. We have

$$\align
\frac{\bar{P}_{{\bold l},{\bold q}}}{C_nP_{{\bold l},{\bold q}^{(n)}}}&=
\frac{\bar{c}_{{\bold l},{\bold q}}}{c_{{\bold l},{\bold q}^{(n)}}}
\frac{2\,(2q_n)\,!}{(2x+2l_n+2)(x+l_n-q_n+2)_{2q_n-1}}\\
&\ \ \ \ \ \cdot\frac{\prod_{i=1}^{n}\prod_{j=i}^{l_i-1}(x+l_n-j)(x+l_n+j+1)}
{\prod_{i=1}^{n}\prod_{j=i}^{l_i-1}(x+l_n-j)(x+l_n+j+1)}\\
&\ \ \ \ \ \cdot\frac{\prod_{i=1}^{n}\prod_{j=i}^{q_i-1}(x+l_n-j)(x+l_n+j+1)}
{\prod_{i=1}^{n-1}\prod_{j=i}^{q_i-1}(x+l_n-j)(x+l_n+j+1)}\\
&\ \ \ \ \ \cdot\frac{\bar{B}_{n,n}(x+l_n-n)}{B_{n,n-1}(x+l_n-n)}\\
&\ \\
&=\frac{2^{{0\choose2}-n+1}}{2^{{-1\choose 2}-n}}\,(2q_n)\,!\\
&\ \ \ \ \ \cdot\frac{{\displaystyle \prod_{i=1}^{n}\frac{1}{(2l_i-1)\,!}
\prod_{i=1}^{n}\frac{1}{(2q_i)\,!}   
\frac{\prod_{1\leq i<j\leq n}(l_j-l_i)\prod_{1\leq i<j\leq n}(q_j-q_i)}
{\prod_{i=1}^{n}\prod_{j=1}^{n}(l_i+q_j)}  }  }
{{\displaystyle \prod_{i=1}^{n}\frac{1}{(2l_i)\,!}\prod_{i=1}^{n-1}\frac{1}{(2q_i-1)\,!}   
\frac{\prod_{1\leq i<j\leq n}(l_j-l_i)\prod_{1\leq i<j\leq n-1}(q_j-q_i)}
{\prod_{i=1}^{n}\prod_{j=1}^{n-1}(l_i+q_j)}  }  }\\
&\ \ \ \ \ \cdot\frac{\prod_{j=n}^{q_n-1}(x+l_n-j)(x+l_n+j+1)}
{(2x+2l_n+2)(x+l_n-q_n+2)_{2q_n-1}}\\
&\ \ \ \ \ \cdot(x+l_n-n+1)_{2n},\tag5.10
\endalign$$
where we have used the equality

$$\frac{\bar{B}_{n,n}(x)}{B_{n,n-1}(x)}=(x+1)_{2n},\tag5.11$$
which can be readily checked. Simplifying and rearranging in (5.10) leads to

$$
\frac{\bar{P}_{{\bold l},{\bold q}}}{C_nP_{{\bold l},{\bold q}^{(n)}}}=
\frac{l_1\cdots l_n}{q_1\cdots q_{n-1}}
\frac{\prod_{i=1}^{n-1}(q_n-q_i)}{\prod_{i=1}^{n}(l_i+q_n)}
\frac{x+l_n-q_n+1}{x+l_n+1}.\tag5.12
$$

By (5.7), (5.9) and (5.12), it follows that the statement of the Lemma is equivalent to 
the following rational function identity:

$$\align
\frac{(2q_n+m-n)\,!}{(2q_n-1)\,!}
&\frac{\prod_{i=1}^{n-1}(q_n-q_i)}{\prod_{i=1}^m(l_i+q_n)}\\
\cdot&\frac{(2x+2l_m-m+n+2)_{n-m}}{(2x+2l_m-m+n+2)}
\frac{\prod_{i=1}^{m}(x+l_m+l_i-m+n+1)}{\prod_{i=1}^{n}(x+l_m-q_i-m+n+1)}\\
&=\sum_{k=1}^{n}(-1)^{n-k}\frac{\alpha_k}{x+l_m-q_k-m+n+1}+
\delta_{mn}\frac{\beta}{x+l_n-q_n+1},\tag5.13
\endalign$$
where $\alpha_k$ and $\beta$ are given by

$$\align
\alpha_k&=
\frac{(2q_k-1)\,!}{(2q_n-1)\,!}\frac{(2q_n+m-n)\,!}{(2q_k+m-n)\,!}
\frac{\prod_{i=1}^m(l_i+q_k)}{\prod_{i=1}^m(l_i+q_n)}
\frac{\prod_{i=1,\,i\ne k}^{n-1}(q_n-q_i)}
{\prod_{i=1}^{k-1}(q_k-q_i)\prod_{i=k+1}^{n-1}(q_i-q_k)}\\
\beta&=\frac{l_1\cdots l_n}{q_1\cdots q_{n-1}}
\frac{\prod_{i=1}^{n-1}(q_n-q_i)}{\prod_{i=1}^{n}(l_i+q_n)}.
\endalign$$

However, (5.13) can be regarded as the partial fraction decomposition of the rational 
function on its left hand side. Since the denominator is factored in distinct linear
factors, the values of the constants $\alpha_k$ and $\beta$ that make (5.13) true can be
immediately obtained by suitable specializations of $x$. It is easy to see that these
values are precisely the ones in the above two equalities. $\square$

\proclaim{Lemma 5.2} Let $D_k$, $\bar{C}_k$ and $\bar{D}_k$ be given by $(4.7)$, 
$(4.11)$ and $(4.13)$.

$(a)$ For $m>n$, we have

$$P_{{\bold l},{\bold q}}(x)=\sum_{k=1}^{m-1}(-1)^{m-k}D_k
P_{{\bold l}^{(k)},{\bold q}}(x-1)+
D_mP_{{\bold l}^{(m)},{\bold q}}(x+l_m-l_{m-1}-1).\tag5.14$$

$(b)$ For $m<n$, we have

$$\bar{P}_{{\bold l},{\bold q}}(x)=\sum_{k=1}^{n}(-1)^{n-k}\bar{C}_k
\bar{P}_{{\bold l},{\bold q}^{(k)}}(x).\tag5.15$$

$(c)$ For $m\geq n$, we have

$$\align
\bar{P}_{{\bold l},{\bold q}}(x)=\sum_{k=1}^{m-1}(-1)^{m-k}\bar{D}_k
\bar{P}_{{\bold l}^{(k)},{\bold q}}(x-1&)+
\bar{D}_m\bar{P}_{{\bold l}^{(m)},{\bold q}}(x+l_m-l_{m-1}-1)\\
+(-1)^m\delta_{mn}P_{{\bold l},{\bold q}}(x-1&).\tag5.16
\endalign$$
\endproclaim

\pf The identities (5.14)--(5.16) can be proved using the same approach as in
the proof of Lemma 5.1. For the first two, we divide through by the last term 
on the right hand side, and for the third by the second to last. In each of the
three cases, the terms on the right hand sides become ratios of linear polynomials, all
with the same numerator. Dividing through by this common numerator, the identity to be
proved becomes the partial fraction decomposition of the rational function on the left hand
side. In all three cases, the denominator of this rational function is a product of
distinct linear factors, and thus, as in the proof of Lemma 5.1, the identity is easily 
checked by suitable specializations of $x$.

This can be carried out in a straightforward way. The shifts in the arguments in 
(5.14)--(5.16) are so that the division mentioned above results in great simplification
of the involved expressions. The only difference that needs to be pointed out is that
we need the following relations, analogous to (5.6) and (5.11), involving the monic 
polynomials $B_{m,n}$ and $\bar{B}_{m,n}$:

$$\align
\frac{B_{m,n}(x)}{B_{m-1,n}(x)}=&\frac{1}{2^{m+n-1}}(x+m+n)(x+m+n+1)(2x+m+n+2)_{m+n-1}\\
&\ \ \ \cdot\prod_{i=0}^{n-1}\frac{x+m+i}{x+m+i+1/2}\\
\frac{\bar{B}_{m,n}(x)}{\bar{B}_{m,n-1}(x)}=&\frac{1}{2^{m+n}}(x+m+n)(2x+m+n+1)_{m+n}\\
&\ \ \ \cdot\prod_{i=0}^{m-1}\frac{x+m+n-i-1}{x+m+n-i-1/2}\\
\frac{\bar{B}_{m,n}(x)}{\bar{B}_{m-1,n}(x)}=&\frac{1}{2^{m+n}}(2x+m+n+1)_{m+n}
\prod_{i=0}^{n-1}\frac{x+m+i+1}{x+m+i+1/2}\\
&\ \ \ \ \ \ \ \ \ \ \ \ \ \ \ \ 
\frac{B_{n,n}(x-1)}{\bar{B}_{n-1,n}(x)}=\frac{(x)_{2n+1}}{x+n}.
\endalign$$
These are all easily verified directly using the defining formulas for $B_{m,n}$ and 
$\bar{B}_{m,n}$. $\square$

\proclaim{Lemma 5.3} We have

$$\align
P_{{\bold l},{\bold q}}(q_n-l_m-n+m-1)&=
\frac{1}{2}P_{{\bold l},{\bold q}^{(n)}}(q_n-l_m-n+m-1)\\
\bar{P}_{{\bold l},{\bold q}}(q_n-l_m-n+m)&=
\frac{1}{2}\bar{P}_{{\bold l},{\bold q}^{(n)}}(q_n-l_m-n+m),
\endalign$$

and for ${\bold l}\ne\emptyset$ we have

$$\align
P_{{\bold l},{\bold q}}(0)&=P_{{\bold l}^{(m)},{\bold q}}(l_{m}-l_{m-1}-1)\\
\bar{P}_{{\bold l},{\bold q}}(0)&=
\bar{P}_{{\bold l}^{(m)},{\bold q}}(l_{m}-l_{m-1}-1).
\endalign$$
\endproclaim

\pf For $m\leq n$, from (5.7) and (5.4) one obtains a simple expression for the ratio 
$P_{{\bold l},{\bold q}}(x)/P_{{\bold l},{\bold q}^{(n)}}(x)$. It is readily seen that
for $x=q_n-l_m-n+m-1$ this expression becomes 1/2, thus proving the first equality in
the statement of the Lemma in this case. For $m>n$, one obtains from (5.5)
an expression for $P_{{\bold l},{\bold q}}(x)/(C_nP_{{\bold l},{\bold q}^{(n)}}(x))$
similar to (5.7). Again, it follows easily that the ratio
$P_{{\bold l},{\bold q}}(x)/P_{{\bold l},{\bold q}^{(n)}}(x)$ specializes to 1/2 for
$x=q_n-l_m-n+m-1$.

The remaining three equalities are proved similarly, just as the proof of Lemma 5.2 is
similar to that of Lemma 5.1. $\square$

\medskip
{\it Proof of Proposition 2.1.} We proceed by induction on $l_m+n+x$. For 
$R_{{\bold l},{\bold q}}(x)$, the minimum possible value of $l_m+n+x$ is $-1$. 
In this case
${\bold l}$ and ${\bold q}$ are necessarily empty, and (2.1) is true by our
definitions. For $\bar{R}_{{\bold l},{\bold q}}(x)$, the minimum value of $l_m+n+x$ is 0.
Again, this implies ${\bold l}={\bold q}=\emptyset$, and (2.2) follows from our definitions.

Let $N>0$ and suppose that (2.1) and (2.2) are true for all ${\bold l}$, 
${\bold q}$ and $x$ with $l_m+n+x<N$. Consider ${\bold l}$, ${\bold q}$ and $x$ so
that $l_m+n+x=N$. We claim that (2.2) holds also for these choices of 
${\bold l}$, ${\bold q}$ and $x$.

Indeed, assume first that $x$ is not equal to its minimum possible value.
Suppose $m\geq n$. By (4.12), $\M(\bar{R}_{{\bold l},{\bold q}}(x))$
can be expressed in terms of tiling generating functions of regions from the
$\bar{R}$- and (in case $m=n$) $R$-family. By the induction hypothesis, (2.1) and (2.2)
hold for all the regions on the right hand side of (4.12). Then Lemma 5.2(c) implies 
that (2.2) holds in this case. For $m<n$, (2.2) follows similarly, from (4.10), the 
induction hypothesis and Lemma 5.2(b).

Assume now that $x$ is equal to its minimum possible value. Then $x$ equals the value of 
the argument on the left hand side of either (4.16) or (4.17), and we can express 
$\M(\bar{R}_{{\bold l},{\bold q}})(x)$ in terms of the tiling generating function of
a region from the $\bar{R}$-family for which the induction hypothesis applies. 
Therefore, Lemma 5.3 implies that (2.2) holds for 
$\bar{R}_{{\bold l},{\bold q}}(x)$ as well.

We now show that (2.1) holds for $l_m+n+x=N$. Assume that $x$ is not equal to its 
minimum possible value. Suppose $m\leq n$, and apply Lemma 4.2. By the induction
hypothesis, (2.1) holds for the regions appearing in the terms of the sum on the
right hand side of (4.2). Furthermore, we have seen above that (2.2) holds for 
$\bar{R}_{{\bold l},{\bold q}}(x)$. Therefore, by Lemma 5.1 we obtain that (2.1)
is true for $R_{{\bold l},{\bold q}}(x)$, if $m\leq n$. The case $m>n$ follows
analogously, using Lemma 4.3, the induction hypothesis and Lemma 5.2(a).

Finally, if $x$ takes on its minimum possible value, (2.1) follows by (4.14) and (4.15),
the induction hypothesis and Lemma 5.3. This completes the proof of Proposition 2.1 by 
induction. $\square$ 

\mysec{6. The guessing of $\M(R_{{\bold l},{\bold q}}(x))$ 
and $\M(\bar{R}_{{\bold l},{\bold q}}(x))$} 

\medskip
In this section we describe how we arrived at conjecturing that the tiling generating
functions for the regions $R_{{\bold l},{\bold q}}(x)$ and 
$\bar{R}_{{\bold l},{\bold q}}(x))$ are given by the polynomials 
$P_{{\bold l},{\bold q}}(x)$ and $\bar{P}_{{\bold l},{\bold q}}(x)$.

We discuss here the case of $\bar{R}_{{\bold l},{\bold q}}(x)$. The region
$R_{{\bold l},{\bold q}}(x)$ can be treated similarly. 

It is easy to write down explicitly the entries of the matrix in (4.1). This gives
an easy way to compute the polynomials 
$\M(\bar{R}_{{\bold l},{\bold q}}(x))$ by computer, for specific choices of 
${\bold l}$ and ${\bold q}$. At once, one is striked by the fact that all these
polynomials seem to factor in the form

$$\M(\bar{R}_{{\bold l},{\bold q}}(x))=\bar{c}_{{\bold l},{\bold q}}
F_{{\bold l},{\bold q}}(x),\tag6.1$$
where $\bar{c}_{{\bold l},{\bold q}}$ is a constant and $F_{{\bold l},{\bold q}}(x)$
is the product of linear factors of the form $(x+t)$, where $t$ is either an integer or
a half-integer. 

Moreover, it is easy to conjecture the behaviour of $F_{{\bold l},{\bold q}}(x)$ under
incrementing a single element of ${\bold l}$ or ${\bold q}$. Specifically, after some
experimenting one soon comes up with the guess that

$$\align
&\frac{F_{{\bold l}^{|k\rangle},{\bold q}}(x)}{F_{{\bold l},{\bold q}}(x)}=
(x-l_k+l_m)(x+l_k+l_m-m+n+1),\ \ \ \,\text{for $1\leq k<m$}\tag6.2\\
&\frac{F_{{\bold l}^{|m\rangle},{\bold q}}(x-1)}{F_{{\bold l},{\bold q}}(x)}=
x(x+2l_m-m+n+1)\tag6.3\\
&\frac{F_{{\bold l},{\bold q}^{|k\rangle}}(x)}{F_{{\bold l},{\bold q}}(x)}=
(x+q_k+l_m+1)(x-q_k+l_m-m+n),\ \ \ \text{for $1\leq k\leq n$,}\tag6.4
\endalign$$
where ${\bold l}^{|k\rangle}$ is the list obtained from ${\bold l}$ by increasing its $k$-th
element by 1 (so for $k<m$, ${\bold l}^{|k\rangle}$ is defined only if $l_{k+1}-l_k\geq2$). 
These formulas
provide us with an explicit guess for $F_{{\bold l},{\bold q}}$, in terms of
$F_{[m],[n]}(x)$, where $[m]$ is the list consisting of $1,\dotsc,m$. Furthermore,
some experimentation easily leads one to conjecture that 
$F_{[m],[n]}(x)=\bar{B}_{m,n}(x)$ (where $\bar{B}_{m,n}(x)$ is given by (1.2)).

It is considerably more difficult to guess a formula for the constants
$\bar{c}_{{\bold l},{\bold q}}$ just based on working out examples (this is not
surprising, as factorization of integers gives significantly less information than
factorization of polynomials). However, it turns out that we do not need to guess
the value of these constants: substituting the conjectured expression for  
$F_{{\bold l},{\bold q}}$ in (6.1) and using then Lemma 4.7, we obtain
(conjectured) recurrences satisfied by the $\bar{c}_{{\bold l},{\bold q}}$'s. 

More precisely, suppose $l_m-m\geq q_n-n$. We may then
use (4.16) to relate $\bar{c}_{{\bold l},{\bold q}}$ to
$\bar{c}_{{\bold l}^{(m)},{\bold q}}$. Carrying this out one arrives at

$$\frac{\bar{c}_{{\bold l},{\bold q}}}{\bar{c}_{{\bold l}^{(m)},{\bold q}}}
=2^{m-n-1}\frac{1}{(2l_m-1)\,!}
\frac{\prod_{i=1}^{m-1}(l_m-l_i)}{\prod_{i=1}^{n}(l_m+q_i)}.\tag6.5$$

On the other hand, if $l_m-m\leq q_n-n$, we deduce similarly from (4.17) that

$$\frac{\bar{c}_{{\bold l},{\bold q}}}{\bar{c}_{{\bold l},{\bold q}^{(n)}}}
=2^{n-m-1}\frac{1}{(2q_n)\,!}
\frac{\prod_{i=1}^{n-1}(q_n-q_i)}{\prod_{i=1}^{m}(q_n+l_i)}.\tag6.6$$

We can then repeat this procedure for ${\bold l}$ and ${\bold q}$ replaced by
${\bold l}^{(m)}$ and ${\bold q}$ or ${\bold l}$ and ${\bold q}^{(n)}$, depending
on which of (4.16) and (4.17) applies. 
At each iteration we have to look at the
initial segments of the lists ${\bold l}$ and ${\bold q}$ that are left over, and
compare the difference between their largest entry and their number of entries, before
deciding which of (6.5) or (6.6) to apply. 
However, it is not hard to see that no matter in which order we apply these
recurrences, they lead to an expression of the form 

$$\bar{c}_{{\bold l},{\bold q}}=2^{e({\bold l},{\bold q})}
\prod_{i=1}^{m}\frac{1}{(2l_i-1)\,!}\prod_{i=1}^{n}\frac{1}{(2q_i)\,!}
\frac{\prod_{1\leq i<j\leq m}(l_j-l_i)\prod_{1\leq i<j\leq n}(q_j-q_i)}
{\prod_{i=1}^{m}\prod_{j=1}^{n}(l_i+q_j)},$$
where $e({\bold l},{\bold q})$ is some integer depending only on ${\bold l}$ and ${\bold q}$.
Then after some experimenting one is readily lead to guess that 
$e({\bold l},{\bold q})={n-m\choose2}-m$, and we obtain (1.4). It is remarkable how the 
``weak'' conjectures (6.1) and (6.2)--(6.4) lead us, by the ``forcing''
argument in Lemma 4.7, to the precise (conjectured) formula for 
$\bar{c}_{{\bold l},{\bold q}}$.

A similar analysis can be done for $\M(R_{{\bold l},{\bold q}}(x))$, and one arrives
at formulas (1.5), (1.1) and (1.3).

\bigskip
{\bf Acknowledgments.} 
David Wilson's program Vaxmacs for counting perfect matchings was very 
useful in finding the families of regions $R_{{\bold l},{\bold q}}(x)$ and 
$\bar{R}_{{\bold l},{\bold q}}(x)$.

\mysec{References}
{\openup 1\jot \frenchspacing\raggedbottom
\roster

\myref{1}
  M. Ciucu, Enumeration of perfect matchings in graphs with reflective symmetry, 
{\it J. Combin. Theory Ser. A} {\bf 77} (1997), 67--97.
\myref{2}
  M. Ciucu, Enumeration of lozenge tilings of punctured hexagons, {\it J. Combin. Theory 
Ser. A,} {\bf 83} (1998), 268--272.
\myref{3}
  G. David and C. Tomei, The problem of the calissons, {\it Amer. Math. 
Monthly} {\bf 96} (1989), 429--431.
\myref{4}
  S. Elnitsky, Rhombic tilings of polygons and classes of reduced words in Coxeter 
groups. J. Combin. Theory Ser. A 77 (1997), 193--221. 
\myref{5}
  I. M. Gessel and X. Viennot, Binomial determinants, paths, and hook length formulae,
{\it Adv. in Math.} {\bf 58} (1985), 300--321.
\myref{6} 
  C. Krattenthaler, Generating functions for plane partitions of a given shape, 
{\it Ma\-nu\-scrip\-ta Math.} {\bf 69} (1990), 173--201. 
\myref{7}
  G. Kuperberg, Symmetries of plane partitions and the permanent-de\-ter\-mi\-nant
method, {\it J. Combin. Theory Ser. A} {\bf 68} (1994), 115--151.
\myref{8}
  P. A. MacMahon, Memoir on the theory of the partition of numbers---Part V. Partitions
in two-dimensional space, {\it Phil. Trans. R. S.}, 1911, A.
\myref{9}
  P. A. MacMahon, ``Combinatory Analysis,'' vols. 1--2, Cambridge, 1916, reprinted by 
Chelsea, New York, 1960.
\myref{10}
  J. Propp, Enumeration of matchings: problems and progress, in ``New Perspectives in
Geometric Combinatorics,'' MSRI Publications, vol. 38, 1999, 255--291. 
\myref{11}
  R. P. Stanley, Ordered structures and partitions, {\it Memoirs of the Amer. Math. Soc.,} 
no. 119 (1972). 
\myref{12}
  J. R. Stembridge, Nonintersecting paths, Pfaffians and plane partitions, 
{\it Adv. in Math.} {\bf 83} (1990), 96--131.

\endroster\par}

\enddocument